\documentclass[12pt ]{scrartcl}
\ifx\cs\documentclass \footheight 12pt \fi \footskip 30pt 
\usepackage[latin1]{inputenc}
\usepackage{amssymb}
\usepackage{color}
\usepackage{amsmath}\usepackage{amsmath}
\usepackage{amsfonts}
\usepackage{amssymb}
\usepackage{amsthm}
\usepackage{amsfonts}
\usepackage{exscale}
\usepackage{latexsym}
\usepackage{mathscinet}
\usepackage{overpic}
\newcommand {\Div}{\mathrm{div}}

\newtheorem{defi}{Definition}[section]
\newtheorem{remark}[defi]{Remark}
\newtheorem{lemma}[defi]{Lemma}
\newtheorem{theorem}[defi]{Theorem}

\newtheorem{proposition}[defi]{Proposition}

\def\proof{\noindent {\bf Proof: }}

\def\beq{\begin{equation}}
\def\eeq{\end{equation}}
\def\beqa{\begin{eqnarray*}}
\def\eeqa{\end{eqnarray*}}
\def\beqam{\begin{eqnarray}}
\def\eeqam{\end{eqnarray}}

\def\nz{{\relax\ifmmode I\!\!N\else$I\!\!N$\fi}}
\def\C{{\relax\ifmmode I\!\!\!\!C\else$I\!\!\!\!C$\fi}}

\newcommand{\eps}{\varepsilon}

\begin{document}

\begin{titlepage}
\title{Sharp Interface Limit for the Cahn-Larché System}
\author{Helmut Abels\footnote{Fakult\"at f\"ur Mathematik,  
Universit\"at Regensburg,
93040 Regensburg,
Germany, e-mail: {\sf helmut.abels@mathematik.uni-regensburg.de}}\ \ and Stefan Schaubeck\footnote{Fakult\"at f\"ur Mathematik,  
Universit\"at Regensburg,
93040 Regensburg,
Germany, e-mail: {\sf stefan.schaubeck@mathematik.uni-regensburg.de}}}
%\date{}
\end{titlepage}
\maketitle

\begin{abstract}
We prove rigorously  the convergence of the Cahn-Larch\'e system, which is a Cahn-Hilliard system coupled with the system of linearized elasticity, to a modified Hele-Shaw problem as long as a classical solution of the latter system exists. By matched asymptotic expansion we construct approximate solutions and estimate the difference between the true and the approximate solutions.
\end{abstract}
{\small\noindent
{\bf Mathematics Subject Classification (2000):}
%Primary: 35R35; Secondary  35Q30, 76D27,  76D45, 76T99.\vspace{0.1in}\\
%76T99, %% Two-Phase flows: Others
%76D27, %% Incompressible viscous fluids: Other free-boundary flows; Hele-Shaw flows
%76D03, %% Incompressible viscous fluids: Existence, uniqueness, and regularity theory 
%76D05, %% Incompressible viscous fluids: Navier-Stokes equations
%76D45, %% Incompressible viscous fluids: Capillarity (surface tension)
Primary: 74N15; %% Analysis of microstructure 
Secondary: 35B25, %% Singular perturbations 
35K55, %% Nonlinear parabolic equations 
35Q72, %% Other equations from mechanics 
49J10, %% Free problems in two or more independent variables 
74H10, %% Analytic approximation of solutions (perturbation methods, asymptotic methods, series, etc.) 
82C24.\\ %%  Interface problems; diffusion-limited aggregation 
{\bf Key words:} 
Sharp interface limit, diffuse interface model, Cahn-Hilliard equation, linearized elasticity, Hele-Shaw system, matched asymptotic expansion
%  Mullins-Sekerka equation,% convergence to equilibria. \vspace{0.1in}\\
}

\section{Introduction}

In the present contribution we consider the singular limit $\epsilon \searrow 0$ for the following
Cahn-Larché system
\begin{align}
\partial_t c^\epsilon & = \Delta \mu^\epsilon &&\mbox{in }\Omega \times (0,\infty), \label{system1}\\
\mu^\epsilon & = \epsilon^{-1} f(c^\epsilon) - \epsilon \Delta c^\epsilon + W_{,c}(c^\epsilon, \mathcal{E}(\textbf{u}^\epsilon))&& \mbox{in }\Omega \times (0,\infty), \label{system2} \\
\Div \, \mathcal{S}^\epsilon & = 0 && \mbox{in }\Omega \times (0,\infty), \label{system3} \\
\mathcal{S}^\epsilon & = W_{, \mathcal{E}}(c^\epsilon, \mathcal{E}(\textbf{u}^\epsilon)) && \mbox{in }\Omega \times (0,\infty) \label{system3b} .
\end{align}
 In the following $\Omega \subset \mathbb{R}^d$, $d >1$, is always a bounded domain with smooth boundary $\partial \Omega$. To close the system we require the following boundary and initial values
\begin{align}
\tfrac{\partial}{\partial n}c^\epsilon & = \tfrac{\partial}{\partial n}\mu^\epsilon=0 &&\mbox{on } \partial \Omega \times (0,\infty),\\
\textbf{u}^\epsilon & = 0 && \mbox{on } \partial \Omega \times (0,\infty),\\
\left.c^\epsilon \right|_{t=0}&= c^\epsilon_0 &&\mbox{in } \Omega  \label{system6},
\end{align}
where $n$ denotes the unit normal of $\partial \Omega$. 

The system describes the phase separation of a binary mixture (e.g.\, an alloy consisting of two components) taking the effect of an elastic misfit in the material into account. It was derived by Cahn and Larch\'e~\cite{CahnLarche}. For example, a different lattice structure of the mixture is a reason to consider elasticity in the Cahn-Hilliard model. 
Here  $c^\epsilon\colon \Omega \times (0,\infty) \to \mathbb{R}$ is the concentration difference of the components. The elastic deformations are assumed to be small and are described by the (infinitesimal) displacement $\textbf{u}^\epsilon : \Omega \times (0,\infty) \to \mathbb{R}^d$. Therefore the stress tensor $\mathcal{S}^\epsilon$ depends on the linearised strain tensor
\[
\mathcal{E}(\textbf{u}^\epsilon) = \frac{1}{2} \left( \nabla \textbf{u}^\epsilon + (\nabla \textbf{u}^\epsilon)^T \right) .
\] 
In the model it is assumed that the equation for the linearized elasticity (\ref{system3}) is static because the mechanical equilibrium is attained on a much faster time scale than the concentration changing by diffusion. 

In our case the elastic free energy density is described by 
\[
W(c, \mathcal{E}(\textbf{u})) = \frac{1}{2} \left( \mathcal{E}(\textbf{u}) - \mathcal{E}^\star c \right) : \mathcal{C} \left( \mathcal{E}(\textbf{u}) - \mathcal{E}^\star c \right),
\]
where $\mathcal{C} = (\mathcal{C}_{iji'j'})_{i,j,i',j'=1,\ldots,d}$ is the elasticity strain tensor and $\mathcal{E}^\star c$ is the stress free strain for concentration $c$ with constant matrix $\mathcal{E}^\star \in \mathbb{R}^{d \times d}$. We require that $\mathcal{C}$ is symmetric and positive definite. This form of the elastic free energy is based on the work of Eshelby \cite{Eshelby} and Khachaturyan \cite{Kha}. Then the total energy of the system is given by $E(c,\textbf{u}) = E_1(c) + E_2(c,\textbf{u})$, where
\begin{eqnarray}
E_1(c) & = & \frac{\epsilon}{2} \int_\Omega{\left| \nabla c(x) \right|^2 dx} + \frac{1}{\epsilon} \int_\Omega{F(c(x)) \, dx}  \label{GLenergy}
\end{eqnarray}
is the Ginzburg-Landau energy and 
\begin{eqnarray}
E_2(c,\textbf{u}) & = & \int_\Omega{W(c(x),\mathcal{E}(\textbf{u}(x))) \, dx} \label{elaenergy}
\end{eqnarray}
is the elastic free energy. Here $F(c)$ is a suitable ``double-well'' potential taking its global minimum $0$ at $\pm 1$, for example $F(c)=(1-c^2)^2$. The chemical potential $\mu\colon \Omega \times (0,\infty) \to \mathbb{R}$ is introduced by the first variation of the total energy. 
For a derivation of the Cahn-Larché system we refer to Garcke  \cite{GarckeCH}. Existence and uniqueness results can be found for example in \cite{GarckeCH} and \cite{Garcke00}.

In the Cahn-Larch\'e system above $\eps>0$ is a parameter related to the ``thickness'' of the diffuse interface, i.e., the set $\{|c^\epsilon|<1-\delta\}$ for some small $\delta>0$.
In the following we will study  the sharp interface limit of the diffuse interface model above, i.e., the limit $\epsilon \searrow 0$. We will prove that the solutions of the Cahn-Larché system  converge to a solution of the following sharp interface model, which couples the two-phase Hele-Shaw equation (also called Mullins-Sekerka equation) to the system of linearized elasticity:
\begin{align}
\Delta \mu & = 0 && \mbox{in } \Omega^\pm(t), t>0,\label{sharpsystem1}\\
\Div \, \mathcal{S} &= 0 && \mbox{in }  \Omega^\pm(t), t>0,\\
V & = - \tfrac{1}{2} \left[\nabla \mu\right]_{\Gamma(t)} \cdot \nu && \mbox{on } \Gamma(t), t >0 ,\label{sharpsystem3}\\
\mu &=  \sigma \kappa + \tfrac{1}{2} \nu^T \left[ W \mathrm{Id} - (\nabla \textbf{u})^T \mathcal{S} \right]_{\Gamma(t)} \nu && \mbox{on } \Gamma(t), t >0, \label{sharpsystem4}\\
\left[ \mathcal{S} \nu \right]_{\Gamma(t)} &= \left[ \textbf{u} \right]_{\Gamma(t)} = \left[ \mu \right]_{\Gamma(t)} = 0 && \mbox{on } \Gamma(t), t>0,\label{sharpsystem5}\\
\tfrac{\partial}{\partial n} \mu &= \textbf{u} = 0 && \mbox{on } \partial \Omega, t>0,\\
\Gamma(0) &= \Gamma_{00} && \mbox{for }t=0. \label{sharpsystem7}
\end{align}
Here the components of the alloy fill two disjoint domains $\Omega^+(t),\Omega^-(t) \subset \Omega$ for all times $t \geq 0$. The two domains are separated by a $(d-1)$-dimensional surface $\Gamma(t)$ such that $\Gamma(t) = \partial \Omega^-(t)$. We assume that $\Gamma(t) \subset \Omega$, i.e., we do not consider contact angles. Therefore we obtain $\Omega = \Omega^+(t) \cup \Omega^-(t) \cup \Gamma(t)$. The corresponding elastic energy densities have the form $W_-(\mathcal{E}) : = W(-1,\mathcal{E})$ and $W_+(\mathcal{E}) : = W(1,\mathcal{E})$. Throughout the paper $\nu$ is the unit outer normal of $\partial \Omega^-(t)$, whereas $n$ denotes the unit outer normal of $\partial \Omega$. The normal velocity and the mean curvature of $\Gamma(t)$ are denoted by $V$ and $\kappa$, respectively, taken with respect to $\nu$. The constant $\sigma >0$ describes the surface tension of the interface and $\left[ . \right]_{\Gamma(t)}$ denotes the jump of a quantity across the interface in direction of $\nu$, i.e., $\left[ f \right]_{\Gamma(t)}(x) = \lim_{h \to 0} (f(x+h \nu) - f(x-h\nu))$ for $x \in \Gamma(t)$. The existence of classical solutions is proven in \cite{St}. For the classical Hele-Shaw problem one finds classical solution results in Chen et al.~\cite{ChenEx} and Escher and Simonett \cite{ES}. The global existence of classical solutions and the convergence to spheres are shown in Escher and Simonett \cite{ES2}, provided that the initial value is close to a sphere.

In the case of the Cahn-Hilliard equation, there are two kinds of results for the sharp interface limit. Chen \cite{ChenGlobal} showed the convergence of weak solutions to a varifold solution to the corresponding sharp interface model globally in time. He proved that the family of solutions $\left( c^\epsilon, \mu^\epsilon \right)_{0 < \epsilon \leq 1}$ is weakly compact in some functional spaces. Then he obtained the existence of a convergent subsequence. Garcke and Kwak \cite{GarckeKwak} used this method to show the convergence of the Cahn-Larché system (\ref{system1})-(\ref{system3b}) to the modified Hele-Shaw problem (\ref{sharpsystem1})-(\ref{sharpsystem5}) with Neumann boundary conditions on $\partial \Omega$ and an angle condition for the interface $\Gamma(t)$.  Abels and Röger \cite{AbelsRo} also applied this method to a coupled Navier-Stokes/Cahn-Hilliard system. Recently Abels and Lengeler \cite{AbelsLe} extended this result to fluids with different densities. In all these results convergence holds for arbitrary large times, but the notion of solution to the limit system is rather weak.

 On the other hand, Alikakos et al.\ \cite{ABC} proved convergence of the Cahn-Hilliard equation to the Hele-Shaw problem in strong norms provided that the limit system, which is the two-phase Hele-Shaw system, possesses a smooth solution, which is known to be true at least for a sufficiently small time interval $(0,T)$. By formally matched asymptotic expansions they constructed a family of approximate solutions $\left(c^\epsilon_A,\mu^\epsilon_A \right)_{0<\epsilon \leq 1}$ for the Cahn-Hilliard equation and showed that the difference of the real solution $(c^\epsilon,\mu^\epsilon)$ and approximate solutions converge to $0$ as $\epsilon \searrow 0$ , provided the initial value $c^\epsilon_0$ of the Cahn-Hilliard equation is chosen suitably. Since the zero order expansion of the approximate solutions is based on the solution to the Hele-Shaw problem, they were able to prove the convergence of the Cahn-Hilliard equation to the Hele-Shaw problem as $\epsilon \searrow 0$. Let us mention that Carlen et al.\ \cite{CCO} introduced an alternative method to construct approximate solutions to the Cahn-Hilliard equation. Based on Hilbert expansion they used the ansatz $c^\epsilon(x,t) \approx \sum^N_{i=1}{\epsilon^i c_i(x,\Gamma^{(N)}_t)}$, where $\Gamma^{(N)}_t$ is the $N$th order approximate interface.  An alternative proof of convergence was given by Le~\cite{LeCahnHilliardLimit}, which is based on a gradient flow formulation, $\Gamma$-convergence results and a modified De~Giorgi conjecture, which remains to be verified. For the Cahn-Larché system a formally matched asymptotic expansion was already done in Leo et al.\ \cite{Leo}.

The structure of the article is as follows: In Section \ref{secpre}, we determine the notation and summarize the basic assumptions. Moreover, we mention some spectral analysis results proven by Chen \cite{Chen}. In Section \ref{secdifferenz} we use these results to prove that the difference of approximate and real solutions for the Cahn-Larché system tends to $0$ as $\epsilon \searrow 0$. It is possible to show the convergence in higher norms provided the error terms in the approximate solutions converge to $0$ fast enough. In Section \ref{secapprox}, we construct the approximate solutions by machted asymptotic analysis. Here we follow the method of Alikakos et al.\ \cite{ABC}. In Section \ref{secmain}, we rigorously prove the sharp interface limit for the Cahn-Larché system provided the initial values $c^\epsilon_0$ are chosen suitable. More precisely, we show that the approximate solutions satisfy all the required conditions to apply the statements of the spectral analysis in Section \ref{secpre}. Furthermore the approximate solutions converge to the solution for the modified Hele-Shaw problem (\ref{sharpsystem1})-(\ref{sharpsystem7}) as long as smooth solutions exists for the limit system. Applying the results in Section \ref{secdifferenz}, we finish the proof. \\

\section{Preliminaries}\label{secpre}

\subsection{Notation and Basic Assumptions}

For a sufficient smooth domain $\Omega \subset \mathbb{R}^d$ and an interval $(0,T)$, $T >0$, we define $\Omega_T = \Omega \times (0,T)$ and $\partial_T \Omega = \partial\Omega \times (0,T)$. Moreover, $n$ denotes the exterior unit normal on $\partial \Omega$. For a hypersurface $\Gamma \subset \Omega$ without boundary such that $\Gamma = \partial \Omega^-$ for a reference domain $\Omega^- \subset \Omega$, the interior domain is denoted by $\Omega^-$ and the exterior domain by $\Omega^+ := \Omega \backslash (\Omega^- \cup \Gamma)$, that is $\Gamma$ separates $\Omega$ into an interior and an exterior domain. The exterior unit normal on $\partial \Omega^-$ is denoted by $\nu$. The mean curvature of $\Gamma$ is denoted by $\kappa$ with the sign convention that $\kappa$ is positive, if $\Gamma$ is curved in direction of $\nu$. For a signed distance function $d$ with respect to $\Gamma$, we assume $d < 0$ in $\Omega^-$ and $d>0$ in $\Omega^+ $. By this convention we obtain $\nabla d = \nu$ on $\Gamma$. We denote $a \otimes b = \left( a_i b_j \right)^d_{i,j=1}$ for $a,b \in \mathbb{R}^d$ and $A:B = \sum^d_{i,j=1}{A_{ij} B_{ij}}$ for $A,B \in \mathbb{R}^{d \times d}$.\\
The ``double-well'' potential $F: \mathbb{R} \to \mathbb{R}$ is a smooth function taking its global minimum $0$ at $\pm1$. For its derivative $f(c)=F'(c)$ we assume
\begin{align}
f(\pm 1) & = 0 , & f'(\pm1) & > 0, & \int^u_{-1}{f(s)\, ds} = \int^u_1{f(s)\, ds} & > 0 & \forall u \in (-1,1). \label{bedingungphi1}
\end{align} 
In the following we need an additional assumption
\begin{equation}
c f''(c) \geq 0 \quad \mbox{if } \left| c \right|\geq C_0  \label{eigenschaftf}
\end{equation}
for some constant $C_0>0$. 

\begin{lemma} \label{lemmatheta0}
Let $f \in C^\infty(\mathbb{R})$ be given such that the properties (\ref{bedingungphi1}) hold. Then the problem
\begin{align}
- w'' + f(w) & = 0 \mbox{ in } \mathbb{R}, & w(0) & = 0, & \lim_{z \to \pm\infty}w(z) & = \pm 1  \label{odetheta0}
\end{align}
has a unique solution.\\
In addition, the following properties hold
\begin{align}
w'(z) & > 0 & \forall z & \in \mathbb{R}\,,\label{proptheta1}\\
\left| w^2(z) - 1 \right| + \left| w^{(n)}(z) \right| & \leq C_n e^{-\alpha \left|z\right|} & \forall z & \in \mathbb{R}, \, n \in \mathbb{N} \backslash \left\{0\right\} \label{proptheta2}
\end{align}
for some constants $C_n > 0$, $n \in \mathbb{N} \backslash \left\{0\right\}$, and where $\alpha$ is a fixed constant such that
\[
0 < \alpha < \min \left\{ \sqrt{f'(-1)} ,\sqrt{f'(1)} \right\} \,.
\]
\end{lemma}

\proof See Remark 3.1 in \cite{ABC} or \cite[Lemma 2.6.1]{St}.  
\makebox[1cm]{} \hfill $\Box$\\

The constant elasticity tensor $\mathcal{C} = (\mathcal{C}_{iji'j'})_{i,j,i',j'=1,\ldots,d}$ maps matrices $A \in \mathbb{R}^{d \times d}$ in matrices by the definition
\[
( \mathcal{C} A )_{ij} = \sum^d_{i',j'=1}{\mathcal{C}_{iji'j'} A_{i'j'}}\,.
\]
In addition, we assume the symmetry properties
\[
\mathcal{C}_{iji'j'} = \mathcal{C}_{ijj'i'} = \mathcal{C}_{jii'j'} \quad \mbox{and} \quad \mathcal{C}_{iji'j'} = \mathcal{C}_{i'j'ij}
\]
for all $i,j,i',j'= 1,\ldots ,d$. An important assumption is the positive definiteness of $\mathcal{C}$ on symmetric matrices, that is, there exists some constant $c_2>0$ such that 
\begin{equation}
A:\mathcal{C} A \geq 2 c_2 \left| \mathrm{sym}(A) \right|^2 \quad \forall A \in \mathbb{R}^{d \times d}\,. \label{posdefC}
\end{equation}
An important consequence of the positive definiteness is the following lemma.

\begin{lemma}
Let the tensor $\mathcal{C}$ be defined as above. Then it holds for all $a,b \in \mathbb{R}^d$
\begin{equation}
\left(a \otimes b\right): \mathcal{C} \left( a \otimes b \right) \geq  c_2 \left| a \otimes b \right|^2 \,.
\end{equation}
\end{lemma}

\proof One can easily verify the assertion by direct calculation.
\makebox[1cm]{} \hfill $\Box$\\

\subsection{Spectral Analysis}

In this subsection we summarize some results proven by Chen \cite{Chen}, which will be essential for the convergence proof. \\
Let $\gamma \subset \Omega$ be a smooth $(d-1)$-dimensional manifold without boundary and let $r = r(x)$ be the signed distance function satisfying $r<0$ inside $\gamma$ and $r>0$ outside $\gamma$. Let $s=s(x)$ be the orthogonal projection of $x$ on $\gamma$. Then there exits $\delta_0>0$ such that $\gamma(2\delta_0) := \left\{ x \in \mathbb{R}^d : \left| r(x) \right| < 2 \delta_0 \right\} \subset \Omega$ and such that $\tau : \gamma(2 \delta_0) \to (-2 \delta_0,2\delta_0) \times \gamma$ defined by $\tau(x) = (r(x),s(x))$ is a smooth diffeomorphism where $\delta_0$ only depends on $\gamma$ and $\partial \Omega$, cf. \cite[Kapitel 4.6]{Hildebrandt2}. Let $\phi^\epsilon : \Omega \to \mathbb{R}$ be a given function with the expansion
\begin{eqnarray}
\phi^\epsilon(x) & = & \zeta \! \left( \tfrac{r(x)}{\delta_0} \right) \left( \theta_0 \! \left( \tfrac{r(x)}{\epsilon} \right) + \epsilon p^\epsilon(s(x)) \theta_1 \! \left( \tfrac{r(x)}{\epsilon} \right) + \epsilon^2 q^\epsilon(x) \right) \nonumber \\
&& + \left( 1 - \zeta \! \left( \tfrac{r(x)}{\delta_0} \right) \right) \left( \phi^+_\epsilon(x) \chi_{\left\{ r(x) >0 \right\}} + \phi^-_\epsilon(x) \chi_{\left\{ r(x) < 0 \right\}} \right) , \label{bedingungphi2}
\end{eqnarray}
where $\zeta \in C^\infty_0(\mathbb{R})$ is a cut-off function such that
\begin{align}
\zeta(z) & = 1 \mbox{ if } \left| z \right| < \frac{1}{2} , & \zeta(z) & = 0 \mbox{ if } \left| z \right| > 1 , & z \zeta'(z) \leq 0 \mbox{ in } \mathbb{R} \,, \label{cutoffzeta}
\end{align}
$\theta_0$ is the unique solution to 
\begin{align}
- \theta''_0 + f(\theta_0) & = 0 \mbox{ in } \mathbb{R}, & \theta_0(0) & = 0, & \theta_0(\pm \infty) & = \pm 1 \,, \label{theta0}
\end{align}
$\theta_1 \in C^1(\mathbb{R}) \cap L^\infty(\mathbb{R})$ is any function satisfying
\begin{equation}
\int_{\mathbb{R}}{ \theta_1 (\theta'_0)^2 f''(\theta_0)} = 0 \,, \label{theta1}
\end{equation}
and $p^\epsilon(x)$, $q^\epsilon(x)$, $\phi^+_\epsilon$, and $\phi^-_\epsilon(x)$ are smooth function satisfying 
\begin{alignat}{2}
\sup_{\epsilon \in (0,1]} \left| p^\epsilon \right| + \frac{\epsilon}{\epsilon + \left| r \right|} \left| q^\epsilon \right| & \leq C_\ast && \mbox{in } \gamma(\delta_0) \,,\label{bedingungphi6}\\
\sup_{\epsilon \in (0,1]} \left| \nabla^\gamma \phi^\epsilon \right| & \leq C_\ast && \mbox{in } \gamma(\delta_0) \,,\label{bedingungphi7}\\
\pm \phi^\pm_\epsilon > 0,\, f'(\phi^\pm_\epsilon) & \geq 1/C_\ast \quad && \mbox{in }\Omega  \label{bedingungphi8}
\end{alignat}
for some constant $C_\ast >0$ where $\nabla^\gamma = \nabla - \nabla r (\nabla r \cdot \nabla)$ is the tangential gradient along $\gamma$. With these conditions we obtain the following proposition.

\begin{proposition} \label{propChen}
Assume that (\ref{bedingungphi1}) and (\ref{bedingungphi2})-(\ref{bedingungphi8}) hold. Then, for any given $\gamma_1>0$, there exist constants $\epsilon_0>0$ and $C > 0$ which depend only on $f$, $\theta_1$, $C_\ast$, $\Omega$, $\gamma_1$, and the $C^3$ norm of $\gamma$ such that for every $\epsilon \in (0,\epsilon_0]$, $w \in H^1_{(0)}(\Omega) \backslash \left\{0\right\}$ , and $\Psi \in H^2(\Omega)$ with $-\Delta \Psi = w$ and $\left.\frac{\partial}{\partial n} \Psi \right|_{\partial \Omega} = 0$, the following inequality holds
\begin{equation}
\int_\Omega{\left( \epsilon \left| \nabla w \right|^2 + \epsilon^{-1} f'(\phi^\epsilon) w^2 \right)} \geq - C \left\| \nabla \Psi\right\|^2_{L^2(\Omega)} + \gamma_1 \epsilon \left\| w \right\|^2_{L^2(\Omega)} \,.
\end{equation}
\end{proposition}

\proof Let $\gamma_1$ be any positive constant and let $w$ and $\Psi$ be any given functions as above. First we consider the case $\int_\Omega{\epsilon \left|\nabla w\right|^2 + \epsilon^{-1} f'(\phi^\epsilon) w^2} \leq \gamma_1 \epsilon \left\|w\right\|^2_{L^2(\Omega)}$. Then by \cite[Theorem 3.1.]{Chen}, there exists a constant $C = C(\gamma_1)$ and $\epsilon_1 > 0$ such that for all $\epsilon \in (0,\epsilon_1]$
\[
\epsilon \left\| w \right\|^2_{L^2(\Omega)} \leq C \left\| \nabla \Psi \right\|^2_{L^2(\Omega)}\,.
\]
Due to the spectral estimate \cite[Theorem 1.1.]{Chen}, there exist a constant $C_1$ depending on $f$, $\theta_1$, $C_\ast$, $\Omega$, and the $C^3$ norm of $\gamma$ such that
\begin{eqnarray*}
\int_\Omega{\left( \epsilon \left| \nabla w \right|^2 + \epsilon^{-1} f'(\phi^\epsilon) w^2 \right)} & \geq & - C_1 \left\| \nabla \Psi\right\|^2_{L^2(\Omega)} \\
& \geq & - (C_1 + C \gamma_1) \left\| \nabla \Psi \right\|^2_{L^2(\Omega)} + \gamma_1 \epsilon \left\| w \right\|^2_{L^2(\Omega)} \,.
\end{eqnarray*}
Hence together with the other case $\int_\Omega{\epsilon \left|\nabla w\right|^2 + \epsilon^{-1} f'(\phi^\epsilon) w^2} \geq \gamma_1 \epsilon \left\|w\right\|^2_{L^2(\Omega)}$, the assertion of the lemma follows.
\makebox[1cm]{} \hfill $\Box$\\

\begin{remark} \label{remarkC3norm}
The $C^3$ norm of $\gamma$ is defined as follows: Let $r$ be the signed distance function to $\gamma$ and $M = \sup_{x \in \gamma} \left| D^2r(x) \right|$. Then $r$ is smooth in $\gamma(1/M)$, cf. \cite[Kapitel 4.6]{Hildebrandt2}. So we define $\left\| \gamma \right\|_{C^3} = M + \left\| r \right\|_{C^3(\gamma(\delta))}$ where $\delta= \min\left\{ 1 , 1/2 M \right\}$. 
\end{remark}

\section{Convergence of the Difference of Approximate
and True Solutions} \label{secdifferenz}

In this section we want to show that the difference of approximate solutions having certain properties and true solutions converges to zero as $\epsilon \searrow 0$. 

\begin{theorem} \label{theorem1}
Let $\left\{ c^\epsilon_A, \mu^\epsilon_A, \textbf{u}^\epsilon_A \right\}_{0 < \epsilon \leq 1}$ be a family of functions in the function space $C^\infty (\overline{\Omega_T}) \times C^\infty(\overline{\Omega_T}) \times C^\infty(\overline{\Omega_T}; \mathbb{R}^d)$ satisfying the system of differential equations
\begin{eqnarray}
(c^\epsilon_A)_t = \Delta \mu_A^\epsilon &\mbox{in}&\Omega_T\, \label{systemA1},\\
\mu_A^\epsilon = - \epsilon \Delta c_A^\epsilon + \epsilon^{-1} f(c_A^\epsilon)+ W_{,c}(c^\epsilon_A \, , \mathcal{E}(\textbf{u}^\epsilon_A)) + r_A^\epsilon &\mbox{in}&\Omega_T\, \label{systemA2},\\
\mathrm{div} \, W_{, \mathcal{E}}(c^\epsilon_A, \mathcal{E}(\textbf{u}^\epsilon_A)) = \textbf{s}^\epsilon_A &\mbox{in}&\Omega_T\, ,\label{systemA3}\\
\tfrac{\partial}{\partial n}c_A^\epsilon = \tfrac{\partial}{\partial n}\mu_A^\epsilon=0 &\mbox{on}& \partial_T \Omega \, ,\\
\textbf{u}^\epsilon_A = 0 &\mbox{on}& \partial_T\Omega \, \label{systemA6},
\end{eqnarray}
where $r^\epsilon_A = r^\epsilon_A(x,t)$ and $\textbf{s}^\epsilon_A = \textbf{s}^\epsilon_A(x,t)$ are functions such that
\begin{eqnarray}
\left\| r^\epsilon_A \right\|^2_{L^2(\Omega_T)} + \left\| \textbf{s}^\epsilon_A \right\|^2_{L^2(\Omega_T)} \leq \frac{1}{2} \epsilon^{p k}\,, \label{fehlerabsch}
\end{eqnarray}
$p = \frac{2 \left( d+4 \right)}{d+2}$, and $k \in \mathbb{N}$ such that
\begin{equation}
k > \frac{\left( 4d + 10 \right) \left( d+2 \right)}{4 \left( d+4 \right)}\,. \label{k}
\end{equation}
Also assume that $c^\epsilon_A$ satisfies the boundedness condition
\begin{eqnarray}
\sup_{0 < \epsilon \leq 1} \left\| c^\epsilon_A \right\|_{L^\infty(\Omega_T)} \leq C_0 \label{uniboundcA}
\end{eqnarray}
for some $C_0>0$, the energy density $f$ satisfies (\ref{bedingungphi1}) and (\ref{eigenschaftf}), and 
\begin{eqnarray}
\phi^\epsilon_t(.) := c^\epsilon_A(.,t) \label{phiepsilon}
\end{eqnarray}
has the form (\ref{bedingungphi2}). Let $(c^\epsilon, \mu^\epsilon, \textbf{u}^\epsilon)$ be the unique solution to (\ref{system1})-(\ref{system6}) with $c^\epsilon_0(x) = c^\epsilon_A(x,0)$ in $\Omega$. Then there exists a constant $\epsilon_0 = \epsilon_0(C_0,T,\Omega,k,d) \in (0,1]$ such that if $\epsilon \in (0,\epsilon_0)$, then
\begin{eqnarray*}
\left\| c^\epsilon - c^\epsilon_A \right\|_{L^p(\Omega_T)} + \left\| \textbf{u}^\epsilon - \textbf{u}^\epsilon_A \right\|_{L^2(0,T;W^1_2(\Omega))} \leq C \epsilon^k 
\end{eqnarray*}
for some $C>0$ independent of $\epsilon$.
\end{theorem}

\proof Let $R = c^\epsilon - c^\epsilon_A$ and $\textbf{u}= \textbf{u}^\epsilon - \textbf{u}^\epsilon_A$ be the remainder functions.\\
It holds for all $t \in [0,T]$ 
\begin{eqnarray*}
\int_\Omega{R(.,t)\, dx} = \int^t_0{\int_\Omega{\frac{\partial}{\partial t} \left(c^\epsilon- c^\epsilon_A\right) dx} \, ds} = \int^t_0{\int_\Omega{\Delta \left(\mu^\epsilon - \mu^\epsilon_A\right) dx} \, ds} = 0
\end{eqnarray*}
due to the Neumann boundary condition. Hence there exists a unique smooth solution $\Psi(x,t)$ to the Neumann boundary Problem  
\begin{align*}
- \Delta \Psi(.,t) = R(.,t) \mbox{ in }\Omega, && \frac{\partial}{\partial n} \Psi(.,t)=0 \mbox{ on } \partial \Omega, && \int_\Omega{\Psi(.,t)\, dx} =0 
\end{align*}
for all $t \in [0,T]$. % This can be seen as follows. Applying Poincaré's inequality and the Lax-Milgram theorem, we obtain a unique weak solution $\Psi(.,t)$ in the space $\left\{c \in W^1_2(\Omega) : \int_\Omega{c\,dx} = 0 \right\}$. An easy calculation shows that $\Psi(.,)$ is also a weak solution in $W^1_2(\Omega)$. Then by applying the usual regularity theory, we prove that $\Psi(.,t)$ is smooth, cf. \cite[Chapter 4]{McLean}. Smoothness with respect to time $t$ follows from the smoothness of $R$.  \\
Multiplying the equation $\partial_t R - \Delta \left(\mu^\epsilon-\mu^\epsilon_A\right) = 0 $ by $\Psi$ and integrating over $\Omega$ yields
\begin{eqnarray}
0 & = & \int_\Omega{ \Psi \left( \partial_t R - \Delta \left(\mu^\epsilon-\mu^\epsilon_A\right)\right) dx}\nonumber\\
& = & \int_\Omega{ \Psi \left(- \Delta \partial_t \Psi\right) - \Delta \Psi \left(\mu^\epsilon-\mu^\epsilon_A\right) dx}\nonumber \\
&=& \frac{1}{2} \frac{d}{dt} \int_\Omega{\left|\nabla \Psi\right|^2 dx} - \int_\Omega{ R \left(\epsilon \Delta R - \epsilon^{-1} \left(f(c^\epsilon) - f(c_A^\epsilon)\right) -  W_{,c}(R,\mathcal{\textbf{u}})+ r_A^\epsilon\right) dx}\nonumber\\
&=& \frac{1}{2} \frac{d}{dt} \int_\Omega{\left|\nabla \Psi\right|^2 dx} + \int_\Omega{ \epsilon \left|\nabla R\right|^2 + \epsilon^{-1} f'(c^\epsilon_A) R^2 \, dx }\nonumber\\
&& + \int_\Omega{ \epsilon^{-1} \mathcal{N}(c^\epsilon_A,R) R +  W_{,c}(R,\mathcal{\textbf{u}}) R - r_A^\epsilon R \, dx},\label{gl1a}
\end{eqnarray}
where we have used the Neumann boundary conditions for $\Psi$, $R$ and $\mu^\epsilon-\mu^\epsilon_A$ on $\partial \Omega$, the expressions for $\mu^\epsilon$ and $\mu^\epsilon_A$ and $-\Delta \Psi = R$. Here $\mathcal{N}(.,.)$ is defined by $\mathcal{N}(c,R):= f(c+R) - f(c) - f'(c)R$.\\
Applying Lemma 2.2 in \cite{ABC} yields
\begin{equation}
- \int_{\Omega}{\epsilon^{-1} \mathcal{N}(c^\epsilon_A,R)R \, dx} \leq C \epsilon^{-1} \left\| R \right\|^p_{L^p(\Omega)} \label{resttaylor}
\end{equation}
for some constant $C=C(p)>0$. By integration by parts and using (\ref{system3}) and (\ref{systemA3}), and the symmetry of $\mathcal{C}$, we obtain
\begin{eqnarray*}
\int_\Omega{\mathcal{E}(\textbf{v}) : \mathcal{C} (\mathcal{E}(\textbf{u}) - \mathcal{E}^\star R )\, dx} = \int_\Omega{\textbf{v} \cdot \textbf{s}^\epsilon_A \, dx} \quad \forall \textbf{v} \in H^1_0(\Omega)^d\,.
\end{eqnarray*}	
For $ \textbf{v} = \textbf{u} \in H^1_0(\Omega)^d$ this equation yields since $ W_{,c}(R,\mathcal{E}(\textbf{u})) R = - \mathcal{E}^\star R : \mathcal{C} (\mathcal{E}(\textbf{u}) - \mathcal{E}^\star R)$
\begin{eqnarray}
\lefteqn{ \int_\Omega{ W_{,c}(R,\mathcal{E}(\textbf{u}))} R \, dx} \nonumber \\
& = & \int_\Omega{(\mathcal{E}(\textbf{u}) - \mathcal{E}^\star R) : \mathcal{C} (\mathcal{E}(\textbf{u}) - \mathcal{E}^\star R) \, dx} - \int_\Omega{\textbf{u} \cdot \textbf{s}^\epsilon_A \, dx} \nonumber \\
& \geq & c_2 \left\| \mathcal{E}(\textbf{u}) - \mathcal{E}^\star R \right\|^2_{L^2(\Omega)} - \left( \left\| \mathcal{E}(\textbf{u}) - \mathcal{E}^\star R \right\|_{L^2(\Omega)} + \left\| \mathcal{E}^\star R \right\|_{L^2(\Omega)} \right) \left\| \textbf{s}^\epsilon_A \right\|_{L^2(\Omega)} \nonumber \\
& \geq & \frac{c_2}{2} \left\| \mathcal{E}(\textbf{u}) - \mathcal{E}^\star R \right\|^2_{L^2(\Omega)} - C \left\| \textbf{s}^\epsilon_A \right\|^2_{L^2(\Omega)} - \left\| \mathcal{E}^\star R \right\|_{L^2(\Omega)} \left\| \textbf{s}^\epsilon_A \right\|_{L^2(\Omega)} \, ,
\end{eqnarray}
where we have used the Korn (cf. Roub{\'{\i}}{\v{c}}ek \cite{Roubicek}) and triangle inequality in the first estimate and Young's inequality in the second estimate. Hölder's inequality gives us the estimate
\begin{equation}
\int_\Omega{ r_A^\epsilon R \, dx} \leq \left\|r_A^\epsilon\right\|_{L^2(\Omega)} \left\| R \right\|_{L^2(\Omega)}\,.
\end{equation}
By Proposition \ref{propChen}, there exists some constant $C>0$ such that
\begin{equation}
\int_\Omega{\epsilon \left|\nabla R\right|^2 + \epsilon^{-1} f'(c^\epsilon_A)R^2 \, dx} \geq - C \left\|\nabla \Psi\right\|^2_{L^2(\Omega)} + 2 \epsilon \left\|R\right\|^2_{L^2(\Omega)}. \label{spectralabsch}
\end{equation}
Therefore equation (\ref{gl1a}) provides us with (\ref{resttaylor})-(\ref{spectralabsch}) and Young's inequality 
\begin{eqnarray}
\lefteqn{ \frac{1}{2} \frac{d}{dt} \int_\Omega{\left|\nabla \Psi\right|^2 dx } + \epsilon \left\| R \right\|^2_{L^2(\Omega)} + \frac{c_2}{2} \left\| \mathcal{E}(\textbf{u}) - \mathcal{E}^\star R \right\|^2_{L^2(\Omega)} } \nonumber \\
& \leq & C \left(\left\|\nabla \Psi\right\|^2_{L^2(\Omega)} + \epsilon^{-1} \left\| R \right\|^p_{L^p(\Omega)} + \epsilon^{-1} \left\|r_A^\epsilon\right\|^2_{L^2(\Omega)} + \epsilon^{-1} \left\| \textbf{s}^\epsilon_A \right\|_{L^2(\Omega)}\right)  \label{ab2a}
\end{eqnarray}
for all $t \in(0,T]$ and some $C=C(f,C_0,p)$. Note that $R(.,0)=0$ and therefore $\Psi(.,0)=0$. Then applying Gronwall's inequality yields
\begin{eqnarray}
\sup_{0 \leq \tau \leq t} \left\|\nabla \Psi(.,\tau) \right\|^2_{L^2(\Omega)}  \leq  C \epsilon^{-1} e^{C t} \left(   \left\|r_A^\epsilon\right\|^2_{L^2(\Omega_t)} +  \left\| \textbf{s}^\epsilon_A \right\|^2_{L^2(\Omega_t)} + \left\| R \right\|^p_{L^p(\Omega_t)} \right)  .\label{gron1}
\end{eqnarray}
Integrating inequality (\ref{ab2a}) over $(0,t)$, $t \in [0,T]$, and using (\ref{gron1}) yields
\begin{eqnarray}
 \lefteqn{\epsilon \left\| R \right\|^2_{L^2(\Omega_t)} + \frac{c_2}{2} \left\| \mathcal{E}(\textbf{u}) - \mathcal{E}^\star R \right\|^2_{L^2(\Omega_t)}} \nonumber  \\
 & \leq & C_1 \epsilon^{-1} \left(  \left\|r_A^\epsilon\right\|^2_{L^2(\Omega_t)} +  \left\| \textbf{s}^\epsilon_A \right\|^2_{L^2(\Omega_t)} +  \left\| R \right\|^p_{L^p(\Omega_t)} \right)  \label{L2estimate}
\end{eqnarray}
for some constant $C_1 = C_1(T)>0$ independent of $\epsilon$ and $t \in (0,T]$. Integrating (\ref{gl1a}) over $(0,t)$ yields
\begin{eqnarray}
\epsilon \left\| \nabla R \right\|^2_{L^2(\Omega_t)} \leq C \left( \epsilon \left\| r^\epsilon_A \right\|^2_{L^2(\Omega_t)} + \left\| \textbf{s}^\epsilon_A \right\|^2_{L^2(\Omega_t)} + \epsilon^{-1} \left\| R \right\|^2_{L^2(\Omega_t)} + \epsilon^{-1} \left\| R \right\|^p_{L^p(\Omega_t)} \right). \label{nablaRestimate}
\end{eqnarray}
Since $\int_\Omega{R(.,t)} = 0$, the Gagliardo-Nirenberg inequality (cf. Roub{\'{\i}}{\v{c}}ek \cite{Roubicek}) and a Poincaré's inequality  implies for $p = 2\frac{d+4}{d+2}$ and for all $t \in [0,T]$
\[
\left\| R \right\|^p_{L^p(\Omega)} \leq C \left\| R \right\|^{\frac{8}{d+2}}_{L^2(\Omega)} \left\| \nabla R \right\|^{\frac{2d}{d+2}}_{L^2(\Omega)} \,.
\]
To estimate $\left\| R \right\|_{L^2(\Omega)}$ on the right-hand side, we use integration by parts to get
\[
\left\| R \right\|^2_{L^2(\Omega)} = - \int_\Omega{R \Delta \Psi\, dx} = \int_\Omega{\nabla R \cdot \nabla \Psi \, dx} \leq \left\| \nabla R \right\|_{L^2(\Omega)} \left\| \nabla \Psi \right\|_{L^2(\Omega)} \,,
\]
where we have used the Neumann boundary conditions for $c^\epsilon$ and $c^\epsilon_A$ on $\partial \Omega$. Therefore it holds
\begin{eqnarray}
\left\| R \right\|^p_{L^p(\Omega_t)} & \leq & C \int^t_0{\left\| \nabla \Psi(.,\tau) \right\|^{\frac{4}{d+2}}_{L^2(\Omega)} \left\| \nabla R(.,\tau) \right\|^2_{L^2(\Omega)} \, d\tau} \nonumber \\
& \leq & C \sup_{0 \leq \tau \leq t} \left\| \nabla \Psi(.,\tau) \right\|^{\frac{4}{d+2}}_{L^2(\Omega)} \left\| \nabla R \right\|^2_{L^2(\Omega_t)} \,. \label{Lpinequ}
\end{eqnarray}
In order to treat the term $\left\| R \right\|_{L^2(\Omega_t)}$ in (\ref{nablaRestimate}), we define two sets
\begin{eqnarray*}
A^\epsilon_1 & : = & \left\{ t \in [0,T] : \epsilon \left\| R \right\|^2_{L^2(\Omega_t)} > 2 C_1 \epsilon^{-1} \left\| R \right\|^p_{L^p(\Omega_t)} \right\} , \\
A^\epsilon_2 & : = & \left\{ t \in [0,T] : \epsilon \left\| R \right\|^2_{L^2(\Omega_t)} \leq 2 C_1 \epsilon^{-1} \left\| R \right\|^p_{L^p(\Omega_t)} \right\} ,
\end{eqnarray*}
where $C_1$ is the same constant as in (\ref{L2estimate}). Furthermore, we set 
\[
T^\epsilon : = \sup \left\{ t \in (0,T] : \left\| R \right\|_{L^p(\Omega_t)} \leq \epsilon^k \right\}.
\]
\textbf{1st case:} ``$T^\epsilon \in A^\epsilon_1$'' \\
Then the definition of $A^\epsilon_1$ and (\ref{L2estimate}) yield
\begin{eqnarray*}
\left\| R \right\|^p_{L^p(\Omega_{T^\epsilon})} \leq \frac{\epsilon^2}{2 C_1} \left\| R \right\|^2_{L^2(\Omega_{T^\epsilon})} \leq \frac{1}{2} \left( \left\| r^\epsilon_A \right\|^2_{L^2(\Omega_{T^\epsilon})} + \left\| \textbf{s}^\epsilon_A \right\|^2_{L^2(\Omega_{T^\epsilon})} + \left\| R \right\|^p_{L^p(\Omega_{T^\epsilon})} \right).
\end{eqnarray*}
Therefore we get
\begin{eqnarray}
\left\| R \right\|^p_{L^p(\Omega_{T^\epsilon})} \leq \left\| r^\epsilon_A \right\|^2_{L^2(\Omega_{T^\epsilon})} + \left\| \textbf{s}^\epsilon_A \right\|^2_{L^2(\Omega_{T^\epsilon})} \leq \frac{1}{2} \epsilon^{pk}\,,
\end{eqnarray}
where we have used (\ref{fehlerabsch}). Hence we conclude $T^\epsilon = T$ by definition of $T^\epsilon$.\\
\textbf{2nd case:} ``$T^\epsilon \in A^\epsilon_2$'' \\
We use inequality (\ref{Lpinequ}) and apply (\ref{gron1}), (\ref{nablaRestimate}), and the definition of $A^\epsilon_2$ to obtain
\begin{eqnarray*}
\left\| R \right\|^p_{L^p(\Omega_{T^\epsilon})} & \leq & C \left( \epsilon^{-1} \left\| r^\epsilon_A \right\|^2_{L^2(\Omega_{T^\epsilon})} + \epsilon^{-1} \left\| \textbf{s}^\epsilon_A \right\|^2_{L^2(\Omega_{T^\epsilon})} + \epsilon^{-1} \left\| R \right\|^p_{L^p(\Omega_{T^\epsilon})} \right)^{\frac{2}{d+2}} \nonumber \\
&& \times \left( \left\| r^\epsilon_A \right\|^2_{L^2(\Omega_{T^\epsilon})} + \epsilon^{-1} \left\| \textbf{s}^\epsilon_A \right\|^2_{L^2(\Omega_{T^\epsilon})} + \epsilon^{-4} \left\| R \right\|^p_{L^p(\Omega_{T^\epsilon})} \right).
\end{eqnarray*}
Applying (\ref{fehlerabsch}) and $\left\| R \right\|_{L^p(\Omega_{T^\epsilon})} \leq \epsilon^k$, yields
\begin{eqnarray*}
\left\| R \right\|^p_{L^p(\Omega_{T^\epsilon})} \leq C \epsilon^{\left(-1 +pk\right) \frac{2}{d+2} } \epsilon^{-4 + pk} \leq  C \epsilon^{pk} \epsilon^{ \frac{4\left(d+4\right)}{\left(d+2\right)^2} \left(k - \frac{(4d+10)(d+2)}{4(d+4)}\right)} \,,
\end{eqnarray*}
where we have used the definition of $p$ in the second inequality. By assumption we have $k > \frac{(4d+10)(d+2)}{4(d+4)}$. Therefore there exists $\epsilon_0 \in (0,1]$ such that for all $\epsilon \in (0,\epsilon_0]$, it holds
\[
\left\| R \right\|^p_{L^p(\Omega_{T^\epsilon})} \leq \frac{1}{2} \epsilon^{pk}\,,
\]
provided $T^\epsilon \in A^\epsilon_2$. As in the 1st case, it follows $T^\epsilon = T$. \\
The estimate for $\textbf{u}$ follows from
\[
\left\| \mathcal{E}(\textbf{u}) \right\|_{L^2(\Omega_T)} \leq \left\| \mathcal{E}(\textbf{u}) - \mathcal{E}^\star R \right\|_{L^2(\Omega_T)} + \left\| \mathcal{E}^\star R \right\|_{L^2(\Omega_T)} \leq C \epsilon^{\frac{pk-1}{2}} + C \epsilon^k \leq C \epsilon^k
\]
by (\ref{L2estimate}) and since $pk-1 \geq 2k$ (use the definitions of $p$ and $k$). Thus the assertion follows due to the Korn inequality.
\makebox[1cm]{} \hfill $\Box$\\

We even get assertions in stronger norms.

\begin{theorem} \label{hoeregR}
Let the assumptions of Theorem \ref{theorem1} hold. Let $m>0$ be any fixed integer and assume $\left\| c^\epsilon_A \right\|_{W^{m+l+1}_2(\Omega_T)} + \left\| \mu^\epsilon_A \right\|_{W^{m+l-1}_2(\Omega_T)} + \left\| \textbf{u}^\epsilon_A \right\|_{W^{m+l+1}_2(\Omega_T)} \leq \epsilon^{-K(m)}$ for $l>\frac{d+1}{2}$, some integer $K(m)$, and all small $\epsilon>0$. If $k$ in (\ref{fehlerabsch}) is large enough, then
\[
\left\|c^\epsilon - c^\epsilon_A \right\|_{C^{m}(\Omega_{T})} + \left\|\mu^\epsilon - \mu^\epsilon_A \right\|_{C^{m-2}(\Omega_{T})} + \left\|\textbf{u}^\epsilon - \textbf{u}^\epsilon_A \right\|_{C^{m+1}(\Omega_{T})} \leq \epsilon 
\]
for all sufficiently small $\epsilon >0$.
\end{theorem}

\proof We show the assertion in the same way as in \cite[Theorem 2.3.]{ABC}. \\
%Since for every $ m \in \mathbb{N}$
%\[
%W^{m+l}_2(\Omega_T) \hookrightarrow C^m(\overline{\Omega_T}) \quad \mbox{if }l > \frac{d+1}{2}
%\]
%and since 
%\[
%\left( L^2(\Omega_T), W^{m+l+1}_2(\Omega_T) \right)_{\theta,2} = W^{m+l}_2(\Omega_T) \quad \mbox{for } \theta= \frac{m+l}{m+l+1} \, ,
%\]
%cf. Section \ref{secinterpol}, it follows
%\begin{equation}
%\left\| c^\epsilon - c^\epsilon_A \right\|_{C^m(\Omega_T)} \leq C \left\| c^\epsilon -c^\epsilon_A \right\|^{\theta}_{L^2(\Omega_T)} \left\| c^\epsilon -c^\epsilon_A \right\|^{1-\theta}_{W^{m+l+1}_2(\Omega_T)}  \label{c-cACm}
%\end{equation}
%for some $C>0$. 
This means by an interpolation argument and Theorem \ref{theorem1}, it is sufficient to show
\[
\left\| c^\epsilon \right\|_{W^{l}_2(\Omega_T)} \leq \epsilon^{-K(l)}
\]
for some integer $l,K(l)>0$ if $k$ in (\ref{fehlerabsch}) is large enough. % The estimate for $\mu^\epsilon - \mu^\epsilon_A$ follows from the equations for the chemical potential (\ref{system2}) and (\ref{systemA2}) (see end of the proof). \\
%We replace $f$ by $\bar{f}$ such that $f = \bar{f}$ in $(-\frac{3}{2} C_0, \frac{3}{2} C_0)$ and $\bar{f}(c)$ is linear when $\left|c\right|>2 C_0$ where $C_0$ is the same constant as in (\ref{uniboundcA}). Denote by $(\bar{c}^\epsilon, \bar{\textbf{u}}^\epsilon) $ the solution to the modified system with $\bar{f}$. Define $A : D(A) \to L^p(\Omega)$ by $A= -\Delta + \mathrm{Id}$ with $ D(A) = \left\{ c \in W^{2}_p(\Omega) : \left.\frac{\partial}{\partial n} c \right|_{\partial \Omega} = 0 \right\}$. W.l.o.g. we assume that in Lemma \ref{semigrle} the constant $\tau = 1$ since we only consider a finite number of different $p$'s. Otherwise we replace $A$ by $-\Delta +c \mathrm{Id}$ for some $c \in \mathbb{R}$. Therefore $-A^2$ is sectorial with domain $D(-A^2) = \left\{ c \in W^{4}_p(\Omega) : \left.\frac{\partial}{\partial n} c \right|_{\partial \Omega} = \left.\frac{\partial}{\partial n} A c \right|_{\partial \Omega} = 0 \right\}$ and $W^1_p(\Omega) = D(A^{\frac{1}{2}})$ with equivalent norms, see Section \ref{secsemigr}.\\
The main difference to \cite[Theorem 2.3]{ABC} is that in general $c^\epsilon \notin D(-A^2)= \left\{ c \in W^{4}_p(\Omega) : \left.\frac{\partial}{\partial n} c \right|_{\partial \Omega} = \left.\frac{\partial}{\partial n} A c \right|_{\partial \Omega} = 0 \right\}$ where $A = -\Delta + \mathrm{Id}$ with $ D(A) = \left\{ c \in W^{2}_p(\Omega) : \left.\frac{\partial}{\partial n} c \right|_{\partial \Omega} = 0 \right\}$. Therefore we add a function $\Psi$ such that $c^\epsilon + \epsilon^{-1} \Psi \in D(-A^2)$. Define $\Psi(x,t)$ as the unique solution to the Neumann boundary problem
\begin{align}
- \Delta \Psi(.,t) & =  W_{,c}(c^\epsilon, \mathcal{E}(\textbf{u}^\epsilon))(.,t) - \frac{1}{\left| \Omega \right|} \int_{\Omega}{ W_{,c}(c^\epsilon, \mathcal{E}(\textbf{u}^\epsilon))(x,t) \, dx } & \mbox{in }\Omega \,, \label{equPsi}\\
 \frac{\partial}{\partial n} \Psi(.,t)& = 0 \; \mbox{ on } \partial \Omega\,, \quad  \int_\Omega{\Psi(x,t) \,dx} = 0 \,.
\end{align}
Then it is not difficult to verify that $c^\epsilon + \epsilon^{-1} \Psi \in D(-A^2)$. By semigroup theory we can estimate $c^\epsilon + \epsilon^{-1} \Psi \in D(-A^2)$ in the same way as in \cite[Theorem 2.3]{ABC}. Also we can control $c^\epsilon$ by the definition of $\Psi$. For more details we refer to \cite{St}. 
\makebox[1cm]{} \hfill $\Box$\\

\section{Approximate Solution} \label{secapprox}

In this section we use a matched asymptotic expansion as in \cite{ABC} to construct a family of approximate solutions $\left\{c^\epsilon_A,\mu^\epsilon_A,\textbf{u}^\epsilon_A\right\}_{0 < \epsilon \leq 1}$ satisfying (\ref{systemA1})-(\ref{phiepsilon}). \\
Unlike \cite{ABC} we first define a family of approximate solutions and then we verify that this family satisfies (\ref{systemA1})-(\ref{phiepsilon}) and motivate the definition.

\subsection{Definition of an Approximate Solution} \label{secapp}

Let $\Gamma_{00} \subset\subset \Omega$ be a given smooth hypersurface without boundary and assume that the Hele-Shaw problem (\ref{sharpsystem1})-(\ref{sharpsystem7}) admits a smooth solution $(\mu,\textbf{u}, \Gamma)$ in the time interval $[0,T]$. Let $d^0$ be the signed distance from $x$ to $\Gamma_t$ such that $d^0<0$ inside of $\Gamma_t$, and let $\delta$ be a small positive constant such that $\mbox{dist}\left(\Gamma_t, \partial \Omega\right)> 2\delta$ for all $t \in [0,T]$. The distance function $d^0$ is smooth in $\Gamma (2\delta):=\left\{(x,t) \in \Omega_T : \left|d^0\right|<2\delta\right\}$ and $\mu^\pm := \left.\mu\right|_{Q^\pm_0}$ and $\textbf{u}^\pm := \left.\textbf{u}\right|_{Q^\pm_0}$ have a smooth extension to $Q^\pm_0 \cup \Gamma(2\delta)$ where $Q^\pm_0:=\left\{(x,t) \in \Omega_T : \pm d^0 > 0 \right\}$. Furthermore, $d_B(x)$ is the signed distance function from $x$ to $\partial \Omega$ such that $d_B<0$ in $\Omega$ and $S_B : \left\{x \in \mathbb{R}^d : \left|d_B\right| \leq \delta \right\} \to \partial \Omega$ is the orthogonal projection of $x$ to $\partial \Omega$. Moreover, we define $\partial_T \Omega (\delta) = \left\{ (x,t) \in \Omega_T : -\delta < d_B(x) <0 \right\}$.\\
Define the hypersurface $\Gamma^0$ by
\begin{equation}
\Gamma^0 := \Gamma\,. \label{0thinterface}
\end{equation}
First we define the approximate solutions near the Interface $\Gamma^0$ $(c^K_I,\mu^K_I,\textbf{u}^K_I)$ (inner approximate solutions), near the boundary $\partial \Omega$ $(c^K_\partial,\mu^K_\partial,\textbf{u}^K_\partial)$ (boundary-layer approximate solutions), and away from the Interface $\Gamma^0$ and the boundary $\partial \Omega$ $(c^K_O,\mu^K_O,\textbf{u}^K_O)$ (outer approximate solutions). Then we ``glue'' together these solutions. Let $K$ be a positive integer. Then we define $c^K_I$ by 
\begin{align*}
c^K_I(x,t) & := \sum^K_{i=0}{\epsilon^i \left. c^i(z,x,t) \right|_{z = d^K_\epsilon / \epsilon}} & \forall (x,t) \in \Gamma^0(\delta) \, , \\
\end{align*}
where we define $c^i$, $i=0,\ldots,K$, below. Here
\begin{eqnarray}
d^K_\epsilon(x,t) & := & \sum^K_{i=0}{\epsilon^i d^i(x,t)} \,, \quad \forall (x,t) \in \Gamma^0(\delta),
\end{eqnarray}
for some $d^i: \Gamma^0(\delta) \to \mathbb{R}$, $i=1,\ldots,K$, defined below. Note that we have already determined $d^0$ above. Define $c^K_O$ by
\begin{align*}
c^K_O(x,t) & := \sum^K_{i=0}{\epsilon^i \left( c^+_i(x,t) \chi_{\overline{Q^+_0}} + c^-_i(x,t) \chi_{Q^-_0} \right) } & \forall (x,t) \in \Omega_T \, ,
\end{align*}
where $c^\pm_i$, $i=0,\ldots,K$, are defined below. We use an analogous definition for $(\mu^K_I,\textbf{u}^K_I)$ and $(\mu^K_O,\textbf{u}^K_O)$. Furthermore, we define $(c^K_\partial,\mu^K_\partial,\textbf{u}^K_\partial)$ by
\begin{align*}
c^K_\partial(x,t) & := \sum^K_{i=0}{\epsilon^i \left. c^i_B(z,x,t) \right|_{z = d_B / \epsilon}} - \epsilon^K c^K_B(0,x,t) & \forall (x,t) \in  \overline{ \partial_T \Omega(\delta)} \, , \\
\mu^K_\partial(x,t) & := \sum^K_{i=0}{\epsilon^i \left. \mu^i_B(z,x,t) \right|_{z = d_B / \epsilon}} - \epsilon^K \mu^K_B(0,x,t) & \forall (x,t) \in  \overline{ \partial_T \Omega(\delta)} \, , \\
\textbf{u}^K_\partial(x,t) & := \sum^K_{i=0}{\epsilon^i \left. \textbf{u}^i_B(z,x,t) \right|_{z = d_B / \epsilon}} & \forall (x,t) \in  \overline{ \partial_T \Omega(\delta)} \, ,
\end{align*}
where we define $(c^i_B, \mu^i_B, \textbf{u}^i_B)$, $i=0,\ldots,K$, below. The next step is to combine these functions. Let $\zeta \in C^\infty_0(\mathbb{R})$ be a smooth cut-off function such that
\begin{align*}
\zeta(z) & = 1 \mbox{ if } \left| z \right| < \frac{1}{2} , & \zeta(z) & = 0 \mbox{ if } \left| z \right| > 1 , & z \zeta'(z) \leq 0 \mbox{ in } \mathbb{R} \,. 
\end{align*}
Then we define
\begin{eqnarray}
c^K_A := \left\{ \begin{array}{l@{\quad}l}
c^K_\partial & \mbox{in } \overline{\partial_T \Omega(\delta/2)} \, , \\
c^K_\partial \zeta(d_B/\delta) + c^K_O (1 - \zeta(d_B / \delta)) & \mbox{in } \partial_T \Omega(\delta) \backslash \overline{ \partial_T \Omega(\delta/2)} \, , \\
c^K_O & \mbox{in }\Omega_T \backslash (\partial_T \Omega(\delta) \cup \Gamma^0(\delta)) \, , \\
c^K_I \zeta(d^0 / \delta) + c^K_O ( 1 - \zeta(d^0 / \delta)) & \mbox{in } \Gamma^0(\delta) \backslash \Gamma^0(\delta/2) \, , \\
c^K_I & \mbox{in } \Gamma^0(\delta/2) \, .
\end{array} \right. \label{cKA}
\end{eqnarray}
We similarly define $\mu^K_A$ and $\textbf{u}^K_A$. It remains to define the inner expansion functions, $(c^i,\mu^i, \textbf{u}^i)$, the outer expansion functions $(c^\pm_i,\mu^\pm_i, \textbf{u}^\pm_i)$, and the boundary-layer expansion functions $(c^i_B,\mu^i_B, \textbf{u}^i_B)$ for $i=0,\ldots,K$, and $d^i$ for $i=1,\ldots,K$ recursively. Additionally, let $\eta(z) \in C^\infty(\mathbb{R})$ be an arbitrary fixed function satisfying
\begin{eqnarray}
&& \eta(z)= \left\{
\begin{array}{r@ {\quad \mbox{if }} l}
0 & z \leq -1 \\
1 & z \geq 1
 \end{array}\right. , \: \eta'(z) \geq 0 \quad \forall z \in \mathbb{R} ,  \label{eta}\\
& \mbox{and}& \int_{\mathbb{R}}{(\eta(z) - \tfrac{1}{2}) \theta'_0(z) \, dz} = \int_\mathbb{R}{z \eta'(z)\theta'_0(z) \, dz} = 0, \label{eigeneta}
\end{eqnarray}
where $\theta_0$ is the unique solution to (\ref{theta0}). 

\begin{remark}
If $\theta'_0$ is axisymmetric, then we choose $\eta$ such that $\eta - 1/2$ is point symmetric and $\eta'$ is axisymmetric. So both equalities of (\ref{eigeneta}) are fulfilled. This holds for example for $f(c)= c^3-c$ or more generally $f(-c) = f(c)$ for all $c \in \mathbb{R}$. Then $\theta'_0$ is an even function.
\end{remark}

For $i=0$ we define the outer expansion functions $(c^\pm_0, \mu^\pm_0, \textbf{u}^\pm_0)$ in $Q^\pm_0 \cup \Gamma^0(\delta)$ by
\begin{align}
c^\pm_0(x,t) & = \pm 1, & \mu^\pm_0(x,t) & = \mu^\pm(x,t), & \textbf{u}^\pm_0(x,t) & = \textbf{u}^\pm(x,t) , \label{0thouter}
\end{align}
where $\mu^\pm = \left. \mu \right|_{\Omega^\pm},  \textbf{u}^\pm = \left. \textbf{u} \right|_{\Omega^\pm}$ and $(\mu,\textbf{u})$ is the solution to (\ref{sharpsystem1})-(\ref{sharpsystem7}).
The inner expansions functions $(c^0, \mu^0, \textbf{u}^0)$ are defined in $\mathbb{R} \times \Gamma^0(\delta)$ by
\begin{equation}
\begin{aligned}
c^0(z,x,t) & = \theta_0(z), \quad \mu^0(z,x,t) = \mu^+(x,t) \eta(z) + \mu^-(x,t) (1 - \eta(z)),\\
\textbf{u}^0(z,x,t) & = \textbf{u}^+(x,t)  \eta(z) + \textbf{u}^-(x,t) (1 - \eta(z)),
\end{aligned}
\end{equation}
and the boundary expansions functions $(c^0_B, \mu^0_B, \textbf{u}^0_B)$ in $(-\infty,0] \times \overline{\partial_T \Omega(\delta)}$ by
\begin{align}
c^0_B(z,x,t) & = 1 , & \mu^0_B(z,x,t) & = \mu^+(x,t), & \textbf{u}^0_B(z,x,t) & = \textbf{u}^+(x,t). \label{0thboundary}
\end{align}
Note that we have already defined $d^0$ as the signed distance function from $x$ to $\Gamma^0_t$.\\
Now let $j \in \{1,\ldots,K\}$ be any integer.\\
Next we assume that $(c^i, \mu^i, \textbf{u}^i,c^\pm_i, \mu^\pm_i, \textbf{u}^\pm_i,c^i_B, \mu^i_B, \textbf{u}^i_B)$ are defined for every $i=0,\ldots,j-1$. By Taylor series expansion we obtain functions $f^i$, $i = 1, \ldots,K$, such that 
\[
f(c_0 + \epsilon c_1 + \ldots + \epsilon^K c_K) = f(c_0) + f'(c_0) \sum^K_{i=1}{\epsilon^i c_i } + \sum^{K+1}_{i=2}{\epsilon^i f^{i-1}(c_0, \ldots , c_{i-1})} \, 
\]
for all $c_1, \ldots, c_K \in \mathbb{R}$. We define for $j \geq 1$
\begin{equation}
c^\pm_j(x,t) := \frac{\mu^\pm_{k-1} - f^{k-1}(c_0^\pm, \ldots , c_{k-1}^\pm) + \mathcal{E}^\star : \mathcal{C} \left(\nabla\textbf{u}^\pm_{k-1} - \mathcal{E}^\star c^\pm_{k-1} \right) + \Delta c^\pm_{k-2}}{f'(c_0^\pm)} \label{defcpmj}
\end{equation}
for all $(x,t) \in Q^\pm_0 \cup \Gamma^0(\delta)$ and where $c^\pm_{-1}=0$. \\
For $(z,x,t) \in \mathbb{R} \times \Gamma^0(\delta)$ we define $c^j(z,x,t)$ as the unique solution to the following system 
\begin{equation}
\begin{aligned}
& c^j_{zz}(z,x,t) - f'(\theta_0(z))\, c^j(z,x,t) = a^1_{j-1}(z,x,t) \quad \forall z \in \mathbb{R},\\
& c^j(0,x,t)=0, \quad c^j(.,x,t) \in L^\infty(\mathbb{R}) \,, \label{defsystemck}
\end{aligned}
\end{equation}
where $a^1_{j-1}(z,x,t)$ is a given function depending only on functions of order lower than $j-1$. The existence of a unique solution to system (\ref{defsystemck}) is proved in Subsection \ref{subsecasyana} where we define the exact form of $a^1_{j-1}$.\\
For $(z,x,t) \in (-\infty,0] \times \overline{\partial_T\Omega(\delta)}$ the function $c^j_B(z,x,t)$ is defined as the unique solution to the system
\begin{equation}
\begin{aligned}
& c^j_{B,zz}(z,x,t) - f'(1)\, c^j_B(z,x,t) = B^{j-1}_B(z,x,t) \quad \forall z \in (- \infty,0],\\
& c^j_{B,z}(0,x,t) = - \nabla d_B \cdot \nabla \mu^{j-1}_B(0,S_B(x),t), \quad c^j_B(.,x,t) \in L^\infty((-\infty,0]) \,, \label{systemckB}
\end{aligned}
\end{equation}
where 
\begin{eqnarray}
B^{j-1}_B(z,x,t) & = & - \mathcal{E}^\star : \mathcal{C} (\textbf{u}^j_{B,z} \otimes \nabla d_B) - \mathcal{E}^\star : \mathcal{C}(\nabla \textbf{u}^{j-1}_B - \mathcal{E}^\star c^{j-1}_B) - \mu^{j-1}_B \nonumber \\
&&  + f^{j-1}(c^0_B, \ldots, c^{j-1}_B) - 2 \nabla d_B \cdot \nabla c^{j-1}_{B,z} - \Delta d_B c^{j-1}_{B,z} - \Delta c^{j-2}_B  , \:\:\:\: \label{B_B}
\end{eqnarray}
where we have defined $\textbf{u}^{-2}_B = \textbf{u}^{-1}_B = c_B^{-2} = c^{-1}_B = \mu^{-2}_B = \mu^{-1}_B = 0$ and $c^0_B = 1$. The existence of a unique solution to system (\ref{systemckB}) is proved in Subsection \ref{subsecasyana} below. 

\begin{remark}
Due to equation (\ref{DefuB}) below, we will see that $B^{k-1}_B$ is independent of the functions of order $k$, that is independent of $(c^k, \mu^k, \textbf{u}^k,c^\pm_k, \mu^\pm_k, \textbf{u}^\pm_k,c^k_B, \mu^k_B, \textbf{u}^k_B)$. More precisely, we get the identity $\textbf{u}^k_{B,z} = \int^z_{-\infty}{A^{k-1}_B(y)\, dy}$. So it is possible to construct $c^k_B$ by functions of order lower than $k$.
\end{remark}

So we can assume that $(c^j,c^\pm_j,c^j_B)$ are known functions which only depend on functions of order lower than $j-1$.\\
We define $(\mu^\pm_j, \textbf{u}^\pm_j, d^j)$ as the unique solution to the linearized Hele-Shaw problem 
\begin{align}
\Delta \mu^\pm_j & = \partial_t c^\pm_j & \mbox{in }& Q^\pm_0, \label{linHele1}\\
\mu^\pm_j & =  - \sigma \Delta d^j - \tfrac{1}{2} \mathcal{E}^\star : \mathcal{C} \left( \nabla \textbf{u}^+_j + \nabla \textbf{u}^-_j \right) \nonumber \\
& \quad \,- a^{2 \pm}_{j-1} d^j + a^{3 \pm}_{j-1} & \mbox{on } & \Gamma^0,\\
\tfrac{\partial }{\partial n} \mu^+_j & = G^{j-1} & \mbox{on }& \partial_T \Omega, \label{linHele3}\\
\Div \left( \mathcal{C} \mathcal{E}(\textbf{u}^\pm_j) \right) &= \Div \left( \mathcal{C} \mathcal{E}^\star c^\pm_j \right)  & \mbox{in }& Q^\pm_0 , \label{linHele4} \\
\left( \mathcal{C} \left[ \nabla \textbf{u}^\pm_j \right]_{\Gamma^0_t} \right) \nu_{\Gamma^0_t} &= a^{4}_{j-1} \nabla d^j + a^{5}_{j-1} d^j + a^{6}_{j-1} & \mbox{on }& \Gamma^0, \\
\left[ \textbf{u}^\pm_j \right]_{\Gamma^0_t} &= a^{7}_{j-1} d^j + a^{8}_{j-1} & \mbox{on } & \Gamma^0, \\
\textbf{u}^+_j & = \textbf{F}^{j-1} & \mbox{on }& \partial_T \Omega, \label{linHele7}\\
\partial_t d^j &= a^{9}_{j-1} d^j + \tfrac{1}{2} \left[\tfrac{\partial }{\partial \nu} \mu^\pm_j \right] + a^{10}_{j-1} & \mbox{on }& \Gamma^0, \label{linHele8}\\
\nabla d^0 \cdot \nabla d^j &= - \tfrac{1}{2} \sum\nolimits^{j-1}_{i=1}{\nabla d^i \cdot \nabla d^{j-i}} & \mbox{in }& \Gamma^0(\delta), \label{linHele9}\\
d^j(x,0) &= 0 & \mbox{in }& \Gamma_0, \label{linHele10}
\end{align}
where $a^{i}_{j-1}$, $i=2,\ldots,10$, only depends on known functions of order less than or equal to $j-1$. Here we define $\textbf{F}^{-1} = G^{-2} = G^{-1} = 0$ and
\begin{eqnarray}
\textbf{F}^{j-1}(x,t) & = & - \int^0_{-\infty}{\int^z_{-\infty}{\textbf{A}^{j-1}_B(w,x,t) \, dw \,} dz} \,, \label{randbedu}\\
G^{j-1}(x,t) & = & (\Delta d_B(x) + \nabla d_B(x) \cdot \nabla) \int^0_{-\infty}{\int^z_{-\infty}{C^{j-1}_B(w,x,t) \, dw \, dz}} \nonumber \\
&& + \int^0_{-\infty}{(\Delta \mu^{j-1}_B - c^{j-1}_{B,t})(z,x,t) \, dz}  \label{randbedmu}
\end{eqnarray}
for all $j\geq 1$ and for all $(x,t) \in \overline{\partial_T\Omega(\delta)}$ and where
\begin{eqnarray}
\textbf{A}^{j-1}_B(z,x,t) & = & M_B^{-1} \left[ - (\mathcal{C}_{iji'j'} \partial_j (\textbf{u}^{j-1}_{B,i'})_z \partial_{j'} d_B)_{i=1,\ldots,d} -(\mathcal{C}_{iji'j'} (\textbf{u}^{j-1}_{B,i'})_z \partial_{jj'} d_B)_{i=1,\ldots,d} \right. \nonumber \\
&& - (\mathcal{C} : \nabla \textbf{u}^{j-1}_{B,z}) \nabla d_B - \Div \left(\mathcal{C} : \nabla \textbf{u}^{j-2}_{B}\right) \nonumber \\
&& \left. + c^{j-1}_{B,z} (\mathcal{C}:\mathcal{E}^\star) \nabla d_B + (\mathcal{C}:\mathcal{E}^\star) \nabla c^{j-2}_B \right], \label{A_B}\\
C^{j-1}_B(z,x,t) & = & -\Delta d_B \mu^{j-1}_{B,z} - 2 \nabla d_B \cdot \nabla \mu^{j-1}_{B,z} - \Delta \mu^{j-2}_B + c^{j-2}_{B,t} \label{C_B}
\end{eqnarray}
for all $(z,x,t) \in (-\infty,0] \times \overline{\partial_T\Omega(\delta)}$ and where $M_B(x) = \left(\mathcal{C}_{iji'j'} \, \partial_j d_B \, \partial_{j'} d_B \right)^d_{i,i'=1} $. The invertibility of $M_B$ can be shown by standard arguments.\\
In $\mathbb{R} \times \Gamma^0(\delta)$ the inner expansion functions $\mu^j$ and $\textbf{u}^j$ are defined by
\begin{align}
 \mu^j(z,x,t) & := \mu^+_j(x,t) \eta(z) + \mu^-_j(x,t) (1 - \eta(z)) + a^{11}_{j-1}(z,x,t), \label{Defmuj}\\
\textbf{u}^j(z,x,t) & := \textbf{u}^+_j(x,t)  \eta(z) + \textbf{u}^-_j(x,t) (1 - \eta(z)) + a^{12}_{j-1}(z,x,t), \label{Defuj}
\end{align}
where $a^{11}_{j-1}$ and $a^{12}_{j-1}$ are given function depending only on functions of order less than or equal to $j-1$. The exact definition of $a^i_{j-1}$, $i=11,12$, is given in Subsection \ref{subsecjth} below.\\
For $(z,x,t) \in (-\infty,0] \times \overline{\partial_T \Omega(\delta)}$ we define
\begin{eqnarray}
\mu^j_B(z,x,t) & := & \mu^+_j(x,t) +\int^z_{-\infty}{\int^y_{-\infty}{C^{j-1}_B(w,x,t)\, dw\,}dy}, \label{DefcB} \\
\textbf{u}^j_B(z,x,t) & := & \textbf{u}^+_j(x,t) + \int^z_{-\infty}{\int^y_{-\infty}{\textbf{A}^{j-1}_B(w,x,t)\, dw\,}dy}\,. \label{DefuB}
\end{eqnarray}
By a suitable choice of $a^i_{j-1}$, $i=1,\ldots,12$ (see Subsection \ref{subsecasyana} below), we can verify that the family $\left\lbrace c^K_A, \mu^K_A, \textbf{u}^K_A\right\rbrace_{0< \epsilon \leq 1} $ satisfies the required properties of Theorem \ref{theorem1} for $K>0$ large enough. In the next Subsection we define $a^i_{j-1}$, $i=1,\ldots,12$, and motivate the definition of $\left\lbrace c^K_A, \mu^K_A, \textbf{u}^K_A\right\rbrace_{0< \epsilon \leq 1} $.\\
Moreover, define 
\begin{eqnarray}
\Gamma^K_\epsilon  := \left\{(x,t) \in \Gamma^0(\delta) : d^K_\epsilon(x,t) = 0\right\} \label{GammaKepsilon} .
\end{eqnarray}
Because of (\ref{linHele9}), we obtain that $d^K_\epsilon$ is a $K$-th order approximate distance function, i.e. $d^K_\epsilon$ vanishes on $\Gamma^K_\epsilon$ and 
\begin{eqnarray}
\left|\nabla d^K_\epsilon \right|^2 = 1 + \sum_{\substack{1 \leq i,j \leq K \\ i+j \geq K+1}}{ \epsilon^{i+j} \nabla d^j \cdot \nabla d^i } = 1 + \mathcal{O}(\epsilon^{K+1}) \quad \forall (x,t) \in \Gamma^0(\delta) \, . \label{nabladK}
\end{eqnarray}
Here and in the following the Landau symbol $\mathcal{O}$ is with respect to the $C^0$-norm.\\
After a small modification of the approximate solution (see Subsection \ref{subsecproof}), we will prove the following result: 

\begin{theorem} \label{theoremcA-c0}
Let $\Gamma_{00} \subset \Omega$ be a given smooth hypersurface without boundary and let $(\mu_0,\textbf{u}_0,\Gamma^0)$ be a smooth solution to the Hele-Shaw problem (\ref{sharpsystem1})-(\ref{sharpsystem7}) for $t \in [0,T]$ with initial value $\Gamma_{00}$ such that $\Gamma^0 \subset \Omega \times [0,T]$. Then for every $K>3$, there exists a positive constant $\epsilon_0$ such that for every $\epsilon \in (0,\epsilon_0]$ there exists an approximate solution $(c^\epsilon_A,\mu^\epsilon_A,\textbf{u}^\epsilon_A)$ satisfying 
\begin{align}
(c^\epsilon_A)_t - \Delta \mu^\epsilon_A & = 0 && \mbox{in } \Omega_T\, , \label{apprsol1} \\
\mu^\epsilon_A + \epsilon \Delta c^\epsilon_A - \epsilon^{-1} f(c^\epsilon_A) - W_{,c}(c^\epsilon_A,\mathcal{E}(\textbf{u}^\epsilon_A)) & = \mathcal{O}(\epsilon^{K-3}) && \mbox{in } \Omega_T \, ,\\
\Div \left( \mathcal{C} \mathcal{E}(\textbf{u}^\epsilon_A) \right) - \Div \left(\mathcal{C} \mathcal{E}^\star c^\epsilon_A\right) & = \mathcal{O}(\epsilon^{K-2}) && \mbox{in } \Omega_T \, \label{apprsol3} ,\\
\tfrac{\partial}{\partial n} c^\epsilon_A = \tfrac{\partial}{\partial n} \mu^\epsilon_A & = 0 && \mbox{on } \partial_T \Omega \, , \\
\textbf{u}^\epsilon_A & = 0 && \mbox{on } \partial_T \Omega  \, . \label{apprsol5}
\end{align}
Additionally, it holds as $\epsilon \searrow 0$
\begin{eqnarray*}
\left\| \mu_A^\epsilon - \mu_0 \right\|_{C^0(\Omega_T)} & = & \mathcal{O}(\epsilon) \, , \\
\left\| c^\epsilon_A(x,t) - \theta_0(d^0(x,t)/\epsilon + d^1(x,t)) \right\|_{C^0(\Gamma^0(\delta))} & = & \mathcal{O}(\epsilon) \, , \\
\left\| c^\epsilon_A \mp 1 \right\|_{C^0(Q^\pm_0 \backslash \Gamma^0(\delta/2))} & = & \mathcal{O}(\epsilon) \, ,\\
\left\| \textbf{u}_A^\epsilon - \textbf{u}_0 \right\|_{C^0(\Omega_T)} & = & \mathcal{O}(\epsilon) \, . 
%\left\| \nabla (\textbf{u}^\epsilon_A - \textbf{u}_0) \right\|_{C^0(Q^\pm_0 \backslash \Gamma^0(\delta/2))} & = & \mathcal{O}(\epsilon) \, ,\textcolor{red}{nicht notwendig} \\
%\left\| \nabla (\textbf{u}^\epsilon_A(.,t) - \textbf{u}_0(.,t) \right\|_{L^p(\Omega)} & = & \mathcal{O}(\epsilon^{\frac{1}{p}}) \,, \textcolor{red}{nicht notwendig}
\end{eqnarray*}
\end{theorem}

\subsection{Formal Derivation of the Expansion} \label{subsecasyana}

Away from the interface $\Gamma^\epsilon = \left\lbrace(x,t) : c^\epsilon(x,t)=0 \right\rbrace $ we use the original variables to determine the expansion of the solutions $(c^\epsilon,\mu^\epsilon, \textbf{u}^\epsilon)$. This is called the outer expansion. Near the interface $\Gamma^\epsilon$ we expect that $\nabla d^\epsilon \cdot \nabla c^\epsilon \approx \frac{C}{\epsilon}$ for some constant $C>0$ and where $d^\epsilon$ is the spatial signed distance function to the interface. Therefore we introduce the new variable $z = \frac{d^\epsilon(x,t)}{\epsilon}$ to describe the sharp change near the interface. This is called the inner expansion. We also use a boundary-layer expansion to satisfy the boundary conditions. By the so-called matching conditions we connect the inner and outer expansion  and the outer and boundary-layer expansion to obtain suitable approximate solutions $(c^\epsilon_A,\mu^\epsilon_A,\textbf{u}^\epsilon_A)$ for all $\epsilon \in (0,1]$.\\
We use the same convention as in \cite{ABC}, more precisely that means the notation $\sum^\infty_{i=0}$ is formal and should be understood only as a finite sum $\sum^K_{i=0}$ plus an error term of order $\mathcal{O}(\epsilon^{K+1})$ where $K$ is a large integer depending on the order of approximation needed. Also the word ``smooth'' should be understood to mean that all the needed derivatives exist and are continuous and bounded. Similarly, the phrase ``for all natural integers'' should be understood to mean ``for all integers needed''.\\
In the following we use the symmetry assumption $\mathcal{C}_{iji'j'} = \mathcal{C}_{ijj'i'}$ for all $i,j,j',i' \in \left\{1,\ldots,d\right\}$. Therefore we obtain $\mathcal{C} A = \mathcal{C} A^T$ for all $A \in \mathbb{R}^{d \times d}$ and in particular, $\mathcal{C} \mathcal{E}(\textbf{u}) = \mathcal{C} \nabla \textbf{u}$. For the sake of clarity in the asymptotic expansion we use $\mathcal{C} \nabla \textbf{u}$ instead of $\mathcal{C} \mathcal{E}(\textbf{u})$. 

\begin{remark}
In this Subsection we motivate the definition of $\left\{ c^K_A, \mu^K_A, \textbf{u}^K_A \right\}_{0 < \epsilon \leq 1}$ in Subsection \ref{secapp}. Therefore we use for some functions in this Subsection the same expressions as in Subsection \ref{secapp}. We will see that these functions coincide.
\end{remark}

\subsubsection{Representation of the Interface} \label{repinterface}

Let $(c^\epsilon,\mu^\epsilon,\textbf{u}^\epsilon)$ be the solution to (\ref{system1})-(\ref{system6}). For the formal derivation we assume that 
\[
\Gamma^\epsilon : = \left\{ (x,t) \in \Omega_T : c^\epsilon (x,t) = 0 \right\} = \bigcup_{0 < t < T} \left( \Gamma^\epsilon_t \times \left\{t\right\} \right)
\]
is a smooth hypersurface. $\Gamma^\epsilon$ is called the interface. Let $Q^-_\epsilon$ be the interior of $\Gamma^\epsilon$ and $Q^+_\epsilon : = \Omega_T \backslash \left( \Gamma^\epsilon \cup Q^-_\epsilon \right)$. Furthermore, let $d^\epsilon(x,t)$ be the spatial signed distance function to $\Gamma^\epsilon_t$ such that $d^\epsilon < 0$ in $Q^-_\epsilon$. Then $d^\epsilon$ is a smooth function and $\left| \nabla d^\epsilon \right|=1$ in a neighborhood of $\Gamma^\epsilon$, which depends on the curvature of $\Gamma^\epsilon$. Also we assume that $d^\epsilon$ has the expansion 
\[
d^\epsilon(x,t) = \sum^\infty_{i=0}{\epsilon^i d^i(x,t)} \,,
\]
where $d^0$ is defined in $\overline{\Omega_T}$ and $d^i$, $i \geq 1$, is defined in a neighborhood of $\Gamma^\epsilon$. Since $d^i$ is independent of $\epsilon$ for all $i \geq 0$, the equation $\left| \nabla d^\epsilon \right| =1 $ is equivalent to 
\begin{eqnarray}
\nabla d^0 \cdot \nabla d^k = \left\{ \begin{array}{l@{\quad}l}
1 & \mbox{if }k=0\,, \\
0 & \mbox{if }k=1\,, \\
- \tfrac{1}{2} \sum^{k-1}_{i=1}{ \nabla d^i \cdot \nabla d^{k-i}} & \mbox{if }k \geq 2\,,
\end{array} \right. \label{djbedingung}
\end{eqnarray}
where all the equations are satisfied in a neighborhood of $\Gamma^\epsilon$.  Note that this equation coincides with (\ref{linHele9}) and the definition of $d^0$ in Subsection \ref{secapp}. Then we can assume that $d^0$ is a spatial signed distance function, and we define
\begin{eqnarray}
\Gamma^0 & = & \left\{ (x,t) \in \Omega_T : d^0(x,t) = 0 \right\} \,,\label{defgamma0} \\
\Gamma^0(\delta) & = & \left\{ (x,t) \in \Omega_T : \left| d^0(x,t) \right| < \delta \right\} \,, \\
Q^\pm_0 & = & \left\{ (x,t) \in \Omega_T :  \pm d^0(x,t) > 0 \right\}  \label{defQ-}
\end{eqnarray}
for some constant $\delta >0$. 

\begin{remark}
Only to motivate the construction of the approximate solutions, we need the assumptions that $\Gamma^\epsilon$ is a smooth hypersurface and $d^\epsilon$ has a series expansion. But these assumptions are not necessary to show the convergence of the Cahn-Larché system (\ref{system1})-(\ref{system6}) to the modified Hele-Shaw problem (\ref{sharpsystem1})-(\ref{sharpsystem7}). More precisely, in Theorem \ref{hoeregR}, \ref{theoremcA-c0}, and \ref{mainresult} these assumptions do not occur. Moreover, the assumption that for all $\epsilon \in (0,1]$ the solution $(c^\epsilon,\mu^\epsilon,\textbf{u}^\epsilon)$ has a series expansion near the interface $\Gamma^\epsilon$, near the boundary $\partial \Omega$, and away from the interface $\Gamma^\epsilon$ and the boundary $\partial \Omega$ (see Subsection \ref{secouterexp}, \ref{secinnerexp}, and \ref{boundaryexp}) is also not necessary for the proof of the convergence.
\end{remark}

\subsubsection{Outer Expansion}\label{secouterexp}

We assume that $\Gamma^0$ is known. Then so are $d^0$, $Q^+_0$ and $Q^-_0$. Also we assume that away from the interface $\Gamma^\epsilon$ the  solution function $c^\epsilon$ has the expansion
\begin{align*}
c^\epsilon(x,t) & = c_0^\pm(x,t) + \epsilon c^\pm_1(x,t) + \epsilon^2 c^\pm_2(x,t)+... & \mbox{in } & Q^\pm_0 \backslash \Gamma^0 \! \left(\tfrac{\delta}{2}\right), 
\end{align*}
where $c^\pm_i$ are appropriate functions defined in $Q^\pm_0$ and $\delta > 0$ is a fixed constant independent of $\epsilon$ which is to be determined later. We assume an analogous expansion for $\mu^\epsilon$ and $\textbf{u}^\epsilon$.\\
We substitute the expansion  for $(c^\epsilon,\mu^\epsilon, \textbf{u}^\epsilon)$ into (\ref{system1})-(\ref{system3b}) and match the terms with the same power of $\epsilon$. Then for all $(x,t) \in Q^\pm_0$, we obtain $c^\pm_0=\pm 1$ and (\ref{defcpmj}) for $k \geq 1$, $\Delta \mu^\pm_0=0$ and (\ref{linHele1}) for $k \geq 1$, and $\Div \left( \mathcal{C} \nabla \textbf{u}^\pm_0\right) =0 $ and (\ref{linHele4}) for $k \geq 1$. These are the so-called outer expansion equations.

\begin{remark} \label{remarkoutexp}
For the outer expansion let us mention the following points.
\begin{enumerate}
	\item When $k=0$ we only obtain $\Delta \mu^\pm_0=\partial_t c^\pm_0$, $f(c^\pm_0) = 0$, and $\Div \left( \mathcal{C} \nabla \textbf{u}^\pm_0\right) = \Div \left( \mathcal{C} \mathcal{E}^\star c^\pm_0 \right) $. But as in \cite{ABC} we require $c^\pm_0 = \pm 1$ for definiteness. 
	\item In order to construct an approximate solution we have required that all outer expansion equations be satisfied in $Q^\pm_0$ instead of $Q^\pm_0 \backslash \Gamma^0(\delta/2)$.
	\item Because of the definition of $c^\pm_k$ we can not require boundary conditions for $c^+_k$ on $\partial_T \Omega$. To avoid this problem, we use a boundary-layer expansion. For details see Subsection \ref{boundaryexp}.
	\item To determine $\mu^\pm_k$ and $\textbf{u}^\pm_k$ uniquely, we need boundary conditions on $\partial \Omega$ and $\Gamma^0$. We obtain these conditions by the boundary-layer expansion and the inner expansion, see Subsection \ref{subsec0th} and \ref{subsecjth}. 
	\item For the inner expansion it is necessary to define $(c^\pm_k,\mu^\pm_k,\textbf{u}^\pm_k)$ not only in the domain $Q^\pm_0$, but also in $\Gamma^0(\delta) \backslash Q^\pm_0$. But it is sufficient to choose any $C^{\bar{K}}$ extension of $(c^\pm_k,\mu^\pm,\textbf{u}^\pm_k)$ from $Q^\pm_0$ to $Q^\pm_0 \cup \Gamma^0(\delta)$ where $\bar{K}$ is an integer depending on the order of approximations needed. For the extension of $c^\pm_k$ we can use equation (\ref{defcpmj}). We can extend $\mu^\pm_k$ and $\textbf{u}^\pm_k$ as in \cite[Remark 4.1]{ABC}.
\end{enumerate}
\end{remark}

In order to have bounded solutions in the inner expansion, we need the following definitions
\begin{align}
O^\pm_k(x,t) & :=(c^\pm_k)_t - \Delta \mu^\pm_k, & O^\pm & := \sum^\infty_{i=0}{\epsilon^i O^\pm_i} & \mbox{in }Q^\pm_0\cup\Gamma^0(\delta)\,, \label{Opm} \\
\textbf{P}^\pm_k(x,t) & := \Div \left( \mathcal{C} \nabla \textbf{u}^\pm_k \right) - \left( \mathcal{C} \mathcal{E}^\star \right) \nabla c^\pm_k\,, & \textbf{P}^\pm & := \sum^\infty_{i=0}{\epsilon^i \textbf{P}^\pm_i} & \mbox{in }Q^\pm_0\cup\Gamma^0(\delta)\,. \label{Ppm}
\end{align}
Due to the definition of $(c^\pm_k,\mu^\pm_k,\textbf{u}^\pm_k)$, it holds that $\textbf{P}^\pm_k = 0$ and $O^\pm_k = 0$ in $\overline{Q^\pm_0}$ for all $k \geq 0$.\\

\subsubsection{Inner Expansion} \label{secinnerexp}

As above we assume that $d^0(x,t)$ is known. To understand the behavior of the solution $(c^\epsilon,\mu^\epsilon,\textbf{u}^\epsilon)$ near the interface $\Gamma^0_t$, we assume that in $\Gamma^0(\delta)$ the solution $(c^\epsilon,\mu^\epsilon, \textbf{u}^\epsilon)$ has the expansion 
\begin{align*}
c^\epsilon(x,t) & = \tilde{c}^\epsilon\left(\frac{d^\epsilon(x,t)}{\epsilon}, x,t\right),& \tilde{c}^\epsilon(z,x,t) & = \sum^\infty_{i=0}{\epsilon^i c^i(z,x,t)},
\end{align*}
where $\tilde{c}^\epsilon$ and $c^i$ are appropriate functions defined in $\mathbb{R} \times \Gamma^0(\delta)$. We assume an analogous expansion for $\mu^\epsilon$ and $\textbf{u}^\epsilon$.\\
To construct a solution in the whole domain $\Omega_T$, we require some matching conditions for the inner and outer expansion. Since $ \frac{d^0(x,t)}{\epsilon}  \to  \infty$ as $\epsilon \to 0$ in $Q^+_0$, we require the following inner-outer matching conditions as $ z \to \infty$
\begin{eqnarray}
D_x^m D_t^n D_z^l \left[(c^k,\mu^k,\textbf{u}^k)(\pm z,x,t) - (c^\pm_k,\mu^\pm_k,\textbf{u}^\pm_k)(x,t)\right] & = & \mathcal{O}(e^{-\alpha z}) \label{match1}
\end{eqnarray}
for all $(x,t) \in \Gamma^0(\delta)$ and all $k,m,n,l \in \left\lbrace  0, \ldots, \bar{K}\right\rbrace $ where $\bar{K}>0$ depends of the order of expansion. Here and in the following $\alpha>0$ is the same constant as in Lemma \ref{lemmatheta0}.

\begin{remark}
In Subsection \ref{secapp} we glue together the outer and inner expansions. Then it is necessary that the matching conditions holds for $m,n,l\in \left\{0,1,2\right\}$ for each order $k$. Since the equations for $(c^\pm_k, c^k, \mu^\pm_k, \mu^k, \textbf{u}^\pm_k, \textbf{u}^k)$ depend on space and time derivatives and derivatives with respect to $z$ of functions of lower order, it is necessary and sufficient that $ m,n,l\in \left\{0,\ldots,\bar{K}\right\}$ where $\bar{K}$ is a constant  depending on the order of expansion. To verify the matching conditions we consider the inner expansion equations for $(c^j,\mu^j,\textbf{u}^j)$, which we obtain below, and then we use the results of Subsection \ref{seccomp}. One can even verify that the matching conditions are true for all $m,n,l \in \mathbb{N}$.
\end{remark}

Since $c^\epsilon = 0$ on $\Gamma^\epsilon$, it is natural that we require in our construction 
\[
c^k(0,x,t) = 0 \quad \forall (x,t) \in \Gamma^0(\delta)
\]
for all $k \geq 0$.\\
As above we substitute the expansion of $(c^\epsilon,\mu^\epsilon, \textbf{u}^\epsilon)$ into (\ref{system1}), which yields
\begin{eqnarray*}
 \epsilon^{-1} \tilde{c}^\epsilon_z d^\epsilon_t + \tilde{c}^\epsilon_t = \epsilon^{-2} \tilde{\mu}^\epsilon_{zz} + 2 \epsilon^{-1} \nabla \tilde{\mu}^\epsilon_z \cdot \nabla d^\epsilon + \epsilon^{-1} \tilde{\mu}^\epsilon_z \Delta d^\epsilon + \Delta \tilde{\mu}^\epsilon,
\end{eqnarray*}
into (\ref{system2}), which yields
\begin{eqnarray*}
 \tilde{\mu}^\epsilon & = & \epsilon^{-1} \left( f(\tilde{c}^\epsilon) - \tilde{c}^\epsilon_{zz} \right) - \tilde{c}^\epsilon_{z} \Delta d^\epsilon - 2 \nabla \tilde{c}^\epsilon_z \cdot \nabla d^\epsilon - \epsilon \Delta \tilde{c}^\epsilon \\
 &&- \epsilon^{-1} \mathcal{E}^\star : \mathcal{C} \left(\tilde{\textbf{u}}^\epsilon_z \otimes \nabla d^\epsilon  \right) - \mathcal{E}^\star : \mathcal{C} \left( \nabla \tilde{\textbf{u}}^\epsilon - \mathcal{E}^\star \tilde{c}^\epsilon \right) \,,
\end{eqnarray*}
and into (\ref{system3}), which yields
\begin{align*}
& \epsilon^{-2} \left( \mathcal{C} \left( \tilde{\textbf{u}}_{zz}^\epsilon \otimes \nabla d^\epsilon \right) \right) \nabla d^\epsilon + \epsilon^{-1} \left(\mathcal{C}_{iji'j'} \, \partial_j (\tilde{\textbf{u}}^\epsilon_{i'})_z \partial_{j'} d^\epsilon \right)_{i=1, \ldots , d} \\
& + \epsilon^{-1} \left( \mathcal{C}_{iji'j'} \, (\tilde{\textbf{u}}^\epsilon_{i'})_z \partial_{jj'} d^\epsilon \right)_{i=1, \ldots , d} + \epsilon^{-1} \left( \mathcal{C} \nabla \tilde{\textbf{u}}^\epsilon_z \right) \nabla d^\epsilon + \left( \mathcal{C}_{iji'j'} \, \partial_{jj'} \tilde{\textbf{u}}^\epsilon_{i'} \right)_{i=1, \ldots , d} \\
& =  \epsilon^{-1} \tilde{c}^\epsilon_z \left( \mathcal{C} \mathcal{E}^\star \right) \nabla d^\epsilon + \left( \mathcal{C} \mathcal{E}^\star \right) \nabla \tilde{c}^\epsilon 
\end{align*}
for $(z,x,t) \in S^\epsilon:= \left\{ (z,x,t) \in \mathbb{R} \times \Gamma^0(\delta) : z= d^\epsilon(x,t)/\epsilon \right\}$. We can consider these equations as a system of ordinary differential equations for $(c^i,\mu^i,\textbf{u}^i)$ with independent variable $z \in \mathbb{R}$, whereas $(x,t)$ are considered as parameters. Of course, it is not clear that the solutions to these ordinary differential equations satisfy the inner-outer matching conditions. Note that these equations have to be satisfied only on $S^\epsilon$. So we can add any terms which vanishes on $S^\epsilon$ to enforce the inner-outer matching conditions. We denote these terms by $g^\epsilon$, $k^\epsilon$, $h^\epsilon$, $L^\epsilon$, $\textbf{l}^\epsilon$, $\textbf{j}^\epsilon$ and $\textbf{K}^\epsilon$. We will determine them later. \\
Furthermore, we set
\begin{equation}
\eta^\pm_N(z) = \eta(-N\pm z) , \quad z \in \mathbb{R}, \label{etaN}
\end{equation}
where $\eta$ is defined as in (\ref{eta}), (\ref{eigeneta}) and $N>0$ is a constant to be determined.\\
From now we consider the following modified equations for $(\tilde{c}^\epsilon, \tilde{\mu}^\epsilon, \tilde{\textbf{u}}^\epsilon)$
\begin{align}
\tilde{c}^\epsilon_{zz} - f(\tilde{c}^\epsilon) = & - \mathcal{E}^\star : \mathcal{C} \left( \tilde{\textbf{u}}^\epsilon_z \otimes \nabla d^\epsilon \right) \nonumber \\
& -\epsilon \left(\tilde{\mu}^\epsilon + \Delta d^\epsilon \tilde{c}^\epsilon_z + 2 \nabla d^\epsilon \cdot \nabla \tilde{c}^\epsilon_z + \mathcal{E}^\star : \mathcal{C} \left( \nabla \tilde{\textbf{u}}^\epsilon - \mathcal{E}^\star \tilde{c}^\epsilon \right) \right) \, \nonumber\\ 
& - \epsilon^2 \Delta \tilde{c}^\epsilon + g^\epsilon \eta' \left(d^\epsilon - \epsilon z \right) + k^\epsilon \eta' \left(d^\epsilon - \epsilon z \right), \label{modinnere1} \\
\tilde{\mu}^\epsilon_{zz} ={}& \epsilon \left(\tilde{c}^\epsilon_z d^\epsilon_t - 2 \nabla \tilde{\mu}^\epsilon_z \cdot \nabla d^\epsilon - \tilde{\mu}^\epsilon_z \Delta d^\epsilon \right)+ \epsilon^2 \left(\tilde{c}^\epsilon_t - \Delta \tilde{\mu}^\epsilon \right)\nonumber\\
& + \left(h^\epsilon \eta'' + L^\epsilon \eta'\right)\left(d^\epsilon - \epsilon z\right) - \epsilon^2 \left(O^+ \eta^+_N+O^- \eta^-_N\right) ,\label{modinnere2} \\
\left(\mathcal{C} \left( \tilde{\textbf{u}}^\epsilon_{zz} \otimes \nabla d^\epsilon \right) \right) \nabla d^\epsilon  = & - \epsilon \left(\mathcal{C}_{iji'j'} \, \partial_j (\tilde{\textbf{u}}^\epsilon_{i'})_z \partial_{j'} d^\epsilon \right)_{i=1, \ldots , d} \nonumber \\
& - \epsilon \left( \mathcal{C}_{iji'j'} \, (\tilde{\textbf{u}}^\epsilon_{i'})_z \partial_{jj'} d^\epsilon \right)_{i=1, \ldots , d} - \epsilon \left( \mathcal{C} \nabla \tilde{\textbf{u}}^\epsilon_z \right) \nabla d^\epsilon \nonumber \\
& + \epsilon \tilde{c}^\epsilon_z \left( \mathcal{C} \mathcal{E}^\star \right) \nabla d^\epsilon - \epsilon^2 \left( \mathcal{C}_{iji'j'} \, \partial_{jj'} \tilde{\textbf{u}}^\epsilon_{i'} \right)_{i=1, \ldots , d} \nonumber \\
& + \epsilon^2 \left( \mathcal{C} \mathcal{E}^\star \right) \nabla \tilde{c}^\epsilon + M \left(\textbf{l}^\epsilon \eta'' + \textbf{K}^\epsilon \eta' \right) \left(d^\epsilon - \epsilon z\right) \nonumber \\
& + \textbf{j}^\epsilon \eta'' \left( d^\epsilon - \epsilon z\right) + \epsilon^2 \left( \textbf{P}^+ \eta^+_N + \textbf{P}^- \eta^-_N \right)  \label{modinnere3}
\end{align}
for $z \in \mathbb{R}$ and $(x,t) \in \Gamma^0(\delta)$ and 
\begin{equation}
M(x,t) = \left( \mathcal{C}_{iji'j'} \, \partial_{j'} d^0(x,t) \partial_j d^0(x,t) \right)^d_{i,i'=1}\,. \label{defM}
\end{equation}
We have added the terms $\epsilon^2 ( O^+ \eta^+_N+O^- \eta^-_N)$ and $\epsilon^2 ( \textbf{P}^+ \eta^+_N + \textbf{P}^- \eta^-_N)$ to satisfy the compatibility conditions for $\mu^i$ and $\textbf{u}^i$ where $O^\pm$ and $\textbf{P}^\pm$ are defined as in (\ref{Opm}) and (\ref{Ppm}). We will see more details in Subsection \ref{seccomp}. 

\begin{remark} \label{remarkN}
It remains to fix the constant $N$ such that the terms $\epsilon^2 ( O^+ \eta^+_N+O^- \eta^-_N)$ and $\epsilon^2 ( \textbf{P}^+ \eta^+_N + \textbf{P}^- \eta^-_N)$ do not affect the equations needed for $\tilde{\mu}^\epsilon(\frac{d^\epsilon}{\epsilon},x,t)$ and $\tilde{\textbf{u}}^\epsilon(\frac{d^\epsilon}{\epsilon},x,t)$. To this end we will see in Subsection \ref{subsecjth} that we can determine $d^0$ and $d^1$ independent of $\epsilon^2 ( O^+ \eta^+_N+O^- \eta^-_N)$ and $\epsilon^2 ( \textbf{P}^+ \eta^+_N + \textbf{P}^- \eta^-_N)$. So we can set 
\[
N := \left\|  d^1\right\| _{C^0(\Gamma^0(\delta))} + 2\,.
\]
The so-defined $N$ satisfies the required property. This can be seen as in \cite[Remark 4.2]{ABC}. 
\end{remark}

Moreover, we assume that for $(x,t) \in \Gamma^0(\delta)$ the terms $g^\epsilon$, $L^\epsilon$, $h^\epsilon$, $\textbf{l}^\epsilon$, and $\textbf{K}^\epsilon$ have the expansion
\begin{align*}
g^\epsilon(x,t) &= \sum^\infty_{i=0}{\epsilon^{i+1} g^i(x,t)}, & k^\epsilon(x,t) &= \sum^\infty_{i=0}{\epsilon^{i} k^i(x,t)} ,\\
 L^\epsilon(x,t) &= \sum^\infty_{i=0}{\epsilon^{i+1} L^i(x,t)}, & h^\epsilon(x,t) &= \sum^\infty_{i=0}{\epsilon^{i} h^i(x,t)}, \\
\textbf{l}^\epsilon(x,t) & = \sum^\infty_{i=0}{\epsilon^{i} \textbf{l}^i(x,t)}, & \textbf{j}^\epsilon(x,t) & = \sum^\infty_{i=1}{\epsilon^{i} \textbf{j}^i(x,t)}, \\
 \textbf{K}^\epsilon(x,t) & = \sum^\infty_{i=0}{\epsilon^{i+1} \textbf{K}^i(x,t)} .&
\end{align*}
As for the outer expansion we substitute these expansions and the expansions for $\tilde{c}^\epsilon$, $\tilde{\mu}^\epsilon$, $\tilde{\textbf{u}}^\epsilon$, and $d^\epsilon$ into (\ref{modinnere1})-(\ref{modinnere3}) and match the terms with the same power of $\epsilon$. Note that $\mathcal{C} \left( \textbf{u}^0_{zz} \otimes \nabla d^0 \right) \nabla d^0 = M \textbf{u}^0_{zz}$. We show that $M$ is an invertible matrix. Let $\textbf{v} \neq 0$ be an arbitrary vector, then it holds for all $(x,t) \in \Gamma^0(\delta)$
\begin{eqnarray*}
\textbf{v} \cdot (M \textbf{v}) = \left(\textbf{v} \otimes \nabla d^0 \right) : \mathcal{C} \left(\textbf{v} \otimes \nabla d^0 \right) \geq c_2 \left|\textbf{v} \otimes \nabla d^0 \right|^2 > 0 \,. 
\end{eqnarray*}
So we get the following ordinary differential equations
\begin{eqnarray}
\left.\begin{array}{r@{\,=\,}l@{\:}l}
\left(\textbf{u}^0 - \textbf{l}^0 d^0 \eta \right)_{zz} & 0 & \\
\left( \textbf{u}^k - \left( \textbf{l}^k d^0 + \textbf{l}^0 d^k\right) \eta \right)_{zz} & D^{k-1}(z,x,t), & k\geq 1 
\end{array}\right\} z \in \mathbb{R}, \: (x,t) \in \Gamma^0(\delta) , \label{innere3}\\
\left.\begin{array}{r@{\,=\,}l@{\:}l}
c^0_{zz} - f(c^0) & E^0(z,x,t) & \\
c^k_{zz} - f'(c^0) c^k & (E^k + A^{k-1})(z,x,t), & k\geq 1 
\end{array}\right\} z \in \mathbb{R}, \: (x,t) \in \Gamma^0(\delta) , \label{innere1}\\
\left.\begin{array}{r@{\,=\,}l@{\:}l}
\left(\mu^0 - h^0 d^0 \eta\right)_{zz} & 0 & \\
\left(\mu^k - \left(h^k d^0 + h^0 d^k\right) \eta \right)_{zz} & B^{k-1}(z,x,t), & k \geq 1
\end{array}\right\} z \in \mathbb{R}, \: (x,t) \in \Gamma^0(\delta) , \label{innere2}
\end{eqnarray}
where in $\mathbb{R} \times \Gamma^0(\delta)$ $D^0$, $E^0$, $A^0$, and $B^0$ have the following form
\begin{eqnarray}
D^0 & = & M^{-1} \left[- \left(\mathcal{C} \left( \textbf{u}^0_{zz} \otimes \nabla d^1 \right)\right) \nabla d^0 - \left(\mathcal{C} \left( \textbf{u}^0_{zz} \otimes \nabla d^0 \right)\right) \nabla d^1 \right. \nonumber \\
&& - \left(\mathcal{C}_{iji'j'} \, \partial_j (\textbf{u}^{0}_{i'})_z \partial_{j'} d^0 \right)_{i=1, \ldots , d} - \left( \mathcal{C}_{iji'j'} \, (\textbf{u}^{0}_{i'})_z \partial_{jj'} d^0 \right)_{i=1, \ldots , d} \nonumber \\
&& \left.- \left( \mathcal{C} \nabla \textbf{u}^{0}_z \right) \nabla d^0 + c^{0}_z \left(\mathcal{C} \mathcal{E}^\star \right) \nabla d^0 + \textbf{j}^1 d^0 \eta''\right]  - z \textbf{l}^{0} \eta'' +  \textbf{K}^{0} d^0 \eta' \,,\\
E^0 & = & - \mathcal{E}^\star : \mathcal{C} \left( \textbf{u}^0_z \otimes \nabla d^{0} \right) + k^0 d^0 \eta' \,, \\
A^0 & = & - \mu^0 - \Delta d^0 c^0_z - 2 \nabla d^0 \cdot \nabla c^0_z  - \mathcal{E}^\star : \mathcal{C} \left( \nabla \textbf{u}^0 - \mathcal{E}^\star c^0 \right) \nonumber \\
&& + g^0 d^0 \eta' - z k^0 \eta' \,,\\
B^0 & = & d^{0}_t c^0_z - \Delta d^0 \mu^{0}_z - 2 \nabla d^0 \cdot \nabla \mu^{0}_z - z h^{0} \eta'' + L^{0} d^0 \eta' \,,
\end{eqnarray}
and where $E^k$ for $k \geq 1$ and $D^{k-1}, A^{k-1}, B^{k-1}$ for $k \geq 2$ have the following form 
\begin{eqnarray}
D^{k-1} & = & M^{-1} \left[- \left(\mathcal{C} \left( \textbf{u}^0_{zz} \otimes \nabla d^k \right)\right) \nabla d^0 - \left(\mathcal{C} \left( \textbf{u}^0_{zz} \otimes \nabla d^0 \right)\right) \nabla d^k \right. \nonumber \\
&& - \left(\mathcal{C} \left( \textbf{u}^{k-1}_{zz} \otimes \nabla d^0 \right)\right) \nabla d^1 - \left(\mathcal{C} \left( \textbf{u}^{k-1}_{zz} \otimes \nabla d^1 \right)\right) \nabla d^0 \nonumber \\
&& - \left(\mathcal{C} \left( \textbf{u}^{1}_{zz} \otimes \nabla d^{k-1} \right)\right) \nabla d^0 - \left(\mathcal{C} \left( \textbf{u}^{1}_{zz} \otimes \nabla d^0 \right)\right) \nabla d^{k-1} \nonumber \\
&& - \left(\mathcal{C} \left( \textbf{u}^{0}_{zz} \otimes \nabla d^1 \right)\right) \nabla d^{k-1} - \left(\mathcal{C} \left( \textbf{u}^{0}_{zz} \otimes \nabla d^{k-1} \right)\right) \nabla d^1 \nonumber \\
&& - \left(\mathcal{C}_{iji'j'} \, \partial_j (\textbf{u}^{k-1}_{i'})_z \partial_{j'} d^0 \right)_{i=1, \ldots , d} - \left(\mathcal{C}_{iji'j'} \, \partial_j (\textbf{u}^{0}_{i'})_z \partial_{j'} d^{k-1} \right)_{i=1, \ldots , d} \nonumber \\
&& - \left( \mathcal{C}_{iji'j'} \, (\textbf{u}^{k-1}_{i'})_z \partial_{jj'} d^0 \right)_{i=1, \ldots , d} - \left( \mathcal{C}_{iji'j'} \, (\textbf{u}^{0}_{i'})_z \partial_{jj'} d^{k-1} \right)_{i=1, \ldots , d} \nonumber \\
&& - \left( \mathcal{C} \nabla \textbf{u}^{k-1}_z \right) \nabla d^0 - \left( \mathcal{C} \nabla \textbf{u}^{0}_z \right) \nabla d^{k-1} + c^{k-1}_z \left(\mathcal{C} \mathcal{E}^\star \right) \nabla d^0 \nonumber \\
&& \left. + c^{0}_z \left(\mathcal{C} \mathcal{E}^\star \right) \nabla d^{k-1} + \left( \textbf{j}^k d^0 + \textbf{j}^{k-1} d^1 + \textbf{j}^1 d^{k-1} \right) \eta'' - z \textbf{j}^{k-1} \eta'' \right] \nonumber \\
&& + \left( \textbf{l}^{k-1} d^1 + \textbf{l}^1 d^{k-1} \right) \eta'' - z \textbf{l}^{k-1} \eta'' + \left( \textbf{K}^{k-1} d^0 + \textbf{K}^0 d^{k-1} \right) \eta' + \mathcal{D}^{k-2}, \makebox[1cm]{}\\
E^k & = & - \mathcal{E}^\star : \mathcal{C} \left( \textbf{u}^0_z \otimes \nabla d^{k} + \textbf{u}^{k}_z \otimes \nabla d^0 \right) + \left( k^k d^0 + k^0 d^k \right) \eta'\,, \\
A^{k-1} &=& - \mu^{k-1} - \left(\Delta d^0 c^{k-1}_z + \Delta d^{k-1} c^0_z\right) \nonumber\\
&&-2 \left(\nabla d^0 \cdot \nabla c^{k-1}_z + \nabla d^{k-1} \cdot \nabla c^0_z\right)
 + f^{k-1}(c^0, \ldots , c^{k-1}) \nonumber \\
&& - \mathcal{E}^\star : \mathcal{C} \left( \textbf{u}^1_z \otimes \nabla d^{k-1} + \textbf{u}^{k-1}_z \otimes \nabla d^1 \right) - \mathcal{E}^\star : \mathcal{C} \left( \nabla \textbf{u}^{k-1} - \mathcal{E}^\star c^{k-1} \right) \nonumber \\
&& + \left(g^{k-1} d^0 + g^0 d^{k-1}\right) \eta' + \left( k^{k-1} d^1 + k^1 d^{k-1} \right) \eta' - z k^{k-1} \eta' + \mathcal{A}^{k-2} \label{Ak-1},\\
B^{k-1} &=& \left(d^{k-1}_t c^0_z + d^0_t c^{k-1}_z\right) - \left(\Delta d^0 \mu^{k-1}_z + \Delta d^{k-1} \mu^0_z\right) \nonumber\\ 
&& - 2 \left(\nabla d^0 \cdot \nabla \mu^{k-1}_z + \nabla d^{k-1} \cdot \nabla \mu^0_z\right) + \left(d^1 h^{k-1} + d^{k-1} h^1\right) \eta'' \nonumber\\
&& - z h^{k-1} \eta'' + \left(L^{k-1} d^0 + L^0 d^{k-1}\right) \eta' + \mathcal{B}^{k-2} \,. \label{Bk-1}
\end{eqnarray}
Here $\mathcal{D}^{k-2}$, $\mathcal{A}^{k-2}$, and $\mathcal{B}^{k-2}$ have the following form 
\begin{eqnarray}
\mathcal{D}^{k-2} & = & M^{-1} \Bigg[ - \sum^{k-2}_{\substack{i,j=0 \\ 2 \leq i+j \leq k}}{\left( \mathcal{C} \left( \textbf{u}^i_{zz} \otimes \nabla d^j \right) \right) \nabla d^{k-i-j}} \nonumber \\
&& - \sum^{k-2}_{l=1}{\left(\mathcal{C}_{iji'j'} \, \partial_j (\textbf{u}^{l}_{i'})_z \partial_{j'} d^{k-1-l} \right)_{i=1, \ldots , d}} \nonumber \\
&& - \sum^{k-2}_{l=1}{\left[\left( \mathcal{C}_{iji'j'} \, (\textbf{u}^{l}_{i'})_z \partial_{jj'} d^{k-1-l} \right)_{i=1, \ldots , d} + \left( \mathcal{C} \nabla \textbf{u}^{l}_z \right) \nabla d^{k-1-l}\right]} \nonumber \\
&&\left. + \sum^{k-2}_{l=1}{c^l_z \left( \mathcal{C} \mathcal{E}^\star \right) \nabla d^{k-1-l}} + \sum^{k-2}_{l=2}{\textbf{j}^l d^{k-l} \eta''}\right] + \sum^{k-2}_{l=1}{\textbf{K}^l d^{k-1-l} \eta'} \nonumber \\
&& + \sum^{k-2}_{l=2}{\textbf{l}^l d^{k-l} \eta''} - z \textbf{K}^{k-2} \eta' - M^{-1} \left( \mathcal{C}_{iji'j'} \, \partial_{jj'} \textbf{u}^{k-2}_{i'} \right)_{i=1,\ldots,d} \nonumber \\
&& + M^{-1} \left( \left( \mathcal{C} \mathcal{E}^\star \right) \nabla c^{k-2} \right) + M^{-1} \left( \textbf{P}^+_{k-2} \eta^+_N + \textbf{P}^-_{k-2} \eta^-_N \right) \label{Dk-2}, \\
\mathcal{A}^{k-2} &=& \sum^{k-2}_{i=1}{\left(- \Delta d^i c^{k-1-i}_z - 2 \nabla d^i \cdot \nabla c^{k-1-i}_z + d^i g^{k-1-i} \eta' \right)}\nonumber\\
&& - \sum^{k-2}_{i=2}{\mathcal{E}^\star : \mathcal{C} \left( \textbf{u}^i_z \otimes \nabla d^{k-i} \right)} + \sum^{k-2}_{i=2}{k^i d^{k-i} \eta'}\nonumber \\
&& - \Delta c^{k-2} - z g^{k-2} \eta' ,\label{Ak-2}\\
\mathcal{B}^{k-2} &=& \sum^{k-2}_{i=1}{\left(-d^i_t c^{k-1-i}_z - \Delta d^i \mu^{k-1-i}_z - 2 \nabla d^i \cdot \nabla \mu^{k-1-i}_z \right)}\nonumber\\
&& + \sum^{k-2}_{i=1}{d^i L^{k-1-i} \eta'} + \sum_{i=2}^{k-2}{d^i h^{k-i} \eta''} - z L^{k-2} \eta' \nonumber\\
&& + \left(c^{k-2}_t - \Delta \mu^{k-2}\right) - O^+_{k-2} \eta^+_N - O^-_{k-2} \eta^-_N \,. \label{Bk-2}
\end{eqnarray}
Here we have used the conventions that if the upper limit of the summation is less than the lower limit, then the summation is zero, that $ \mathcal{D}^{-1} = \mathcal{A}^{-1} = \mathcal{B}^{-1} = 0$, that $a^{k-1} b^1 + a^1 b^{k-1} = a^1 b^1$ when $k=2$, and that $a^{k-1} b^0 c^1 + a^1 b^0 c^{k-1} = a^1 b^0 c^1$ when $k=2$. We also use these conventions in the following.\\
Observe that $D^{k-1}$ depends on $d^k$ and $\textbf{j}^k$. This would result in difficulties in the construction of $d^k$. To avoid this, we set 
\begin{eqnarray}
\textbf{j}^k := \left(\mathcal{C} \left( \textbf{l}^0 \otimes \nabla d^k \right)\right) \nabla d^0 + \left(\mathcal{C} \left( \textbf{l}^0 \otimes \nabla d^0 \right)\right) \nabla d^k \quad \mbox{for } k\geq 1 \,. \label{jk}
\end{eqnarray}
Since $\textbf{u}^0_{zz} = d^0 \textbf{l}^0 \eta''$, we obtain
\begin{eqnarray}
D^0 & = & M^{-1} \big[ - \left(\mathcal{C}_{iji'j'} \, \partial_j (\textbf{u}^{0}_{i'})_z \partial_{j'} d^0 \right)_{i=1, \ldots , d} - \left( \mathcal{C}_{iji'j'} \, (\textbf{u}^{0}_{i'})_z \partial_{jj'} d^0 \right)_{i=1, \ldots , d} \nonumber \\
&& \left.- \left( \mathcal{C} \nabla \textbf{u}^{0}_z \right) \nabla d^0 + c^{0}_z \left(\mathcal{C} \mathcal{E}^\star \right) \nabla d^0  \right] - z \textbf{l}^{0} \eta'' +  \textbf{K}^{0} d^0 \eta' \,, \label{D0}
\end{eqnarray}
and since $\textbf{u}^{k-1}_{zz} = ( d^{k-1} \textbf{l}^0 + d^0 \textbf{l}^{k-1} ) \eta'' + D^{k-2}$ for $k \geq 2$, we obtain for $k \geq 2$
\begin{eqnarray}
D^{k-1} & = & M^{-1} \left[ - \left(\mathcal{C}_{iji'j'} \, \partial_j (\textbf{u}^{k-1}_{i'})_z \partial_{j'} d^0 \right)_{i=1, \ldots , d} - \left(\mathcal{C}_{iji'j'} \, \partial_j (\textbf{u}^{0}_{i'})_z \partial_{j'} d^{k-1} \right)_{i=1, \ldots , d} \right. \nonumber \\
&& - \left( \mathcal{C}_{iji'j'} \, (\textbf{u}^{k-1}_{i'})_z \partial_{jj'} d^0 \right)_{i=1, \ldots , d} - \left( \mathcal{C}_{iji'j'} \, (\textbf{u}^{0}_{i'})_z \partial_{jj'} d^{k-1} \right)_{i=1, \ldots , d} \nonumber \\
&& - \left( \mathcal{C} \nabla \textbf{u}^{k-1}_z \right) \nabla d^0 - \left( \mathcal{C} \nabla \textbf{u}^{0}_z \right) \nabla d^{k-1} + c^{k-1}_z \left(\mathcal{C} \mathcal{E}^\star \right) \nabla d^0 + c^{0}_z \left(\mathcal{C} \mathcal{E}^\star \right) \nabla d^{k-1} \nonumber \\
&& - d^0 \left(\mathcal{C} \left( \textbf{l}^{k-1} \otimes \nabla d^0 \right)\right) \nabla d^1 \eta'' - d^0 \left(\mathcal{C} \left( \textbf{l}^{k-1} \otimes \nabla d^1 \right)\right) \nabla d^0 \eta'' \nonumber \\
&&- d^0 \left(\mathcal{C} \left( \textbf{l}^{1} \otimes \nabla d^{k-1} \right)\right) \nabla d^0 \eta''- d^0 \left(\mathcal{C} \left( \textbf{l}^{1} \otimes \nabla d^0 \right)\right) \nabla d^{k-1} \eta'' \nonumber \\
&&- d^0 \left(\mathcal{C} \left( \textbf{l}^{0} \otimes \nabla d^{1} \right)\right) \nabla d^{k-1} \eta'' - d^0 \left(\mathcal{C} \left( \textbf{l}^{0} \otimes \nabla d^{k-1} \right)\right) \nabla d^{1} \eta'' \nonumber \\
&& - z \textbf{j}^{k-1} \eta'' \Big]  + \left( \textbf{l}^{k-1} d^1 + \textbf{l}^1 d^{k-1} \right) \eta'' - z \textbf{l}^{k-1} \eta'' + \left( \textbf{K}^{k-1} d^0 + \textbf{K}^0 d^{k-1} \right) \eta' \nonumber\\
&& - M^{-1} \left[ \left(\mathcal{C} \left( D^{k-2} \otimes \nabla d^0 \right)\right) \nabla d^1 + \left(\mathcal{C} \left( D^{k-2} \otimes \nabla d^1 \right)\right) \nabla d^0 \right]  \nonumber \\
&& - M^{-1} \left[ \left(\mathcal{C} \left( D^{0} \otimes \nabla d^{k-1} \right)\right) \nabla d^0 + \left(\mathcal{C} \left( D^{0} \otimes \nabla d^0 \right)\right) \nabla d^{k-1} \right] + \mathcal{D}^{k-2}, \label{Dk-1} 
\end{eqnarray}
where $\mathcal{D}^{k-2}$ is defined as in (\ref{Dk-2}). \\
To get bounded solutions we will see in the next subsection that it is necessary to require $D^{k-1} = \mathcal{O}(e^{-\alpha \left| z \right|})$ and $B^{k-1} = \mathcal{O}(e^{-\alpha \left| z \right|})$. This is the reason why we add $\epsilon^2\left( O^+ \eta^+_N + O^- \eta_N \right)$ and $\epsilon^2 \left( \textbf{P}^+ \eta^+_N + \textbf{P}^- \eta_N^- \right)$ to handle the terms $c_t^{k-2} - \Delta \mu^{k-2}$ and $- (\mathcal{C}_{iji'j'} \partial_{jj'} \textbf{u}^{k-2}_{i'}) + \left( \mathcal{C} \mathcal{E}^\star \right) \nabla c^{k-2}$ in $\mathcal{B}^{k-2}$ and $\mathcal{D}^{k-2}$ (see Lemma \ref{loesle5} and \ref{loesle4}). 

\subsubsection{Compatibility Conditions} \label{seccomp}

In this part we study the compatibility conditions of the ordinary differential equations (\ref{innere3})-(\ref{innere2}) in order to have bounded solutions. Additionally, we investigate the behavior of the solutions as $z \to \pm \infty$. It will come out that there exists bounded solutions $(c^k,\mu^k,\textbf{u}^k)$ for every $k \in \mathbb{N}$. All the assertions of this subsection can be shown as in \cite{ABC}.

\begin{lemma} \label{loesle5}
Let $D^{k-1}$ and $\mathcal{D}^{k-2}$ be defined as in (\ref{D0}), (\ref{Dk-1}), and (\ref{Dk-2}). Then (\ref{innere3}b) has a bounded solution for $k=1$ in $\Gamma^0(\delta)$ if and only if for all $(x,t) \in \Gamma^0(\delta)$, it holds 
\begin{eqnarray}
0 & = & - \left( \mathcal{C}_{iji'j'} \left[ \partial_j \textbf{u}^{0}_{i'} \right] \partial_{j'} d^0 \right)_{i=1,\ldots,d} - \left( \mathcal{C}_{iji'j'} \left[ \textbf{u}^{0}_{i'} \right] \partial_{jj'} d^0 \right)_{i=1,\ldots,d} \nonumber \\
&& - \left( \mathcal{C} \left[ \nabla \textbf{u}^{0} \right] \right) \nabla d^0 + \left[ c^{0} \right] \left(\mathcal{C} \mathcal{E}^\star \right) \nabla d^0 + M \textbf{l}^{0} + M \textbf{K}^{0} d^0 \,, \label{compu1}
\end{eqnarray}
and for $k \geq 2$ it has a solution in $\Gamma^0(\delta)$ if and only if for all $(x,t) \in \Gamma^0(\delta)$, it holds 
\begin{eqnarray}
M \tilde{\mathcal{D}}^{k-2} & = & - \left( \mathcal{C}_{iji'j'} \left[ \partial_j \textbf{u}^{k-1}_{i'} \right] \partial_{j'} d^0 \right)_{i=1,\ldots,d} - \left( \mathcal{C}_{iji'j'} \left[ \partial_j \textbf{u}^{0}_{i'} \right] \partial_{j'} d^{k-1} \right)_{i=1,\ldots,d} \nonumber \\
&& - \left( \mathcal{C}_{iji'j'} \left[ \textbf{u}^{k-1}_{i'} \right] \partial_{jj'} d^0 \right)_{i=1,\ldots,d} - \left( \mathcal{C}_{iji'j'} \left[ \textbf{u}^{0}_{i'} \right] \partial_{jj'} d^{k-1} \right)_{i=1,\ldots,d} \nonumber \\
&& - \left( \mathcal{C} \left[ \nabla \textbf{u}^{k-1} \right] \right) \nabla d^0 - \left( \mathcal{C} \left[ \nabla \textbf{u}^{0} \right] \right) \nabla d^{k-1} + \left[ c^{k-1} \right] \left(\mathcal{C} \mathcal{E}^\star \right) \nabla d^0 \nonumber \\
&& + \left[ c^{0} \right] \left(\mathcal{C} \mathcal{E}^\star \right) \nabla d^{k-1} + \textbf{j}^{k-1} + M \textbf{l}^{k-1} + M \left( \textbf{K}^{k-1} d^0 + \textbf{K}^0 d^{k-1} \right), \makebox[1cm]{}\label{compu}
\end{eqnarray}
where $\left[ . \right] = \left. . \right|^{z=+\infty}_{z=-\infty}$ and 
\[
\tilde{\mathcal{D}}^{k-2}(x,t) = - \int_\mathbb{R}{\mathcal{D}^{k-2}(z,x,t) \, dz} \, .
\]
Furthermore, for $k=0$ every bounded solution to (\ref{innere3}a) has the form
\begin{eqnarray}
\textbf{u}^0(z,x,t) = \tilde{\textbf{u}}^{0}(x,t) + \textbf{l}^0(x,t) d^0(x,t) \left(\eta(z) - \tfrac{1}{2}\right)  \label{uka}\, ,
\end{eqnarray}
and if (\ref{compu}) is satisfied, then for $k \geq 1$ every solution to (\ref{innere3}b) has the form
\begin{eqnarray}
\textbf{u}^k(z,x,t) = \tilde{\textbf{u}}^{k}(x,t) + \left( \textbf{l}^k d^0 + \textbf{l}^0 d^k \right)(x,t) \left(\eta(z) - \tfrac{1}{2}\right) + \textbf{u}^{k-1}_\ast(z,x,t)  \label{uk}\, ,
\end{eqnarray}
where $\tilde{\textbf{u}}^{k}(x,t)$, $k \geq 0$, is an arbitrary function and $\textbf{u}^{k-1}_\ast(z,x,t)$ is a special solution depending only on functions of order lower than $k$ and is uniquely determined by the normalization 
\begin{eqnarray}
\int_\mathbb{R}{\textbf{u}^{k-1}_\ast(z,x,t) \theta'_0(z) \, dz} = 0 \quad \forall (x,t) \in \Gamma^0(\delta)\, .\label{intuast}
\end{eqnarray}
In addition, there exists some $\textbf{u}^{\pm}_{\ast (k-1)}(x,t)$ depending only on functions of order lower than $k$ such that
\begin{eqnarray}
D^m_x D^n_t D^l_z \left( \textbf{u}^{k-1}_\ast(\pm z,x,t) - \textbf{u}^{\pm}_{\ast (k-1)}(x,t) \right) = \mathcal{O}(e^{-\alpha z}) \quad \mbox{as } z \to \infty  \label{uast}
\end{eqnarray}
for all $m,n,l \geq 0$ and for all $(x,t) \in \Gamma^0(\delta)$ provided $\left(c^i, c^\pm_i, \textbf{u}^i, \textbf{u}^\pm_i \right)$, $i=1, \ldots,k-1$, satisfy the matching condition (\ref{match1}).
\end{lemma}

\proof The first assertion of the lemma follows from \cite[Lemma 4.3]{ABC} and the identities $\int_\mathbb{R}{z \eta'' \, dz} = - \int_\mathbb{R}{\eta' \, dz} = -1$ and $\int_\mathbb{R}{\eta'' \,dz} = 0$ and the fact that all terms involving the derivatives with respect to $z$ tend to zero exponentially fast.\\
The second assertion of the lemma follows from the second assertion of \cite[Lemma 4.3]{ABC}, the inner-outer matching condition (\ref{match1}), and the definition of $\textbf{P}^\pm_{k-2}$ (therefore $D^m_x D^n_t D^l_z \mathcal{D}^{k-2} = \mathcal{O}(e^{-\alpha \left|z\right|})$ as $z \to \pm \infty$). 
\makebox[1cm]{} \hfill $\Box$\\

From now we set
\begin{eqnarray}
k^k(x,t) := \left\{
\begin{array}{l@{\quad\mbox{ if }}l}
 \mathcal{E}^\star : \mathcal{C} \left( \textbf{l}^0 \otimes \nabla d^0 \right) & k=0\,, \\
 \mathcal{E}^\star : \mathcal{C} \left( \left( \textbf{l}^k \otimes \nabla d^0 \right) + \left( \textbf{l}^0 \otimes \nabla d^k \right) \right) & k\geq 1 
\end{array}\right. \label{kk}
\end{eqnarray}
for all $(x,t) \in \Gamma^0(\delta)$. By equation (\ref{uka}) and (\ref{uk}) it is not difficult to obtain for all $k \geq 0$
\begin{eqnarray}
E^{k}(z,x,t) = - \mathcal{E}^\star : \mathcal{C} \left( (\textbf{u}_\ast^{k-1})_z \otimes \nabla d^0 \right)  \label{Ek}
\end{eqnarray}
for all $(z,x,t) \in \mathbb{R} \times \Gamma^0(\delta)$ and where $\textbf{u}_\ast^{-1} = 0$. In particular, it holds $E^0=0$. Thus  $c^0(z,x,t)=  \theta_0(z)$ is the unique solution to (\ref{innere1}a) satisfying $\lim_{z \to \infty} c^0(\pm z) = \pm 1 $ and the initial condition $c^0(0,x,t)=0$.

\begin{lemma}\label{loesle2}
Let $k\geq 1$ be any integer and $E^k$, $A^{k-1}$, and $\mathcal{A}^{k-2}$ be defined as in (\ref{Ek}),(\ref{Ak-1}), and (\ref{Ak-2}). Then the system
\begin{equation}
\begin{aligned}
& c^k_{zz}(z,x,t) - f'(\theta_0(z))\, c^k(z,x,t) = (E^k + A^{k-1})(z,x,t) \quad \forall z \in \mathbb{R},\\
& c^k(0,x,t)=0, \quad c^k(.,x,t) \in L^\infty(\mathbb{R}) \label{systemck}
\end{aligned}
\end{equation}
has a unique solution for $k=1$ in $\Gamma^0(\delta)$ if and only if 
\begin{eqnarray}
0 &=& - \overline{\mu}^{0} - \sigma \Delta d^{0} - \mathcal{E}^\star : \mathcal{C} \left( \nabla \overline{\textbf{u}}^{0} \right) + \eta_0 g^{0} d^0 \quad \mbox{in } \Gamma^0(\delta)\,, \label{loes1a}
\end{eqnarray}
and it has a unique solution for $k\geq 2$ in $\Gamma^0(\delta)$ if and only if for all $(x,t) \in \Gamma^0(\delta)$
\begin{eqnarray}
\widetilde{\mathcal{A}}^{k-2} &=& - \overline{\mu}^{k-1} - \sigma \Delta d^{k-1} - \mathcal{E}^\star : \mathcal{C} \nabla \overline{\textbf{u}}^{k-1} + \eta_0 \left(g^{k-1} d^0 + g^0 d^{k-1}\right) \nonumber \\
&& + \eta_0 d^0 \mathcal{E}^\star : \mathcal{C} \left(M^{-1} \left[  \left(\mathcal{C} \left( \textbf{l}^{k-1} \otimes \nabla d^0 \right)\right) \nabla d^1 + \left(\mathcal{C} \left( \textbf{l}^{k-1} \otimes \nabla d^1 \right)\right) \nabla d^0 \right.\right. \nonumber \\
&&  + \left(\mathcal{C} \left( \textbf{l}^{1} \otimes \nabla d^{k-1} \right)\right) \nabla d^0 + \left(\mathcal{C} \left( \textbf{l}^{1} \otimes \nabla d^0 \right)\right) \nabla d^{k-1} + \left(\mathcal{C} \left( \textbf{l}^{0} \otimes \nabla d^{1} \right)\right) \nabla d^{k-1} \nonumber \\
&& \left.\left. + \left(\mathcal{C} \left( \textbf{l}^{0} \otimes \nabla d^{k-1} \right)\right) \nabla d^{1} \right] \otimes \nabla d^0   - \textbf{l}^1 \otimes \nabla d^{k-1} - \textbf{l}^{k-1} \otimes \nabla d^1 \right) \,, \label{loes1}
\end{eqnarray}
where
\begin{align*}
\overline{\mu}^{k-1}(x,t) &= \frac{1}{2} \int_\mathbb{R}{\mu^{k-1}(z,x,t) \theta'_0\left(z\right) dz}\,,\\
\sigma &= \frac{1}{2} \int_\mathbb{R}{\left(\theta'_0\left(z\right)\right)^2 dz}\,,\\
\overline{\textbf{u}}^{k-1}(x,t) & =  \frac{1}{2} \int_\mathbb{R}{\textbf{u}^{k-1}(z,x,t) \theta'_0(z) \, dz} \, ,\\
\eta_0 &= \frac{1}{2} \int_\mathbb{R}{\eta'(z) \theta'_0\left(z\right) dz}\,,\\
\widetilde{\mathcal{A}}^{k-2}(x,t) &= \frac{1}{2} \int_\mathbb{R}\big[\Delta d^0 c^{k-1}_z + 2 \nabla d^0\cdot \nabla c^{k-1}_z - f^{k-1}(c^0,\ldots,c^{k-1}) - \mathcal{E}^\star : \mathcal{C} \mathcal{E}^\star c^{k-1}\\
& \quad - \mathcal{A}^{k-2}\big] \theta'_0 - \left[  c^{k-1}_z \mathcal{E}^\star : \mathcal{C} \left( (\mathcal{C} \mathcal{E}^\star) \nabla d^0 \otimes \nabla d^0 \right) - \mathcal{E}^\star : \mathcal{C} \left( \mathcal{D}^{k-2} \otimes \nabla d^0 \right)\right. \\
& \quad + \mathcal{E}^\star : \mathcal{C} \left(M^{-1} \left[ (\mathcal{C}_{iji'j'} \partial_j (\textbf{u}^{k-2}_{\ast,i'})_z \partial_{j'} d^0)_{i=1,\ldots,d} + (\mathcal{C}_{iji'j'} (\textbf{u}^{k-2}_{\ast,i'})_z \partial_{jj'} d^0)_{i=1,\ldots,d} \right.\right. \\
& \quad + (\mathcal{C} \nabla (\textbf{u}^{k-2}_{\ast})_z) \nabla d^0 + \left( \mathcal{C} \left( D^{k-2} \otimes \nabla d^0 \right) \right) \nabla d^1 \\
& \quad \left.\left.\left. + \left( \mathcal{C} \left( D^{k-2} \otimes \nabla d^1 \right) \right) \nabla d^0 \right] \otimes \nabla d^0\right) + \mathcal{E}^\star : \mathcal{C} \left( D^{k-2} \otimes \nabla d^1 \right)  \right] \theta_0 \, dz\,.
\end{align*}
In addition, if (\ref{loes1a}) for $k=1$ or (\ref{loes1}) for $k \geq 2$ is satisfied and\\ $\left(c^0, c^\pm_0, \mu^0, \mu^\pm_0, \textbf{u}^0, \textbf{u}^\pm_0 \right)$ ,$\ldots$ ,$\left(c^{k-1}, c^\pm_{k-1}, \mu^{k-1}, \mu^\pm_{k-1}, \textbf{u}^{k-1}, \textbf{u}^\pm_{k-1}\right)$ satisfy the matching condition (\ref{match1}), then the unique solution $c^k$ to (\ref{innere1}b) satisfies the matching condition (\ref{match1}) where $c^\pm_k$ is given by (\ref{defcpmj}).
\end{lemma}

\begin{remark}
In the proof we verify that $\int_{\mathbb{R}}{D^0(z) \theta_0(z) \, dz} = 0$. Hence $\widetilde{\mathcal{A}}^{0}$ is independent of $d^1$ for $k=2$. Note that $\widetilde{\mathcal{A}}^{k-2}$, $k\geq 2$, also depends on $c^{k-1}$. But later we will see that $c^{k-1}$ can be considered as a quantity only depending on functions of lower order.
\end{remark}

\proof We prove the first assertion only for $k \geq 2$. For $k=1$ one can apply the same techniques as in the case $k \geq 2$. Due to \cite[Lemma 4.1]{ABC} the system (\ref{systemck}) has a unique solution if and only if
\[
0 = \frac{1}{2} \int_\mathbb{R}{ \left(E^k(z,x,t) + A^{k-1}(z,x,t) \right) \theta'_0(z) \, dz}
\]
for all $(x,t) \in \Gamma^0(\delta)$. Using the equation (\ref{Ek}) and integration by parts yields
\begin{eqnarray*}
\frac{1}{2} \int_{\mathbb{R}}{E^k(z) \theta'_0(z) \, dz} = \frac{1}{2} \int_\mathbb{R}{ \mathcal{E}^\star : \mathcal{C} \left( D^{k-1}(z) \otimes \nabla d^0 \right) \theta_0(z) \, dz} \
\end{eqnarray*}
since $\lim_{z \to \pm\infty} (\textbf{u}^{k-1}_\ast)_z(z) = 0$. For computing the right-hand side we use (\ref{Dk-1}). It holds for all $k \geq 1$
\[
\int_\mathbb{R}{\textbf{u}^{k-1}_z(z) \theta_0(z) \, dz} = \int_\mathbb{R}{(\textbf{u}^{k-2}_\ast)_z(z) \theta_0(z) \, dz} \,,
\]
where we have used (\ref{uk}) and 
\[
\int_\mathbb{R}{\eta'(z) \theta_0(z) \, dz} = \left. \left( \eta - \tfrac{1}{2} \right) \theta_0 \right|^{z= +\infty}_{z=-\infty} - \int_\mathbb{R}{\left( \eta(z) - \tfrac{1}{2} \right) \theta'_0(z) \,dz} =0 
\]
due to (\ref{eigeneta}). Furthermore, we apply the integrals
\begin{align*}
\int_\mathbb{R}{c^0_z(z) \theta_0(z) \, dz} & = \int_\mathbb{R}{ \theta'_0(z) \theta_0(z) \, dz} = 0 , & \int_\mathbb{R}{\eta''(z) \theta_0(z) \, dz} & = -2 \eta_0 \,,\\
\int_\mathbb{R}{z \eta''(z) \theta_0(z) \, dz}  & = 0 , & \int_\mathbb{R}{D^0(z) \theta_0(z) \, dz} & = 0\,,
\end{align*}
where the third identity follows from integration by parts and (\ref{eigeneta}) and where the last identity follows from (\ref{D0}) and the identities above. Therefore we get
\begin{eqnarray*}
\lefteqn{ \frac{1}{2} \int_{\mathbb{R}}{E^k(z) \theta'_0(z) \, dz}} \\
& = & \eta_0 d^0 \mathcal{E}^\star : \mathcal{C} \left(M^{-1} \left[  \left(\mathcal{C} \left( \textbf{l}^{k-1} \otimes \nabla d^0 \right)\right) \nabla d^1 + \left(\mathcal{C} \left( \textbf{l}^{k-1} \otimes \nabla d^1 \right)\right) \nabla d^0 \right.\right. \\
&&  + \left(\mathcal{C} \left( \textbf{l}^{1} \otimes \nabla d^{k-1} \right)\right) \nabla d^0 + \left(\mathcal{C} \left( \textbf{l}^{1} \otimes \nabla d^0 \right)\right) \nabla d^{k-1} + \left(\mathcal{C} \left( \textbf{l}^{0} \otimes \nabla d^{1} \right)\right) \nabla d^{k-1} \\
&& \left.\left. + \left(\mathcal{C} \left( \textbf{l}^{0} \otimes \nabla d^{k-1} \right)\right) \nabla d^{1} \right] \otimes \nabla d^0 \right) - \eta_0 \mathcal{E}^\star : \mathcal{C} \left(\left( \textbf{l}^{k-1} d^1 + \textbf{l}^1 d^{k-1} \right) \otimes \nabla d^0 \right) \\
&& + \frac{1}{2} \int_\mathbb{R}{c^{k-1}_z \theta_0 \, dz} \, \mathcal{E}^\star : \mathcal{C} \left( (\mathcal{C} \mathcal{E}^\star) \nabla d^0 \otimes \nabla d^0 \right) + \frac{1}{2} \mathcal{E}^\star : \mathcal{C} \left( \int_\mathbb{R}{ \mathcal{D}^{k-2} \theta_0 \, dz} \otimes \nabla d^0 \right) \\
&& -  \frac{1}{2} \mathcal{E}^\star :  \mathcal{C} \int_{R}{\left[ M^{-1} \left[ (\mathcal{C}_{iji'j'} \partial_j (\textbf{u}^{k-2}_{\ast,i'})_z \partial_{j'} d^0)_{i=1,\ldots,d} + (\mathcal{C}_{iji'j'} (\textbf{u}^{k-2}_{\ast,i'})_z \partial_{jj'} d^0)_{i=1,\ldots,d} \right. \right. }\\
&& \left. \left. + (\mathcal{C} \nabla (\textbf{u}^{k-2}_\ast)_z) \nabla d^0 \right] \otimes \nabla d^0 \right] \theta_0(z) \, dz \\&& - \frac{1}{2} \mathcal{E}^\star : \mathcal{C} \left( M^{-1} \int_{\mathbb{R}}\left[ \left( \mathcal{C} \left( D^{k-2} \otimes \nabla d^0 \right) \right) \nabla d^1 \right.\right.\\ 
& & \left.+ \left( \mathcal{C} \left( D^{k-2} \otimes \nabla d^1 \right) \right) \nabla d^0 \right] \theta_0 \, dz \otimes \nabla d^0 \bigg)\,. 
\end{eqnarray*}
To compute the integral $\frac{1}{2} \int_\mathbb{R}{A^{k-1} \theta'_0 \, dz}$, we use the identity (\ref{Ak-1}). We apply the definition (\ref{kk}) and the identity (\ref{uk}) to simplify $A^{k-1}$
\begin{eqnarray*}
\lefteqn{- \mathcal{E}^\star : \mathcal{C} \left( \textbf{u}^1_z \otimes  \nabla d^{k-1} + \textbf{u}^{k-1}_z \otimes  \nabla d^{1} \right) + \left( k^{k-1} d^1 + k^1 d^{k-1} \right) \eta' }\\
& = & - \mathcal{E}^\star : \mathcal{C} \left( (\textbf{u}^0_\ast)_z \otimes \nabla d^{k-1} + (\textbf{u}^{k-2}_\ast)_z \otimes \nabla d^{1} \right) \makebox[2cm]{}\\
&& - d^0 \mathcal{E}^\star : \mathcal{C} \left( \textbf{l}^1 \otimes \nabla d^{k-1} + \textbf{l}^{k-1} \otimes \nabla d^1 \right) \eta' \\
&& + \mathcal{E}^\star : \mathcal{C} \left( d^1 \textbf{l}^{k-1} \otimes \nabla d^0 + d^{k-1} \textbf{l}^1 \otimes \nabla d^0 \right) \eta' \,.
\end{eqnarray*}
As above we can show that
\[
- \int_{\mathbb{R}}{ \mathcal{E}^\star : \mathcal{C} \left( (\textbf{u}^0_\ast)_z \otimes \nabla d^{k-1} + (\textbf{u}^{k-2}_\ast)_z \otimes \nabla d^{1} \right) \theta'_0 \, dz} = \int_\mathbb{R}{\mathcal{E}^\star : \mathcal{C} \left( D^{k-2} \otimes \nabla d^{1} \right) \theta_0 \, dz}\,.
\]
Hence we get by definition of $\sigma$ and $\eta_0$
\begin{eqnarray*}
\lefteqn{\frac{1}{2} \int_\mathbb{R}{A^{k-1}(z) \theta'_0(z) \, dz}} \\
& = & - \overline{\mu}^{k-1} - \sigma \Delta d^{k-1} - \mathcal{E}^\star : \mathcal{C} \nabla \overline{\textbf{u}}^{k-1} - \eta_0 d^0 \mathcal{E}^\star : \mathcal{C} \left( \textbf{l}^1 \otimes \nabla d^{k-1} + \textbf{l}^{k-1} \otimes \nabla d^1 \right) \\
&& + \eta_0 \mathcal{E}^\star : \mathcal{C} \left(\left( \textbf{l}^{k-1} d^1 + \textbf{l}^1 d^{k-1} \right) \otimes \nabla d^0 \right) + \eta_0 \left( g^{k-1} d^0 + g^0 d^{k-1} \right)\\
&& + \frac{1}{2} \int_\mathbb{R}{\mathcal{E}^\star : \mathcal{C} \left( D^{k-2} \otimes \nabla d^{1} \right) \theta_0 \, dz} + \frac{1}{2} \int_\mathbb{R}\left(- \Delta d^0 c^{k-1}_z - 2 \nabla d^0 \cdot \nabla c^{k-1}_z \right. \\
&& \left. +f^{k-1}(c^0,\ldots,c^{k-1}) + \mathcal{E}^\star : \mathcal{C} \mathcal{E}^\star c^{k-1} + \mathcal{A}^{k-2} \right) \theta'_0 \, dz \,.
\end{eqnarray*}
This shows the first assertion for $k \geq 2$.\\
For all $k \geq 1$ the second assertion of the lemma follows from the second assertion of \cite[Lemma 4.1]{ABC} and
\begin{align*}
 E^k(\pm z) + A^{k-1}(\pm z)  & =  - \mu^{k-1}(\pm z) + f^{k-1}(c^0(\pm z), \ldots, c^{k-1}(\pm z))\\
& \quad  -  \mathcal{E}^\ast : \mathcal{C} \left( \nabla \textbf{u}^{k-1}(\pm z) - \mathcal{E}^\ast c^{k-1}(\pm z) \right) + \mathcal{A}^{k-2}(\pm z) + \mathcal{O}(e^{-\alpha z}) \\
& \to  - \mu^\pm_{k-1} + f^{k-1}(c^\pm_0, \ldots, c^\pm_{k-1})\\
& \quad - \mathcal{E}^\ast : \mathcal{C} \left( \nabla \textbf{u}^\pm_{k-1} - \mathcal{E}^\ast c^\pm_{k-1} \right) - \Delta c^\pm_{k-2} \,,
\end{align*}
as $z \to \infty$, since all the terms involving the derivatives with respect to $z$ tend to zero exponentially fast and 
\[
c^j(\pm z) = c^\pm_j + \mathcal{O}(e^{-\alpha z}), \quad \mu^j(\pm z) = \mu^\pm_j + \mathcal{O}(e^{-\alpha z}), \quad \textbf{u}^j(\pm z) = \textbf{u}^\pm_j + \mathcal{O}(e^{-\alpha z}),
\]
as $z \to \infty$ and for all $j \in \left\lbrace 0, \ldots , k-1\right\rbrace $. Using the outer expansion (\ref{defcpmj}) completes the proof.
\makebox[1cm]{} \hfill $\Box$\\

\begin{lemma}\label{loesle4}
Let $B^{k-1}$ and $\mathcal{B}^{k-2}$ be defined as in (\ref{Bk-1}) and (\ref{Bk-2}). Then, for $k \geq 1$, the ordinary differential equation (\ref{innere2}b) has a bounded solution if and only if for all $(x,t) \in \Gamma^0(\delta)$
\begin{eqnarray}
&& d^{k-1}_t - \frac{1}{2} \left(\Delta d^0 \left[\mu^{k-1}\right] + \Delta d^{k-1}\left[\mu^0\right]\right) - \nabla d^0 \cdot \left[\nabla \mu^{k-1}\right] \nonumber \\
&& - \nabla d^{k-1} \cdot \left[\nabla \mu^0\right] + \frac{1}{2} h^{k-1} +\frac{1}{2} \left( d^0 L^{k-1} + d^{k-1} L^0 \right) = \widetilde{\mathcal{B}}^{k-2}\,,\label{loes2}
\end{eqnarray}
where $[.]= .\left.\right|^{z=+\infty}_{z=-\infty}$ and
\[
\widetilde{\mathcal{B}}^{k-2}(x,t)= - \frac{1}{2} d^0_t(x,t) \left[c^{k-1}(x,t)\right] - \frac{1}{2} \int_\mathbb{R}{\mathcal{B}^{k-2}(z,x,t) \,dz}
\]
if $k\geq 2$ and $\widetilde{\mathcal{B}}^{-1}= 0$. In addition, every bounded solution to (\ref{innere2}a) can be written as
\begin{eqnarray}
\mu^0(z,x,t)= \tilde{\mu}^0(x,t) + d^0(x,t) h^0(x,t) \left(\eta(z)-1/2\right) \,,\label{loes4a}
\end{eqnarray}
and if (\ref{loes2})is satisfied, then every solution to (\ref{innere2}b) can be written as
\begin{eqnarray}
\mu^k(z,x,t)= \tilde{\mu}^k(x,t) + \left(d^0 h^k + d^k h^0\right)\!(x,t)\left(\eta(z)-1/2\right) + \mu_\ast^{k-1}(z,x,t)\,,\label{loes4}
\end{eqnarray}
where $\tilde{\mu}^k(x,t)$ is an arbitrary function and $\mu_\ast^{k-1}(z,x,t)$ is a special solution depending only on $(c^\pm_i, c^i,  \mu^\pm_i, \mu^i, d^i, h^i, g^i)$ for $i\leq k-1$ and is uniquely determined by the normalization
\begin{eqnarray}
\int_{\mathbb{R}}{\mu^{k-1}_\ast(z,x,t) \theta'_0\left(z\right) dz} = 0\quad \forall (x,t) \in \Gamma^0(\delta)\,.\label{loes5}
\end{eqnarray}
Furthermore, there exists some $\mu_{\ast(k-1)}^\pm$ depending only on $(c^\pm_i, c^i, \mu^\pm_i, \mu^i, d^i, h^i, g^i)$ for $i\leq k-1$ such that 
\[
D_x^m D_t^n D_z^l \left[\mu^{k-1}_\ast(\pm z,x,t) - \mu_{\ast(k-1)}^\pm (x,t)\right] = \mathcal{O}(e^{- \alpha z}) \quad \mbox{as }z \to \infty 
\]
for all $m,n,l \geq 0$ and for all $(x,t) \in \Gamma^0(\delta)$, provided $(c^i, c^\pm_i, \mu^i, \mu^\pm_i)$, $i=1,\ldots,k-1$, satisfy the matching condition (\ref{match1}).
\end{lemma}

\proof The first assertion of the lemma follows from \cite[Lemma 4.3]{ABC} and the identities $\int_\mathbb{R}{\theta'_0 \, dz} = 2$, $\int_\mathbb{R}{\eta'' \, dz} = 0$, and $\int_\mathbb{R}{z \eta'' \, dz} = - \int_\mathbb{R}{\eta' \, dz} = -1$. \\
The second assertion follows from the second assertion of \cite[Lemma 4.3]{ABC}, the inner-outer matching condition (\ref{match1}), and the definition of $O^\pm_{k-2}$ (therefore $D^m_x D^n_t D^l_z \mathcal{B}^{k-2} = \mathcal{O}(e^{-\alpha \left|z\right|})$ as $z \to \pm \infty$). 
\makebox[1cm]{} \hfill $\Box$\\

\subsubsection{Boundary-Layer Expansion} \label{boundaryexp}

As for the inner expansion we assume that near the boundary $\partial_T \Omega$ the solutions $(c^\epsilon,\mu^\epsilon,\textbf{u}^\epsilon)$ have for every $\epsilon \in (0,1]$ the expansion
\begin{align*}
c^\epsilon(x,t) & = c_B^\epsilon\left(\frac{d_B^\epsilon(x,t)}{\epsilon}, x,t\right),& c_B^\epsilon(z,x,t) & = 1 + \sum^\infty_{i=1}{\epsilon^i c_B^i(z,x,t)},\\
(\mu^\epsilon, \textbf{u}^\epsilon)(x,t) & = (\mu_B^\epsilon, \textbf{u}^\epsilon_B)\left(\frac{d_B^\epsilon(x,t)}{\epsilon}, x,t\right), & (\mu_B^\epsilon, \textbf{u}^\epsilon_B)(z,x,t) & = \sum^\infty_{i=0}{\epsilon^i (\mu_B^i,\textbf{u}^i_B)(z,x,t)}, 
\end{align*}
where $(x,t) \in \overline{\partial_T \Omega(\delta)}$ and $z \in (-\infty,0]$. \\
To match the boundary-layer and outer expansion, we require as $z \to -\infty$ 
\begin{eqnarray}
D_x^m D_t^n D_z^l \left[(c^k_B, \mu^k_B, \textbf{u}^k_B)( z,x,t) - (c^+_k,\mu^+_k, \textbf{u}^+_k)(x,t)\right] & = & \mathcal{O}(e^{\alpha z}) \, ,\label{matchbound1}
\end{eqnarray}
for all $(x,t) \in \overline{\partial_T \Omega(\delta)}$ and all $k,m,n,l \in \left\{0,\ldots,\bar{K}\right\}$ where $\bar{K}$ depends on the order of the expansion.\\
Similarly to the inner expansion, we define $M_B(x) = \left(\mathcal{C}_{iji'j'} \, \partial_j d_B(x) \, \partial_{j'} d_B(x) \right)^d_{i,i'=1}$ where an analogous proof shows the invertibility of $M_B$. Since $\partial \Omega$ is known, we do not require a series expansion for $d_B$. Also note that $d_B$ is time independent. We substitute the expansion of $(c^\epsilon,\mu^\epsilon,\textbf{u}^\epsilon)$ into (\ref{system1})-(\ref{system3}). Then a calculation similar to the inner expansion gives us 
\begin{align}
\textbf{u}^k_{B,zz}(z,x,t) & = \textbf{A}^{k-1}_B(z,x,t), & k \geq 0 \, ,\label{boundODE1}\\
c^k_{B,zz}(z,x,t) - f'(1) c^k_B(z,x,t) & = B^{k-1}_B(z,x,t), & k \geq 1\, , \label{boundODE2}\\
\mu^k_{B,zz}(z,x,t) & = C^{k-1}_B(z,x,t), & k \geq 0 \label{boundODE3}
\end{align}
for all $(x,t) \in \partial_T \Omega$ and $z \in (-\infty,0)$ where $\textbf{A}^{k-1}_B$, $B^{k-1}_B$ and $C^{k-1}_B$ are defined by (\ref{A_B}), (\ref{B_B}), and (\ref{C_B}) (with $j-1=k-1$).\\
Since $\nabla d_B$ is the unit outer normal to $\partial \Omega$, we obtain $\left. \frac{\partial}{\partial n} \right|_{\partial \Omega} = \epsilon^{-1} \frac{\partial}{\partial z} + \nabla d_B \cdot \nabla $. Also observe that $c^\epsilon_B \! \left(\frac{d_B(x)}{\epsilon},x,t\right) = c^\epsilon_B(0,x,t)$ for all $(x,t) \in \partial_T \Omega$ (analogously for $\mu^\epsilon_B$ and $\textbf{u}^\epsilon_B$). Therefore to satisfy the boundary conditions on $\partial_T \Omega$, we require
\begin{align}
\textbf{u}^k_{B}(0,x,t) & = 0 && \forall (x,t) \in \partial_T \Omega, k \geq 0, \label{boundinitial1}\\
c^k_{B,z}(0,x,t) & = - \nabla d_B \cdot \nabla c^{k-1}_B(0, S_B(x),t) && \forall (x,t) \in \overline{\partial_T \Omega(\delta)}, k \geq 1, \label{boundinitial2}\\
\mu^k_{B,z}(0,x,t) & = - \nabla d_B \cdot \nabla \mu^{k-1}_B(0,x,t) && \forall (x,t) \in \partial_T \Omega, k \geq 0\,. \label{boundinitial3}
\end{align}

\begin{remark}
The boundary condition (\ref{boundinitial2}) is necessary only for $x \in \partial_T\Omega$. We use the same boundary condition as in \cite{ABC}. The boundary condition (\ref{boundinitial2}) has the advantage that we obtain a unique smooth solution in $(x,t)$ and if for all $k = 0,\ldots,j-1$, $\left.\frac{\partial}{\partial n} c^+_k \right|_{\partial_T\Omega} = 0$ and $c^k_B(z,x,t)$, $\mu^k_B(z,x,t)$, and $\textbf{u}^k_B(z,x,t)$ are independent of $z$ (therefore $A^{j-1}(z,x,t)$ is independent of $z$), then so is $c^j_B(z,x,t)$. For $\mu^j_B$ and $\textbf{u}^j_B$ we do not require uniqueness since we specify $\mu^j_B$ and $\textbf{u}^j_B$ directly.
\end{remark}

Now let us show that the ordinary differential equations (\ref{boundODE1})-(\ref{boundODE3}) with initial values (\ref{boundinitial1})-(\ref{boundinitial3}) have bounded solutions.

\begin{lemma} \label{lemmauB}
Let $j \geq 0$ be any integer. Assume that for all $i = 0,\ldots,j-1$, the functions $c^+_i, \textbf{u}^+_i, c^i_B, \textbf{u}^i_B$ are known, smooth, and satisfy the matching condition (\ref{matchbound1}) and the outer-expansion equation $\Div \left( \mathcal{C} \mathcal{E}(\textbf{u}^+_i) \right) = \Div \left( \mathcal{C} \mathcal{E}^\star c^+_i \right) $. Let $\textbf{F}^{j-1}$ be defined as in (\ref{randbedu}) and assume that $\textbf{u}^+_{j}$ satisfies the boundary condition (\ref{linHele7}) (also in the case $j=0$). Also assume that $\textbf{u}^{-2}_B = \textbf{u}^{-1}_B = c_B^{-2} = c^{-1}_B = 0$, and $\textbf{u}^i_B$, $i=0,\ldots,j-1$, are defined by (\ref{DefuB}) (with $j=i$) for all $z \in (-\infty,0]$ and $(x,t) \in \overline{\partial_T \Omega(\delta)}$, and where $\textbf{A}^{i-1}_B$ is defined as in (\ref{A_B}). Then for known smooth $\textbf{u}^+_j$ the function $\textbf{u}^j_B$ defined by (\ref{DefuB}) satisfies for $k=j$ the boundary-expansion equation (\ref{boundODE1}), the boundary condition (\ref{boundinitial1}), and the matching condition (\ref{matchbound1}).
\end{lemma}

\proof First let us show that $\textbf{F}^{j-1}$ and $\textbf{u}^i_B$, $i=0,\ldots,j$, are well-defined and smooth. Since $(c^+_i, \textbf{u}^+_i, c^i_B, \textbf{u}^i_B)$ satisfy the matching conditions (\ref{matchbound1}) and the outer-expansion equation  $\Div \left( \mathcal{C} \mathcal{E}(\textbf{u}^+_i) \right) = \Div \left( \mathcal{C} \mathcal{E}^\star c^+_i \right) $ for $i=0,\ldots,j-1$, we obtain 
\[
\left| \Div \left( \mathcal{C} \nabla \textbf{u}^i_B \right) - (\mathcal{C} \mathcal{E}^\star) \nabla c^i_B \right|  \leq C e^{\alpha z} \quad \forall z \in (-\infty,0], \, \forall (x,t) \in \overline{\partial_T \Omega(\delta)} \, ,
\]
and therefore one concludes by definition of $\textbf{A}^i_B$ and the fact that all terms involving the derivatives with respect to $z$ tend to zero exponentially fast  
\[
\left| \textbf{A}^i_B \right|  \leq C e^{\alpha z} \quad \forall z \in (-\infty,0], \, \forall (x,t) \in \overline{\partial_T \Omega(\delta)} 
\]
for $i=0,\ldots,j-1$ and some $C>0$. Therefore the integrals defining $\textbf{F}^{j-1}$ and $\textbf{u}^i_B$, $i=0,\ldots,j$, are well-defined and smooth. The same arguments as above and the definition of $\textbf{u}^j_B$ yield the matching condition (\ref{matchbound1}) for $\textbf{u}^j_B$ and $\textbf{u}^+_j$. By an easy calculation we obtain (\ref{boundODE1}) for $k=j$. Finally, the boundary condition (\ref{boundinitial1}) immediately follows from the condition (\ref{linHele7}) and the definition of $\textbf{F}^{j-1}$ in (\ref{randbedu}).
\makebox[1cm]{} \hfill $\Box$\\

\begin{remark} \label{remarku0u1}
Since $\textbf{A}^{-1}_B = 0$, it follows $\textbf{u}^0_B(z,x,t) = \textbf{u}^+_0(x,t)$ for all $(x,t) \in \partial_T \Omega(\delta)$ and therefore it holds $\textbf{A}^0_B =0$ since $c^0_B(z,x,t) = 1$. This yields $\textbf{u}^1_B(z,x,t) = \textbf{u}^+_1(x,t)$ for all $(x,t) \in \partial_T \Omega(\delta)$ . 
\end{remark}

\begin{lemma} \label{lemmacB}
Let $j \geq 1$ be any integer. Assume that for all $i=0,\ldots,j-1$, the functions $c^+_i, \mu^+_i, c^i_B, \mu^i_B, \textbf{u}^+_i, \textbf{u}^i_B$ are known, smooth, and satisfy the matching condition (\ref{matchbound1}). Then for $k=j$, the boundary-layer expansion equation (\ref{boundODE2}) subject to the boundary condition (\ref{boundinitial2}) has a unique bounded solution $c^j_B$ for $z \in (-\infty,0]$ and all $(x,t) \in \overline{\partial_T \Omega(\delta)}$. In addition, the solution satisfies the matching condition (\ref{matchbound1}) where $c^+_j$ is defined by $(\ref{defcpmj})$.
\end{lemma}

\proof The assertions can be shown as Lemma 4.5 in \cite{ABC}. Also see Lemma 3.2.14 in \cite{St} for details.
\makebox[1cm]{} \hfill $\Box$\\

\begin{lemma} \label{lemmamuB}
Let $j \geq 0$ be any integer. Assume that for all $i=0,\ldots,j-1$, the functions $c^+_i, \mu^+_i, c^i_B, \mu^i_B$ are known, smooth, and satisfy the matching condition (\ref{matchbound1}) and the outer-expansion equation $\Delta \mu^+_i = \partial_t c^+_i$. Let $G^{j-1}$ be defined as in (\ref{randbedmu}) and assume that $\mu^+_{j-1}$ satisfies the boundary condition (\ref{linHele3}) (also in the case $j-1=0$). Also assume that $\mu^{-2}_B = \mu^{-1}_B = c^{-2}_B = c^{-1}_B = 0$ and $\mu^i_B$, $i=0,\ldots,j-1$, are defined by (\ref{DefcB}) (with $j=i$) for all $z \in (-\infty,0]$ and $(x,t) \in \overline{\partial_T \Omega(\delta)}$ and where $C^{i-1}_B$ is defined as in (\ref{C_B}). Then for known smooth $\mu^+_j$ the function $\mu^j_B$ defined by (\ref{DefcB}) satisfies for $k=j$ the boundary-expansion equation (\ref{boundODE3}), the boundary condition (\ref{boundinitial3}), and the matching condition (\ref{matchbound1}).
\end{lemma}

\proof The assertions can be shown as Lemma 4.6 in \cite{ABC}. Also see Lemma 3.2.15 in \cite{St} for details.
\makebox[1cm]{} \hfill $\Box$\\

Note that for the required boundary conditions (\ref{linHele7}) and (\ref{linHele3}), we only need functions of lower order. Provided the functions of lower order are known, we have given boundary conditions for $\mu^+_k$ and $\textbf{u}^+_k$.

\subsection{Rigorous Construction of Approximate Solutions}

\subsubsection{Overview: Basic Steps for Solving Expansions of each Order} \label{secbasicsteps}

As in \cite{ABC} we define the unknown functions 
\[
\mathcal{V}^j \equiv \left(c^\pm_j, c^j, c^j_B, \mu^\pm_j,\mu^j,\mu^j_B, \textbf{u}^\pm_j, \textbf{u}^j, \textbf{u}^j_B,d^j, g^j, L^j, h^j, \textbf{l}^j, \textbf{K}^j\right)
\]
for each $j \geq 0$ recursively. We call $\mathcal{V}^j$ the $j$th order expansion. \\
We assume that $\mathcal{V}^i$ are known for $i=0,\ldots,j-1$, and the corresponding outer, inner, and boundary-layer expansion equations, the inner-outer matching conditions, and the outer-boundary matching conditions are satisfied for $i=0,\ldots,j-1$. Moreover, we assume that the compatibility conditions (\ref{compu1}) if $j=1$ or (\ref{compu}) if $j>1$, (\ref{loes1a}) if $j=1$ or (\ref{loes1}) if $j>1$, and (\ref{loes2}) are satisfied for $k=j$. In the following we derive the equations for $\mathcal{V}^j$. As in \cite{ABC} we first construct $(c^j,c^\pm_j,c^j_B)$. Then we continue with $(\textbf{u}^j,\textbf{u}^\pm_j,\textbf{u}^j_B)$ and $(\mu^j,\mu^\pm_j,\mu^j_B)$ and finally, we show how to find $d^j$. More precisely, we carry out the following steps:\\
\textbf{Step 1:} After determining $\textbf{u}^{j-1}_\ast$ by the known quantities $\mathcal{V}^i$, $i \leq j-1$, we can determine $(c^j,c^\pm_j,c^j_B)$. Therefore we can consider $(c^j,c^\pm_j,c^j_B)$ as known quantities depending only on $\mathcal{V}^i$, $i \leq j-1$.\\
For Steps 2-9 we assume that $d^j$ is known.\\
\textbf{Step 2:} We obtain $\textbf{u}^j$ by equation (\ref{uk}).\\
\textbf{Step 3:} By the matching condition (\ref{match1}) for $\textbf{u}^j$ and $\textbf{u}^\pm_j$, we determine $\tilde{\textbf{u}}^j$ and $\textbf{l}^j$ in $\Gamma^0(\delta)$.\\
\textbf{Step 4:} The compatibility condition (\ref{compu1}) if $j=0$ or (\ref{compu}) if $j>0$ for $k=j+1$ yields the boundary condition for $\textbf{u}^\pm_j$ on $\Gamma^0$.\\
\textbf{Step 5:} The outer expansion equation $\Div \left( \mathcal{C} \mathcal{E}(\textbf{u}^\pm_j)\right) = \Div \left( \mathcal{C} \mathcal{E}^\star c^\pm_j\right) $, the boundary condition $\textbf{u}^+_j = \textbf{F}^{j-1}$ on $\partial_T \Omega$, and Step 3 and 4 give us an elliptic boundary problem for $\textbf{u}^\pm_j$. So we can determine $\textbf{u}^\pm_j$ uniquely.\\
\textbf{Step 6:} By solving the compatibility condition (\ref{loes1a}) if $j=0$ or (\ref{loes1}) if $j>0$ for $k=j+1$ on $\Gamma^0$, we can determine $\overline{\mu}^j$ on $\Gamma^0$.\\
\textbf{Step 7:} We obtain $\tilde{\mu}^j = \overline{\mu}^j$ in $\Gamma^0(\delta)$. Then by equation (\ref{loes4}), $\mu^j$ is uniquely determined on $\Gamma^0$.\\
\textbf{Step 8:} The matching condition yields the boundary condition for $\mu^\pm_j$ on $\Gamma^0$. Together with the boundary condition $\frac{\partial}{\partial n} \mu^+_j = G^{j-1}$ on $\partial_T \Omega$ and the outer expansion $\Delta \mu^\pm_j = \partial_t c^\pm_j$, we can determine $\mu^\pm_j$ uniquely.\\
\textbf{Step 9:} Again the matching condition and equation (\ref{loes4}) prescribe how to define $\tilde{\mu}^j$ and $h^j$ in $\Gamma^0$. Therefore $\mu^j$ is uniquely determined by (\ref{loes4}). It is not difficult to see that the identity for $\overline{\mu}^j$ on $\Gamma^0$ in Step 6 is satisfied.\\
The following step yields $d^j$.\\
\textbf{Step 10:} By the compatibility condition (\ref{loes2}), we obtain an evolution equation for $d^j$ on $\Gamma^0$. For $j=0$ we require that $d^0$ is a signed spatial distance function and for $j \geq 1$ we require (\ref{djbedingung}). Now we have a system of equations which determine $(\textbf{u}^\pm_j,\mu^\pm_j,d^j)$ uniquely.\\
\textbf{Step 11:} Note that the compatibility conditions (\ref{compu1}) if $j=0$ or (\ref{compu}) if $j>0$, (\ref{loes1a}) if $j=0$ or (\ref{loes1}) if $j>0$, and (\ref{loes2}) are satisfied on $\Gamma^0$ only. We are able to determine $g^j$, $L^j$, and $\textbf{K}^j$ such that these compatibility conditions are satisfied in $\Gamma^0(\delta)$.\\
\textbf{Step 12:} Finally, by (\ref{DefuB}) and (\ref{DefcB}) we get $\textbf{u}^j_B$ and $\mu^j_B$. This completes the construction of $\mathcal{V}^j$.\\

After motivating the construction of $\mathcal{V}^j$ in the Steps 1-12, we verify that $\mathcal{V}^j$ satisfies all the corresponding outer, inner, and boundary-layer expansion equations, the inner-outer matching conditions, and the outer-boundary matching conditions for $k=j$. In addition, we show that the compatibility conditions (\ref{compu1}) if $j=0$ or (\ref{compu}) if $j>0$, (\ref{loes1a}) if $j=0$ or (\ref{loes1}) if $j>0$, and (\ref{loes2}) are also satisfied for $k=j+1$.\\
In the next two subsections we carry out Steps 1-12 in detail.

\subsubsection{The Zero-th Order Expansion} \label{subsec0th}

In this subsection we solve for $\mathcal{V}^0$.\\
\textbf{Step 1:} We already know the leading order term of the outer and boundary-layer expansion for $c^\epsilon$. Since $E^0=0$ we also know the inner expansion. More precisely, it means that $c^\pm_0(x,t) = \pm 1$ for $(x,t) \in Q^\pm_0 \cup \Gamma^0(\delta)$, $c^0_B(z,x,t) = 1$ for $(z,x,t) \in (-\infty,0] \times \overline{\partial_T(\delta)}$, and $c^0(z,x,t)=\theta_0(z)$ for $(z,x,t) \in \mathbb{R} \times \Gamma^0(\delta)$.\\
Now we assume that $\Gamma^0$ and therefore $d^0$ are known. We obtain the construction of $d^0$ below.\\
\textbf{Step 2:} Since the compatibility condition is always satisfied for $\textbf{u}^0$, equation (\ref{uka}) yields 
\begin{eqnarray}
\textbf{u}^0(z,x,t) = \tilde{\textbf{u}}^0(x,t) + d^0(x,t) \textbf{l}^0(x,t) \left( \eta(z) - 1/2 \right) \quad \forall (x,t) \in \Gamma^0(\delta) \,, \label{uzero}
\end{eqnarray}
where we choose $\tilde{\textbf{u}}^0$ and $\textbf{l}^0$ later.\\
\textbf{Step 3:} From equation (\ref{uzero}) we get the condition
\[
\lim_{z \to \infty} \textbf{u}^0(\pm z,x,t) = \tilde{\textbf{u}}^0(x,t) \pm \frac{1}{2} d^0(x,t) \textbf{l}^0(x,t) \quad \forall (x,t) \in \Gamma^0(\delta) \, .
\]
In order to satisfy the matching condition on $\Gamma^0$ (note that $d^0=0$ on $\Gamma^0$), we get the condition
\[
\textbf{u}^+_0(x,t) = \lim_{z \to \infty} \textbf{u}^0(z,x,t) = \tilde{\textbf{u}}^0(x,t) = \lim_{z \to \infty} \textbf{u}^0(-z,x,t) = \textbf{u}^-_0(x,t) \quad \forall (x,t) \in \Gamma^0 \, .
\]
To satisfy the matching condition in $\Gamma^0(\delta) \backslash \Gamma^0$, it is necessary and sufficient to define 
\begin{align}
\tilde{\textbf{u}}^0(x,t) & := \frac{1}{2} \left( \textbf{u}^+_0(x,t) + \textbf{u}^-_0(x,t) \right) & \forall (x,t) & \in \Gamma^0(\delta) \, , \nonumber \\
\textbf{l}^0(x,t) & := \frac{1}{d^0(x,t)} \left( \textbf{u}^+_0(x,t) - \textbf{u}^-_0(x,t) \right) & \forall (x,t) & \in \Gamma^0(\delta) \backslash \Gamma^0\,. \label{l0delta}
\end{align}
The natural way to define $\textbf{l}^0$ on $\Gamma^0$ is
\begin{equation}
\left. \textbf{l}^0 \right|_{\Gamma^0} := \nabla \left( \textbf{u}^+_0 - \textbf{u}^-_0 \right) \nabla d^0 = \frac{\partial}{\partial \nu} \left( \textbf{u}^+_0 - \textbf{u}^-_0 \right), \label{l0}
\end{equation}
where $\nabla d^0 = \nu$ is the unit outward normal of $\Gamma^0_t$.\\
\textbf{Step 4:} We consider the compatibility condition (\ref{compu1}) on $\Gamma^0$ for $k=1$. One gets
\begin{eqnarray}
0 & = & - \left( \mathcal{C}_{iji'j'} \left( \partial_j (\textbf{u}^+_0)_{i'} - \partial_j (\textbf{u}^-_0)_{i'} \right) \partial_{j'}d^0 \right)_{i=1,\ldots,d} - \left( \mathcal{C} \left(\nabla \textbf{u}^+_0 - \nabla \textbf{u}^-_0 \right) \right) \nabla d^0 \nonumber \\
&& + 2 \left( \mathcal{C} \mathcal{E}^\star \right) \nabla d^0 + M \left( \partial_\nu \textbf{u}^+_0 - \partial_\nu \textbf{u}^-_0 \right) , \label{compu0}
\end{eqnarray}
where we have used the definition of $\textbf{l}^0$ on $\Gamma^0$ and $ \left[ \textbf{u}^0 \right] = 0$ on $\Gamma^0$. We can simplify this equation for $(x,t) \in \Gamma^0$. Let $\left\lbrace \tau_1, \ldots , \tau_{d-1}\right\rbrace $ be an orthonormal basis of the tangent space of $\Gamma^0_t$. Then it holds for all $\textbf{u} \in C^1(\Omega;\mathbb{R}^d)$
\begin{eqnarray}
\nabla \textbf{u} = \left( \partial_\nu \textbf{u} \otimes \nu \right) + \sum^{d-1}_{i=1}{\left( \partial_{\tau_i} \textbf{u} \otimes \tau_i \right)} \label{gradu}.
\end{eqnarray}
Since $\textbf{u}^+_0 = \textbf{u}^-_0$ on $\Gamma^0$, we obtain $\left( \partial_{\tau_i} \textbf{u}^+_0 \otimes \tau_i \right) - \left( \partial_{\tau_i} \textbf{u}^-_0 \otimes \tau_i \right) = 0$ on $\Gamma^0$ for all $i=1,\ldots,d-1$, and therefore we have
\begin{equation}
\nabla \textbf{u}^+_0 - \nabla \textbf{u}^-_0 = \left( \partial_\nu \textbf{u}^+_0 \otimes \nu \right) - \left( \partial_\nu \textbf{u}^-_0 \otimes \nu \right) \quad \mbox{on } \Gamma^0. \label{randablu}
\end{equation}
Using this property and $\nabla d^0 = \nu$ on $\Gamma^0$, we obtain on $\Gamma^0$
\begin{eqnarray*}
\lefteqn{ \left( \mathcal{C}_{iji'j'} \left( \partial_j (\textbf{u}^+_0)_{i'} - \partial_j (\textbf{u}^-_0)_{i'} \right) \partial_{j'}d^0 \right)_{i=1,\ldots,d}} \\ 
& = & \left( \mathcal{C}_{iji'j'} \left( \partial_\nu (\textbf{u}^+_0)_{i'} - \partial_\nu (\textbf{u}^-_0)_{i'} \right) \partial_j d^0 \partial_{j'}d^0 \right)_{i=1,\ldots,d} \\
& = & \left( \mathcal{C} \left( \nabla \textbf{u}^+_0 - \nabla \textbf{u}^-_0 \right) \right) \nu \,.
\end{eqnarray*}
By the definition of $M$ (see (\ref{defM})), we have 
\begin{eqnarray*}
M \left( \partial_\nu \textbf{u}^+_0 - \partial_\nu \textbf{u}^-_0 \right) & = & \left( \mathcal{C}_{iji'j'} \left( \partial_\nu (\textbf{u}^+_0)_{i'} - \partial_\nu (\textbf{u}^-_0)_{i'} \right) \partial_j d^0 \partial_{j'}d^0 \right)_{i=1,\ldots,d} \\
& = & \left( \mathcal{C} \left( \nabla \textbf{u}^+_0 - \nabla \textbf{u}^-_0 \right) \right) \nu \, ,
\end{eqnarray*}
where the second equation follows as above. So equation (\ref{compu0}) turns into
\begin{eqnarray*}
0 = \left( 2 \mathcal{C} \mathcal{E}^\star - \mathcal{C} \left( \nabla \textbf{u}^+_0 - \nabla \textbf{u}^-_0 \right) \right) \nu \quad \mbox{on }\Gamma^0\,.
\end{eqnarray*}
\textbf{Step 5:} Note that $\textbf{F}^{-1} = 0 $. In order to apply Lemma \ref{lemmauB}, we prescribe the boundary condition $ \left. \textbf{u}^+_0 \right|_{\partial_T \Omega} = 0$. Moreover, we require $\Div \left(\mathcal{C} \nabla \textbf{u}^\pm_0\right) =0$ in $Q^\pm_0$ due to the outer expansion. Therefore we can determine uniquely $\textbf{u}^\pm_0$ by solving the elliptic boundary problem
\begin{align*}
\Div \left(\mathcal{C} \nabla \textbf{u}^\pm_0\right) & =0 && \mbox{in } Q^\pm_0\, , \\
\left[ (\mathcal{C} \nabla \textbf{u}^\pm_0 - \mathcal{C} \mathcal{E}^\star c^\pm_0) \nu \right]_{\Gamma^0_t} = \left[ \textbf{u}^\pm_0\right]_{\Gamma^0_t} & = 0 && \mbox{on }\Gamma^0_t, t \in [0,T] \,,\\
\textbf{u}^+_0 & = 0 && \mbox{on } \partial \Omega \times [0,T]\, .
\end{align*}
\textbf{Step 6:} Due to the definitions of $\textbf{u}^0$ in Step 2 and $\overline{\textbf{u}}^0$ in Lemma \ref{loesle2} and the property $\int_\mathbb{R}{\left(\eta - 1/2\right) \theta_0'\left(z\right) dz} = 0$, we can conclude that $\tilde{\textbf{u}}^0 = \overline{\textbf{u}}^0 $. Then the compatibility condition (\ref{loes1a}) on $\Gamma^0$ for $k=1$ is
\begin{eqnarray*}
\overline{\mu}^0(x,t) = - \sigma \Delta d^0 - \frac{1}{2} \mathcal{E}^\star : \mathcal{C} \left( \nabla \textbf{u}^+_0 + \nabla \textbf{u}^-_0 \right) \quad \forall (x,t) \in \Gamma^0 \,.
\end{eqnarray*}
It is a known fact that 
\[
\Delta d^0 = \Div_{\Gamma^0_t} \nu = - \kappa_{\Gamma^0_t}\, ,
\]
where $\kappa_{\Gamma^0_t}$ is the mean curvature of $\Gamma^0_t$. Then we get the condition
\begin{eqnarray*}
\overline{\mu}^0(x,t) = \sigma \kappa_{\Gamma^0_t} - \frac{1}{2} \mathcal{E}^\star : \mathcal{C} \left( \nabla \textbf{u}^+_0 + \nabla \textbf{u}^-_0 \right) \quad \forall (x,t) \in \Gamma^0\,.
\end{eqnarray*}
\textbf{Step 7:} We get from equation (\ref{loes4a}) 
\begin{eqnarray}
\mu^0(z,x,t) = \tilde{\mu}^0(x,t) + d^0(x,t) h^0(x,t) \left(\eta(z) - 1/2\right) \quad \forall (x,t) \in \Gamma^0(\delta) \,,\label{0-ko1}
\end{eqnarray}
where we choose $\tilde{\mu}^0$ and $h^0$ later. As in Step 6 we obtain $\tilde{\mu}^0 = \overline{\mu}^0$. Note that $d^0=0$ on $\Gamma^0$, and therefore we obtain
\[
\mu^0(z,x,t) = \sigma \kappa_{\Gamma^0_t} - \frac{1}{2} \mathcal{E}^\star : \mathcal{C} \left( \nabla \textbf{u}^+_0 + \nabla \textbf{u}^-_0 \right) \quad \forall (x,t) \in \Gamma^0\,.
\]
\textbf{Step 8:} The equation above and the matching condition lead to
\begin{equation}
\mu^\pm_0(x,t) = \lim_{z \to \infty} \mu^0(\pm z,x,t) =  \sigma \kappa_{\Gamma^0_t} - \frac{1}{2} \mathcal{E}^\star : \mathcal{C} \left( \nabla \textbf{u}^+_0 + \nabla \textbf{u}^-_0 \right) \quad \forall (x,t) \in \Gamma^0 \,. \label{boundmupm0}
\end{equation}
Since $G^{-1}=0$, we obtain $\left.\frac{\partial}{\partial n } \mu^+_0\right|_{\partial_T \Omega} = 0$ (see Lemma \ref{lemmamuB}), and since $c^\pm_0= \pm 1$, the outer expansion for $k=0$ reads $\Delta \mu^\pm_0 =0$ in $Q^\pm_0$. Together with the boundary condition (\ref{boundmupm0}) on $\Gamma^0_t$, we can determine uniquely $\mu^\pm_0$ by solving the elliptic boundary problem
\begin{align*}
\Delta \mu^\pm_0 & = 0 && \mbox{in }Q^\pm_0 \,, \\
\mu^\pm_0 & = \sigma \kappa_{\Gamma^0_t} - \tfrac{1}{2} \mathcal{E}^\star : \mathcal{C} \left( \nabla \textbf{u}^+_0 + \nabla \textbf{u}^-_0 \right) && \mbox{on }\Gamma^0_t, t\in [0,T] \, ,\\
\tfrac{\partial}{\partial n } \mu^+_0 & = 0 && \mbox{on }\partial \Omega \times [0,T]\,. 
\end{align*}
\textbf{Step 9:} From equation (\ref{0-ko1}) we get the condition
\[ 
\lim_{z \to \infty} \mu^0(\pm z ,x,t) = \overline{\mu}^0(x,t) \pm \frac{1}{2} d^0(x,t) h^0(x,t) \quad \forall (x,t) \in \Gamma^0(\delta)\,.
\]
So in order to satisfy the matching condition $\lim_{z \to \infty} \mu^0(\pm z,x,t) = \mu^\pm_0(x,t)$, it is necessary and sufficient to take
\begin{align}
\overline{\mu}^0(x,t) = \tilde{\mu}^0(x,t) &:= \frac{1}{2} \left(\mu_0^+(x,t) + \mu^-_0(x,t)\right) && \forall (x,t)\in \Gamma^0(\delta), \nonumber \\
h^0(x,t) &:= \frac{1}{d^0(x,t)} \left(\mu_0^+(x,t) - \mu^-_0(x,t)\right) && \forall (x,t)\in \Gamma^0(\delta)\backslash \Gamma^0. \label{h0delta}
\end{align}
The natural way to define $h^0$ on $\Gamma^0$ is
\begin{eqnarray}
h^0\left|_{\Gamma^0}\right. := \nabla d^0 \cdot \nabla \left(\mu_0^+(x,t) - \mu^-_0(x,t)\right)
= \frac{\partial}{\partial \nu}\left(\mu_0^+(x,t) - \mu^-_0(x,t)\right) .\label{h0}
\end{eqnarray}
Note that the definition of $\overline{\mu}^0$ satisfies the identity for $\overline{\mu}^0$ in Step 6.\\
\textbf{Step 10:} On $\Gamma^0$ the compatibility condition (\ref{loes2}) for $k=1$ reads
\begin{eqnarray*}
d^0_t(x,t) & =& \frac{1}{2} \Delta d^0 \left[\mu^0\right] + \nabla d^0 \cdot \left[\nabla \mu^0\right] - \frac{1}{2} h^0 \\
& = & \frac{1}{2} \left(\frac{\partial }{\partial \nu} \mu^+_0 - \frac{\partial }{\partial \nu} \mu^-_0\right) \quad \forall (x,t) \in \Gamma^0,
\end{eqnarray*}
where we have used the equations $\left[\mu^0\right]= \left(\mu^+_0-\mu^-_0\right) = 0$ and $\nabla d^0 \cdot \left[\nabla \mu^0\right]= \nabla d^0 \cdot \left(\nabla \mu^+_0- \nabla \mu^-_0\right) = h^0 $ on $\Gamma^0$. Note that the normal velocity of $\Gamma^0_t$ is given by $-d^0_t$ and the unit outer normal $\nu$ by $\nabla d^0$. Therefore $\Gamma^0$, $\mu_0 := \mu^+_0 \chi_{\left\{d^0\geq0\right\}} + \mu^-_0 \chi_{\left\{d^0<0\right\}}$, and $\textbf{u}_0 := \textbf{u}^+_0 \chi_{\left\{d^0\geq0\right\}} + \textbf{u}^-_0 \chi_{\left\{d^0<0\right\}}$ have to solve the problem (\ref{sharpsystem1})-(\ref{sharpsystem7}). Equation (\ref{sharpsystem4}) is satisfied since $ \nu^T \left[ W \mathrm{Id} -\nabla (\textbf{u}_0)^T \mathcal{S} \right]_{\Gamma^0_t} \nu = - \left( \nabla \textbf{u}^+_0 + \nabla \textbf{u}^-_0 \right) : \mathcal{C} \mathcal{E}^\star$ on $\Gamma^0$. We will show this in the proof of Lemma \ref{lemma0teOrd} below. \\
\textbf{Step 11:} Until now we fulfill the compatibility conditions (\ref{compu1}),(\ref{loes1a}), and (\ref{loes2}) for $k=1$ only on $\Gamma^0$. To satisfy the compatibility conditions in $\Gamma^0(\delta) \backslash \Gamma^0$ we set 
\begin{align}
g^0(x,t) & :=\frac{1}{2 \eta_0 d^0}\left(\mu^+_0+\mu^-_0 + 2\sigma \Delta d^0 + \mathcal{E}^\star : \mathcal{C} \left( \nabla \textbf{u}^+_0 + \nabla \textbf{u}^-_0 \right) \right) , \label{g0delta} \\
L^0(x,t) & := - \frac{1}{d^0}\left(2d^0_t - \left(\Delta d^0+ 2 \nabla d^0 \cdot \nabla \right)\left(\mu^+_0-\mu^-_0\right) + h^0\right) , \label{L0delta}\\
\textbf{K}^0(x,t) & := \frac{1}{d^0} M^{-1} ( \left( \mathcal{C}_{iji'j'} \,  \partial_{j} (\textbf{u}^+_0 - \textbf{u}^-_0 )_{i'} \, \partial_{j'} d^0 \right)_{i=1,\ldots,d} \nonumber \\
& \quad + \left( \mathcal{C}_{iji'j'} \left( \textbf{u}^+_{0} - \textbf{u}^-_0 \right)_{i'} \partial_{jj'} d^0 \right)_{i=1,\ldots,d} \nonumber \\
& \quad + \left( \mathcal{C} \left( \nabla \textbf{u}^+_0 - \nabla \textbf{u}^-_0 \right) \right) \nabla d^0 - 2 \left( \mathcal{C} \mathcal{E}^\star \right) \nabla d^0 - M \textbf{l}^0)  \label{K0delta}
\end{align}
for $(x,t) \in \Gamma^0(\delta) \backslash \Gamma^0$. Since the numerators of $g^0$, $L^0$, and $\textbf{K}^0$ vanish on $\Gamma^0$ we can extend $g^0$, $L^0$, and $\textbf{K}^0$ smoothly to $\Gamma^0$ by
\begin{align}
g^0(x,t) & := \frac{1}{2 \eta_0}\nabla d^0 \cdot \nabla \left( \mu^+_0+\mu^-_0 + 2\sigma \Delta d^0 + \mathcal{E}^\star : \mathcal{C} ( \nabla \textbf{u}^+_0 + \nabla \textbf{u}^-_0 ) \right) , \label{g0}\\
L^0(x,t) & := -\nabla d^0 \cdot \nabla \left(2 d^0_t - \left(\Delta d^0 + 2 \nabla d^0 \cdot \nabla \right)\left(\mu^+_0-\mu^-_0\right) + h^0\right) , \label{L0}\\
\textbf{K}^0(x,t) & := M^{-1} \nabla d^0 \cdot \nabla ( \left( \mathcal{C}_{iji'j'} \,  \partial_{j} (\textbf{u}^+_0 - \textbf{u}^-_0 )_{i'} \, \partial_{j'} d^0 \right)_{i=1,\ldots,d} \nonumber \\
& \quad -\left( \mathcal{C}_{iji'j'} \left( \textbf{u}^+_{0} - \textbf{u}^-_0 \right)_{i'} \partial_{jj'} d^0 \right)_{i=1,\ldots,d} \nonumber \\
& \quad  \left( \mathcal{C} \left( \nabla \textbf{u}^+_0 - \nabla \textbf{u}^-_0 \right) \right) \nabla d^0 - 2 \left( \mathcal{C} \mathcal{E}^\star \right) \nabla d^0 - M \textbf{l}^0) \label{K0}
\end{align}
for $(x,t) \in \Gamma^0$.\\ 
\textbf{Step 12:} Observe that $\textbf{A}^{-1}_B = 0$ and $C^{-1}_B = 0$. Then Lemma \ref{lemmauB} and \ref{lemmamuB} yield $\textbf{u}^0_B(z,x,t) = \textbf{u}^+_0(x,t)$ and $\mu^0_B(z,x,t) = \mu^+_0(x,t)$ for all $(x,t) \in \overline{\partial_T \Omega(\delta)}$.\\  
Note that $(c^\pm_0,c^0, c^0_B, \mu^\pm_0,\mu^0, \mu^0_B,\textbf{u}^\pm_0,\textbf{u}^0, \textbf{u}^0_B, \Gamma^0,d^0)$ coincides with the definition in Subsection \ref{secapp}. After motivating the construction of $\mathcal{V}^0$, we obtain the following result.            
    
\begin{lemma} \label{lemma0teOrd}
Let $\Gamma_{00} \subset\subset \Omega$ be a given smooth hypersurface without boundary and assume that the Hele-Shaw problem (\ref{sharpsystem1})-(\ref{sharpsystem7}) admits a smooth solution $(\mu,\textbf{u}, \Gamma)$ in the time interval $[0,T]$. Let $d^0$ be the signed distance from $x$ to $\Gamma_t$ such that $d^0<0$ inside of $\Gamma_t$, and let $\delta$ be a small positive constant such that $\mbox{dist}\left(\Gamma_t, \partial \Omega\right)> 2\delta$ for all $t \in [0,T]$, $d^0$ is smooth in $\Gamma (2\delta):=\left\{(x,t) \in \Omega_T : \left|d^0\right|<2\delta\right\}$ and $\mu^\pm := \left.\mu\right|_{Q^\pm_0}$ and $\textbf{u}^\pm := \left.\textbf{u}\right|_{Q^\pm_0}$ have a smooth extension to $Q^\pm_0 \cup \Gamma(2\delta)$ where $Q^\pm_0:=\left\{(x,t) \in \Omega_T : \pm d^0 > 0 \right\}$. Define $h^0(x,t)$ by (\ref{h0delta}) and (\ref{h0}), $g^0(x,t)$ by (\ref{g0delta}) and (\ref{g0}), $L^0(x,t)$ by (\ref{L0delta}) and (\ref{L0}), $\textbf{l}^0(x,t)$ by (\ref{l0delta}) and (\ref{l0}), and $\textbf{K}^0(x,t)$ by (\ref{K0delta}) and (\ref{K0}), $\textbf{j}^0(x,t)$ by (\ref{jk}), and $k^0(x,t)$ by (\ref{kk}).
Then, for $k=0$, the hypersurface $\Gamma^0$ defined by (\ref{0thinterface}) and $(c^\pm_0,c^0, c^0_B, \mu^\pm_0,\mu^0, \mu^0_B,\textbf{u}^\pm_0,\textbf{u}^0, \textbf{u}^0_B)$ defined by (\ref{0thouter})-(\ref{0thboundary}) satisfy the outer expansion equations in Subsection \ref{secouterexp}, the inner-expansion equations (\ref{innere3})-(\ref{innere2}), the boundary-layer-expansion equations (\ref{boundODE1}),(\ref{boundODE3}), the inner-outer matching condition (\ref{match1}), the outer boundary matching condition (\ref{matchbound1}), and the boundary conditions (\ref{boundinitial1})-(\ref{boundinitial3}). In addition, the compatibility conditions (\ref{compu1}), (\ref{loes1a}), and (\ref{loes2}) for $k=1$ are also satisfied.
\end{lemma}

\proof We verify all the properties by direct calculation.\\
The outer expansion equations in Subsection \ref{secouterexp} are directly satisfied by definition of $(c^\pm_0,\mu^\pm_0,\textbf{u}^\pm_0)$. \\
\textbf{To (\ref{innere3}):} We consider $\Gamma^0(\delta) \backslash \Gamma^0$ and $\Gamma^0$ separately. Then we obtain
\begin{align*}
\left(\textbf{u}^0 - \textbf{l}^0 d^0 \eta\right)_{zz} &= \left(\textbf{u}_0^+  \eta + \textbf{u}^-_0 (1-\eta) - \left(\textbf{u}^+_0- \textbf{u}^-_0\right) \eta \right)_{zz}& \\
& = (\textbf{u}^+_0)_{zz} = 0 & \mbox{in }& \Gamma^0(\delta) \backslash \Gamma^0 \,,\\
\left(\textbf{u}^0 - \textbf{l}^0 d^0 \eta\right)_{zz} &= \left(\textbf{u}_0^+  \eta + \textbf{u}^-_0 (1-\eta) \right)_{zz} = (\textbf{u}^+_0 - \textbf{u}^-_0) \eta'' = 0 & \mbox{on }&\Gamma^0 
\end{align*}
since $\textbf{u}^+_0 = \textbf{u}^-_0$ on $\Gamma^0$.\\
\textbf{To (\ref{innere1}):} The definitions of $c^0$ and $\theta_0$ yield
\[
c^0_{zz} - f(c^0) = \theta''_0 - f(\theta_0) = 0 = E^0 \quad \mbox{in } \Gamma^0(\delta)\,.
\]
\textbf{To (\ref{innere2}):} The proof for $\mu^0$ is analogous to $\textbf{u}^0$.\\
\textbf{To (\ref{boundODE1}),(\ref{boundODE3}):} The assertions immediately follow from the definitions of $\mu^0_B$ and $\textbf{u}^0_B$.\\
\textbf{To (\ref{match1}):}  For $c^0$ and $c^\pm_0$ the assertion follows from Lemma \ref{lemmatheta0}. Since for $z>0$
\begin{eqnarray*}
\mu^0(+z,x,t)-\mu^+_0(x,t) &=& \left(1- \eta(z)\right) \left(\mu^-_0(x,t)-\mu^+_0(x,t)\right),\\
\mu^0(-z,x,t) -\mu^-_0(x,t) &=& \eta(-z) \left(\mu^+_0(x,t)-\mu^-_0(x,t)\right),
\end{eqnarray*}
the assertion follows from (\ref{eta}) for all $(x,t) \in \Gamma^0(\delta)$. The proof for $\textbf{u}^0$ and $\textbf{u}^\pm_0$ is analogous to $\mu^0$ and $\mu^\pm_0$. \\
\textbf{To (\ref{matchbound1}):} The assertions immediately follow from the definitions of $c^0_B$, $\mu^0_B$, and $\textbf{u}^0_B$.\\
\textbf{To (\ref{boundinitial1})-(\ref{boundinitial3}):} The assertions immediately follow from the definitions of $c^0_B$, $\mu^0_B$, and $\textbf{u}^0_B$.\\
\textbf{To (\ref{compu1}):} Due to the definition of $\textbf{K}^0$ the compatibility condition follows for all $(x,t) \in \Gamma^0(\delta) \backslash \Gamma^0$ immediately.
On $\Gamma^0$ we obtain
\begin{eqnarray*}
\lefteqn{ - \left( \mathcal{C}_{iji'j'} \left[ \partial_{j} \textbf{u}^0_{i'} \right] \partial_{j'} d^0 \right)_{i=1,\ldots,d} - \left( \mathcal{C}_{iji'j'} \left[ \textbf{u}^0_{i'} \right] \partial_{jj'} d^0 \right)_{i=1,\ldots,d} }\\
\lefteqn{ - \left( \mathcal{C} \left[ \nabla \textbf{u}^0 \right] \right) \nabla d^0 + \left[ c^0 \right] \left( \mathcal{C} \mathcal{E}^\star \right) \nabla d^0 + M \textbf{l}^0 + M \textbf{K}^0 d^0 }\\
& = & - \left( \mathcal{C}_{iji'j'} \,  \partial_{j} (\textbf{u}^+_0 - \textbf{u}^-_0 )_{i'} \, \partial_{j'} d^0 \right)_{i=1,\ldots,d} - \left( \mathcal{C}_{iji'j'} \left( \textbf{u}^+_{0} - \textbf{u}^-_0 \right)_{i'} \partial_{jj'} d^0 \right)_{i=1,\ldots,d} \\
&& - \left( \mathcal{C} \left( \nabla \textbf{u}^+_0 - \nabla \textbf{u}^-_0 \right) \right) \nabla d^0 + 2 \left( \mathcal{C} \mathcal{E}^\star \right) \nabla d^0 + M \left( \partial_\nu \textbf{u}^+_0 - \partial_\nu \textbf{u}^-_0\right) \\
& = & \left( 2 \mathcal{C} \mathcal{E}^\star - \mathcal{C} \left( \nabla \textbf{u}^+_0 - \nabla \textbf{u}^-_0 \right) \right) \nu = 0 \, ,
\end{eqnarray*}
where the second equation follows in the same way as the calculation in Step 4 since $\textbf{u}^+_0 = \textbf{u}^-_0$ on $\Gamma^0$ and where the last equation follows from $\left[\mathcal{S} \nu \right]_{\Gamma^0_t}= 0$.\\
\textbf{To (\ref{loes1a}):} In $\Gamma^0(\delta) \backslash \Gamma^0$ we obtain
\begin{eqnarray*}
\lefteqn{ - \overline{\mu}^0 - \sigma \Delta d^0 - \mathcal{E}^\star : \mathcal{C} \left( \nabla \overline{\textbf{u}}^0 \right)+ \eta_0 d^0 g^0} \\
& = & - \frac{1}{2} \mu^+_0 \int_\mathbb{R}{ \eta \theta'_0\;dz} -  \frac{1}{2} \mu^-_0 \int_\mathbb{R}{ \left(1-\eta\right) \theta'_0\;dz} - \sigma \Delta d^0 \\
&& - \frac{1}{2} \mathcal{E}^\star : \mathcal{C} \left( \nabla \left( \textbf{u}^+_0 \int_\mathbb{R}{ \eta \theta'_0\;dz} + \textbf{u}^-_0 \int_\mathbb{R}{ \left(1-\eta\right) \theta'_0\,dz} \right) \right) \\
&& + \frac{\eta_0 d^0}{2 \eta_0 d^0} \left( \mu^+_0+\mu^-_0 + 2\sigma \Delta d^0 + \mathcal{E}^\star : \mathcal{C} \left( \nabla \textbf{u}^+_0 + \nabla \textbf{u}^-_0 \right) \right) = 0 \,,
\end{eqnarray*}
where we have used the the definitions of $\overline{\mu}^0$ and $\overline{\textbf{u}}^0$ in Lemma \ref{loesle2}, the properties of $\theta_0$, and (\ref{eigeneta}). On $\Gamma^0$ we get
\begin{eqnarray*}
\lefteqn{ - \overline{\mu}^0 - \sigma \Delta d^0 - \mathcal{E}^\star : \mathcal{C} \left( \nabla \overline{\textbf{u}}^0 - \mathcal{E}^\star \overline{c}^0 \right) + \eta_0 d^0 g^0 } \\
& = & - \frac{1}{2} \mu^+_0 - \frac{1}{2} \mu^-_0 + \sigma \kappa_{\Gamma^0_t} - \frac{1}{2} \mathcal{E}^\star : \mathcal{C} \left( \nabla \textbf{u}^+_0 + \nabla \textbf{u}^-_0 \right) \\
& = & - \frac{1}{2} \nu^\top \left[ W \mathrm{Id} - (\nabla \textbf{u}_0)^T \mathcal{S}\right]_{\Gamma^0_t} \nu - \frac{1}{2} \mathcal{E}^\star : \mathcal{C} \left( \nabla \textbf{u}^+_0 + \nabla \textbf{u}^-_0 \right) 
\end{eqnarray*}
since $\Delta d^0 = - \kappa_{\Gamma^0_t}$ and $d^0= 0$ on $\Gamma^0$ and where we have applied the boundary condition (\ref{sharpsystem4}). Let us show that the term on the right-hand side vanishes. First note that
\begin{eqnarray*}
\nu^T \left[ W \mathrm{Id} \right]_{\Gamma^0_t} \nu = \frac{1}{2} \mathcal{C} \nabla \textbf{u}^+_0 : \nabla \textbf{u}^+_0 - \frac{1}{2} \mathcal{C} \nabla \textbf{u}^-_0 : \nabla \textbf{u}^-_0 - \mathcal{C} \nabla \textbf{u}^+_0 : \mathcal{E}^\star - \mathcal{C} \nabla \textbf{u}^-_0 : \mathcal{E}^\star \,.
\end{eqnarray*}
Due to the definition of $\mathcal{S}$ we get
\begin{eqnarray*}
\nu^T \left[ (\nabla \textbf{u}_0)^T \mathcal{S} \right]_{\Gamma^0_t} \nu & = & \nu^T (\nabla \textbf{u}^+_0)^T (\mathcal{C} \nabla \textbf{u}^+_0) \nu - \nu^T (\nabla \textbf{u}^-_0)^T (\mathcal{C} \nabla \textbf{u}^-_0) \nu \\
&& - \nu^T (\nabla \textbf{u}^+_0)^T (\mathcal{C} \mathcal{E}^\star) \nu - \nu^T (\nabla \textbf{u}^-_0)^T (\mathcal{C} \mathcal{E}^\star) \nu \\
& = & \left( \partial_\nu \textbf{u}^+_0 \otimes \nu \right) : \mathcal{C} \nabla \textbf{u}^+_0 - \left( \partial_\nu \textbf{u}^-_0 \otimes \nu \right) : \mathcal{C} \nabla \textbf{u}^-_0 \\
&& - \left( \partial_\nu \textbf{u}^+_0 \otimes \nu \right) : \mathcal{C} \mathcal{E}^\star - \left( \partial_\nu \textbf{u}^-_0 \otimes \nu \right) : \mathcal{C} \mathcal{E}^\star \\
& = & \nabla \textbf{u}^+_0 : \mathcal{C} \nabla \textbf{u}^+_0 - \nabla \textbf{u}^-_0 : \mathcal{C} \nabla \textbf{u}^+_0 \\
&& + \left( \partial_\nu \textbf{u}^-_0 \otimes \nu \right) : \mathcal{C} \nabla \textbf{u}^+_0 - \left( \partial_\nu \textbf{u}^-_0 \otimes \nu \right) : \mathcal{C} \nabla \textbf{u}^-_0 \\
&& - \left( \partial_\nu \textbf{u}^+_0 \otimes \nu \right) : \mathcal{C} \mathcal{E}^\star - \left( \partial_\nu \textbf{u}^-_0 \otimes \nu \right) : \mathcal{C} \mathcal{E}^\star \, ,
\end{eqnarray*}
where the last equality follows from (\ref{randablu}) since $\left[ \textbf{u}_0\right]_{\Gamma^0_t} = 0$ on $\Gamma^0_t$. Together we have 
\begin{eqnarray*}
\lefteqn{ \nu^T \left[ W \mathrm{Id} -\nabla (\textbf{u}_0)^T \mathcal{S} \right]_{\Gamma^0_t} \nu} \\
& = & - \frac{1}{2} \left( \nabla \textbf{u}^+_0 - \nabla \textbf{u}^-_0 \right) : \mathcal{C} \left( \nabla \textbf{u}^+_0 - \nabla \textbf{u}^-_0 \right) - \left( \partial_\nu \textbf{u}^-_0 \otimes \nu \right) : \mathcal{C} \left( \nabla \textbf{u}^+_0 - \nabla \textbf{u}^-_0 \right) \\
&& - \left( \nabla \textbf{u}^+_0 + \nabla \textbf{u}^-_0 \right) : \mathcal{C} \mathcal{E}^\star + \left( \partial_\nu \textbf{u}^+_0 \otimes \nu \right) : \mathcal{C} \mathcal{E}^\star + \left( \partial_\nu \textbf{u}^-_0 \otimes \nu \right) : \mathcal{C} \mathcal{E^\star} \\
& = & - \frac{1}{2} \left( \left( \partial_\nu \textbf{u}^+_0 + \partial_\nu \textbf{u}^-_0 \right) \otimes \nu \right) : \mathcal{C} \left( \nabla \textbf{u}^+_0 - \nabla \textbf{u}^-_0 \right) - \left( \nabla \textbf{u}^+_0 + \nabla \textbf{u}^-_0 \right) : \mathcal{C} \mathcal{E}^\star \\
&& + \left( \left( \partial_\nu \textbf{u}^+_0 + \partial_\nu \textbf{u}^-_0 \right) \otimes \nu \right) : \mathcal{C} \mathcal{E}^\star \\
& = & - \left( \nabla \textbf{u}^+_0 + \nabla \textbf{u}^-_0 \right) : \mathcal{C} \mathcal{E}^\star \, ,
\end{eqnarray*}
where we have used (\ref{randablu}) in the second equation and $\left[\mathcal{S} \nu\right]_{\Gamma^0_t} = 0$ in the last equation. Therefore we fulfill the compatibility condition (\ref{loes1a}) on $\Gamma^0$. \\
\textbf{To (\ref{loes2}):} In $\Gamma^0(\delta)\backslash \Gamma^0$ we directly obtain by the definition of $L^0$. On $\Gamma^0$ it holds due to the interface condition (\ref{sharpsystem3}) 
\begin{eqnarray*}
\lefteqn{ d^0_t - \frac{1}{2} \Delta d^0 \left[\mu^0\right] - \nabla d^0 \cdot \left[\nabla \mu^0\right] + \frac{1}{2} h^0 + \frac{1}{2} d^0 L^0} \\
&=& -V - \left(\frac{\partial }{\partial \nu}\mu_0^+ - \frac{\partial}{\partial \nu}  \mu^-_0\right) + \frac{1}{2} \left(\frac{\partial }{\partial \nu}\mu_0^+ - \frac{\partial }{\partial \nu}\mu^-_0\right)=0 
\end{eqnarray*}
since $d^0_t = -V$ on $\Gamma^0_t$.
This completes the proof.
\makebox[1cm]{} \hfill $\Box$

\subsubsection{The Higher-Order Expansions}\label{subsecjth}

Let $j\geq1$ be an integer. Assume that $\mathcal{V}^0,\ldots,\mathcal{V}^{j-1}$ are known and that the matching conditions for $k=0, \ldots, j-1$ and the compatibility conditions of Lemma \ref{loesle5}, \ref{loesle2}, and \ref{loesle4} are satisfied for $k=j$ .\\
\textbf{Step 1:} After determining $\textbf{u}_\ast^{j-1}$, we can calculate $c^j$ in $\mathbb{R} \times \Gamma^0(\delta)$ by equation (\ref{systemck}). Equation (\ref{defcpmj}) gives us $c^\pm_j$ in $Q^\pm_0$ directly. The proof of Lemma \ref{lemmacB} or \cite[Lemma 4.5]{ABC} respectively shows us how we can obtain $c^j_B$ in $(-\infty,0] \times \overline{\partial_T \Omega(\delta)}$. So we can assume that $c^j$, $c^\pm_j$, and $c^j_B$ are known functions which only depend on $\mathcal{V}^0, \ldots , \mathcal{V}^{j-1}$. \\
For Step 2-9 we assume that $d^j$ is known. The construction of $d^j$ is shown afterwards.\\
\textbf{Step 2:} Equation (\ref{uk}) yields for all $(x,t) \in \Gamma^0(\delta)$
\begin{eqnarray}
\textbf{u}^j(z,x,t) = \tilde{\textbf{u}}^j(x,t) + \left( \textbf{l}^j d^0 + \textbf{l}^0 d^j \right)\!(x,t) \left( \eta(z) - 1/2 \right) + \textbf{u}^{j-1}_\ast(z,x,t) \, , \label{uj}
\end{eqnarray}
where we define $\tilde{\textbf{u}}^j$ and $\textbf{l}^j$ later.\\
\textbf{Step 3:} Due to the definition of $\eta$, we obtain from equation (\ref{uj}) and Lemma \ref{loesle5}
\begin{eqnarray*}
\lim_{z \to \infty} \textbf{u}^j(\pm z,x,t) = \tilde{\textbf{u}}^j(x,t) \pm \frac{1}{2} \left( \textbf{l}^j d^0 + \textbf{l}^0 d^j \right) + \textbf{u}^\pm_{\ast(j-1)}(x,t) \quad \forall (x,t) \in \Gamma^0(\delta) .
\end{eqnarray*}
By the inner-outer matching condition we get on $\Gamma^0$
\begin{eqnarray}
\textbf{u}^\pm_j(x,t) = \tilde{\textbf{u}}^j(x,t) \pm \frac{1}{2} \textbf{l}^0 d^j + \textbf{u}^\pm_{\ast(j-1)}(x,t) \,. \label{sprunguj}
\end{eqnarray}
For satisfying the inner-outer matching condition on $\Gamma^0(\delta)\backslash \Gamma^0$, it is necessary and sufficient to define
\begin{eqnarray*}
\tilde{\textbf{u}}^j(x,t) & := & \frac{1}{2} \left( \textbf{u}^+_j + \textbf{u}^-_j - \textbf{u}^+_{\ast(j-1)} - \textbf{u}^-_{\ast(j-1)} \right) \quad \mbox{in } \Gamma^0(\delta) \, ,\\
\textbf{l}^j(x,t) & := & \left\{ \begin{array}{l @{\quad}l }
\frac{1}{d^0} \left(- d^j \textbf{l}^0 + \textbf{u}_j^+ - \textbf{u}^-_j - \textbf{u}^+_{\ast(j-1)} + \textbf{u}^-_{\ast(j-1)}\right) & \mbox{in } \Gamma^0(\delta)\backslash \Gamma^0\\
\nabla d^0 \cdot \nabla \left( - d^j \textbf{l}^0 + \textbf{u}_j^+ - \textbf{u}^-_j - \textbf{u}^+_{\ast(j-1)} + \textbf{u}^-_{\ast(j-1)} \right) & \mbox{on } \Gamma^0.
\end{array} \right. 
\end{eqnarray*}
Note that the numerator in the definition of $\textbf{l}^j$ vanishes on $\Gamma^0$. So the definition of $\textbf{l}^j$ is natural on $\Gamma^0$. \\
\textbf{Step 4:} On $\Gamma^0$ the compatibility condition (\ref{compu}) reads for $k=j+1$
\begin{eqnarray}
M \tilde{\mathcal{D}}^{j-1} & = & - \left( \mathcal{C}_{ili'l'} \left[ \partial_l \textbf{u}^{j}_{i'} \right] \partial_{l'} d^0 \right)_{i=1,\ldots,d} - \left( \mathcal{C}_{ili'l'} \left[ \partial_l \textbf{u}^{0}_{i'} \right] \partial_{l'} d^{j} \right)_{i=1,\ldots,d} \nonumber \\
&& - \left( \mathcal{C}_{ili'l'} \left[ \textbf{u}^{j}_{i'} \right] \partial_{ll'} d^0 \right)_{i=1,\ldots,d} - \left( \mathcal{C}_{ili'l'} \left[ \textbf{u}^{0}_{i'} \right] \partial_{ll'} d^{j} \right)_{i=1,\ldots,d} \nonumber \\
&& - \left( \mathcal{C} \left[ \nabla \textbf{u}^{j} \right] \right) \nabla d^0 - \left( \mathcal{C} \left[ \nabla \textbf{u}^{0} \right] \right) \nabla d^{j} + \left[ c^{j} \right] \left(\mathcal{C} \mathcal{E}^\star \right) \nabla d^0 \nonumber \\
&& + \left[ c^{0} \right] \left(\mathcal{C} \mathcal{E}^\star \right) \nabla d^{j} + \textbf{j}^{j} + M \textbf{l}^{j} + M \textbf{K}^0 d^{j} \nonumber \\
& = & - \left( \mathcal{C}_{ili'l'} \left( \partial_l (\textbf{u}^{+}_j)_{i'} - \partial_l (\textbf{u}^{-}_j)_{i'} \right) \partial_{l'} d^0 \right)_{i=1,\ldots,d} \nonumber \\
&&- \left( \mathcal{C}_{ili'l'} \left( \partial_l (\textbf{u}^{+}_0)_{i'} - \partial_l (\textbf{u}^{-}_0)_{i'} \right) \partial_{l'} d^{j} \right)_{i=1,\ldots,d} \nonumber \\
&& - d^j \left( \mathcal{C}_{ili'l'} \textbf{l}^{0}_{i'} \partial_{ll'} d^0 \right)_{i=1,\ldots,d} - \left( \mathcal{C}_{ili'l'} \left((\textbf{u}^+_{\ast(j-1)})_{i'} - (\textbf{u}^-_{\ast(j-1)})_{i'} \right) \partial_{ll'} d^0 \right)_{i=1,\ldots,d} \nonumber\\
&&- \left( \mathcal{C} \left( \nabla \textbf{u}^{+}_j - \nabla \textbf{u}^{-}_j \right) \right) \nabla d^0 + \left(\mathcal{C} \left( \textbf{l}^0 \otimes \nabla d^j \right)\right) \nabla d^0 + \left( c^+_j - c^-_j \right) \left( \mathcal{C} \mathcal{E}^\star\right) \nabla d^0 \nonumber \\
&& + 2 \left(\mathcal{C} \mathcal{E}^\star \right) \nabla d^{j} - M \left( \partial_\nu (\textbf{l}^0 d^j) \right) + M \left( \partial_\nu \textbf{u}^+_j - \partial_\nu \textbf{u}^-_j \right) \nonumber \\
&& - M \left( \partial_\nu \textbf{u}^+_{\ast(j-1)} - \partial_\nu \textbf{u}^-_{\ast(j-1)} \right) + M \textbf{K}^0 d^j \label{compuj} \, ,
\end{eqnarray}
where $\nu = \nabla d^0$ is the unit outward normal of $\Gamma^0_t$ and since $\left[ \textbf{u}^j \right] = \textbf{l}^0 d^j + \textbf{u}^+_{\ast(j-1)} - \textbf{u}^-_{\ast(j-1)}$ on $\Gamma^0$, $\textbf{u}^+_0 = \textbf{u}^-_0$ on $\Gamma^0$, and
\[
- \left( \mathcal{C} \left[ \nabla \textbf{u}^{0} \right] \right) \nabla d^{j} + \textbf{j}^j = \left(\mathcal{C} \left( \textbf{l}^0 \otimes \nabla d^j \right)\right) \nabla d^0 
\]
due to (\ref{jk}), (\ref{l0}), and (\ref{randablu}). To simplify this equation for $(x,t) \in \Gamma^0$, we use an analogous calculation as in Step 4 in Subsection \ref{subsec0th}. Let $\tau_1, \ldots \tau_{d-1}$ be an orthonormal basis of the tangent space of $\Gamma^0_t$. Then equation (\ref{sprunguj}) yields for all $i=1, \ldots , d-1$
\begin{eqnarray*}
\left( \partial_{\tau_i} \textbf{u}^+_j \otimes \tau_i \right) - \left( \partial_{\tau_i} \textbf{u}^-_j \otimes \tau_i \right) = \partial_{\tau_i} \left( \textbf{l}^0 d^j + \textbf{u}^+_{\ast(j-1)} - \textbf{u}^-_{\ast(j-1)} \right) \otimes \tau_i 
\end{eqnarray*}
and therefore by equation (\ref{gradu}) we get
\begin{eqnarray}
\nabla \textbf{u}^+_j - \nabla \textbf{u}^-_j & = & \left( \partial_\nu \textbf{u}^+_j \otimes \nu \right) - \left( \partial_\nu \textbf{u}^-_j \otimes \nu \right) \nonumber \\
&& + \sum^{d-1}_{i=1}{\partial_{\tau_i} \left( \textbf{l}^0 d^j + \textbf{u}^+_{\ast(j-1)} - \textbf{u}^-_{\ast(j-1)} \right) \otimes \tau_i} \,. \label{sprunggraduj}
\end{eqnarray}
We define the matrix $B(d^j) \in \mathbb{R}^{d \times d}$ by
\begin{eqnarray*}
B(d^j):=\sum^{d-1}_{i=1}{\partial_{\tau_i} \left( \textbf{l}^0 d^j + \textbf{u}^+_{\ast(j-1)} - \textbf{u}^-_{\ast(j-1)} \right) \otimes \tau_i} \, .
\end{eqnarray*}
Using (\ref{sprunggraduj}) and $\nabla d^0 = \nu$ on $\Gamma^0$, we obtain on $\Gamma^0$
\begin{eqnarray*}
\lefteqn{ \left( \mathcal{C}_{ili'l'} \left( \partial_l (\textbf{u}^{+}_j)_{i'} - \partial_l (\textbf{u}^{-}_j)_{i'} \right) \partial_{l'} d^0 \right)_{i=1,\ldots,d}} \\
& = & \left( \mathcal{C} \left( \left( \partial_\nu \textbf{u}^+_j \otimes \nu \right) - \left( \partial_\nu \textbf{u}^-_j \otimes \nu \right) \right)\right) \nu + \left( \mathcal{C}_{ili'l'} B(d^j)_{i'l} \partial_{l'} d^0 \right)_{i=1,\ldots,d} \\
& = & \left( \mathcal{C} \left( \nabla \textbf{u}^+_j - \nabla \textbf{u}^-_j \right) \right) \nu + \left( \mathcal{C}_{ili'l'} B(d^j)_{i'l} \partial_{l'} d^0 \right)_{i=1,\ldots,d} - \left( \mathcal{C} B(d^j) \right) \nu \,. 
\end{eqnarray*}
Due to the definition of $M$ (see (\ref{defM})) and (\ref{sprunggraduj}), we get on $\Gamma^0$ 
\begin{eqnarray*}
M \left( \partial_\nu \textbf{u}^+_j - \partial_\nu \textbf{u}^-_j \right) & = & \left( \mathcal{C}_{ili'l'} \left( \partial_\nu (\textbf{u}^+_j)_{i'} - \partial_\nu ( \textbf{u}^{-}_j )_{i'} \right) \partial_l d^0 \partial_{l'} d^0 \right)_{i=1,\ldots,d} \\
& = & \left( \mathcal{C} \left( \nabla \textbf{u}^+_j - \nabla \textbf{u}^-_j \right) \right) \nu - \left( \mathcal{C} B(d^j) \right) \nu \,.
\end{eqnarray*}
So we obtain
\begin{eqnarray}
\lefteqn{\left( \mathcal{C}_{ili'l'} \left( \partial_l (\textbf{u}^{+}_j)_{i'} - \partial_l (\textbf{u}^{-}_j)_{i'} \right) \partial_{l'} d^0 \right)_{i=1,\ldots,d}} \nonumber \\
\lefteqn{ + \left( \mathcal{C} \left( \nabla \textbf{u}^{+}_j - \nabla \textbf{u}^{-}_j \right) \right) \nabla d^0 - M \left( \partial_\nu \textbf{u}^+_j - \partial_\nu \textbf{u}^-_j \right)} \nonumber \\
& = & \left( \mathcal{C} \left( \nabla \textbf{u}^+_j - \nabla \textbf{u}^-_j \right) \right) \nu + \left( \mathcal{C}_{ili'l'} B(d^j)_{i'l} \partial_{l'} d^0 \right)_{i=1,\ldots,d} \label{sprunggraduj2} \,.
\end{eqnarray}
Furthermore, we use the definition of $B(d^j)$ and $M$ (see (\ref{defM})) to obtain
\begin{eqnarray}
\lefteqn{\left( \mathcal{C}_{ili'l'} B(d^j)_{i'l} \partial_{l'} d^0 \right)_{i=1,\ldots,d} + M \left( \partial_\nu (\textbf{l}^0 d^j) \right) + M \left( \partial_\nu \textbf{u}^+_{\ast(j-1)} - \partial_\nu \textbf{u}^-_{\ast(j-1)} \right)} \nonumber \\
& = & \sum^{d-1}_{k=1}{\partial_{\tau_k} d^j \left( \mathcal{C}_{ili'l'} ( \textbf{l}^0 \otimes \tau_k )_{i'l} \partial_{l'} d^0 \right)_{i=1,\ldots,d}} + \partial_{\nu} d^j \left( \mathcal{C}_{ili'l'} ( \textbf{l}^0 \otimes \nabla d^0 )_{i'l} \partial_{l'} d^0 \right)_{i=1,\ldots,d} \nonumber \\
&& + \sum^{d-1}_{k=1}{ d^j \left( \mathcal{C}_{ili'l'} ( \partial_{\tau_k} \textbf{l}^0 \otimes \tau_k )_{i'l} \partial_{l'} d^0 \right)_{i=1,\ldots,d}} + d^j \left( \mathcal{C}_{ili'l'} ( \partial_{\nu} \textbf{l}^0 \otimes \nabla d^0 )_{i'l} \partial_{l'} d^0 \right)_{i=1,\ldots,d} \nonumber \\
&& + \sum^{d-1}_{k=1}{ \left( \mathcal{C}_{ili'l'} ( \partial_{\tau_k} (\textbf{u}^+_{\ast(j-1)} - \textbf{u}^-_{\ast(j-1)}) \otimes \tau_k )_{i'l} \partial_{l'} d^0 \right)_{i=1,\ldots,d}} \nonumber \\
&& + \left( \mathcal{C}_{ili'l'} ( \partial_{\nu} (\textbf{u}^+_{\ast(j-1)} - \textbf{u}^-_{\ast(j-1)}) \otimes \nabla d^0 )_{i'l} \partial_{l'} d^0 \right)_{i=1,\ldots,d} \nonumber \\
& = & \left( \mathcal{C} \left( \textbf{l}^0 \otimes \nabla d^0 \right) \right) \nabla d^j + d^j \left( \mathcal{C}_{ili'l'} ( \nabla \textbf{l}^0 )_{i'l} \partial_{l'} d^0 \right)_{i=1,\ldots,d} \nonumber \\
&& + \left( \mathcal{C}_{ili'l'} \left( \nabla (\textbf{u}^+_{\ast(j-1)} - \textbf{u}^-_{\ast(j-1)}) \right)_{i'l} \partial_{l'} d^0 \right)_{i=1,\ldots,d}  \label{sprunggraduj3}
\end{eqnarray}
since $\nu = \nabla d^0$ on $\Gamma^0$. Altogether, the compatibility condition (\ref{compuj}) turns on $\Gamma^0$ with (\ref{sprunggraduj2}) and (\ref{sprunggraduj3})
into 
\begin{eqnarray*}
\left( \mathcal{C} \left( \nabla \textbf{u}^+_j - \nabla \textbf{u}^-_j \right) \right) \nu = B_{j-1} \nabla d^j + \textbf{b}_{j-1} d^j + \textbf{c}_{j-1} \quad \mbox{on } \Gamma^0 ,
\end{eqnarray*}
where $B_{j-1} = B_{j-1}(x,t) \in \mathbb{R}^{d \times d}$, $\textbf{b}_{j-1} = \textbf{b}_{j-1}(x,t) \in \mathbb{R}^d$, and $\textbf{c}_{j-1} = \textbf{c}_{j-1}(x,t) \in \mathbb{R}^d$ only depend on the known functions $\mathcal{V}^0, \ldots , \mathcal{V}^{j-1}$.\\
\textbf{Step 5:} Let $\textbf{F}^{j-1}$ be defined as in (\ref{randbedu}). Then equation (\ref{linHele7}) gives us a boundary condition for the outer expansion $\textbf{u}^+_j$ on $\partial_T \Omega$. Together with the outer expansion (\ref{linHele4}) and the conditions on $\Gamma^0$, the functions $\textbf{u}^\pm_j $ have to solve the following boundary value problem for each $t>0$
\begin{align*}
\Div \left(\mathcal{C} \nabla \textbf{u}^\pm_j\right) & = \Div \left( \mathcal{C} \mathcal{E}^\star c^\pm _j \right) && \mbox{in } Q^\pm_0\, , \\
\left[ (\mathcal{C} \nabla \textbf{u}^\pm_j) \nu \right]_{\Gamma^0_t} & = B_{j-1} \nabla d^j + \textbf{b}_{j-1} d^j + \textbf{c}_{j-1} && \mbox{on }\Gamma^0_t,t>0 \,, \\
\left[ \textbf{u}^\pm_j \right]_{\Gamma^0_t} & = \textbf{l}^0 d^j + \big[ \textbf{u}^\pm_{\ast(j-1)} \big]_{\Gamma^0_t} && \mbox{on }\Gamma^0_t, t>0\,,\\
 \textbf{u}^+_j & = \textbf{F}^{j-1} && \mbox{on } \partial \Omega, t>0 \, . 
\end{align*}
\textbf{Step 6:} We consider the compatibility condition (\ref{loes1}) for $k=j+1$. Since $\int_{\mathbb{R}}{(\eta-1/2) \theta'_0 \, dz} = \int_\mathbb{R}{\textbf{u}^{j-1}_\ast \theta'_0 \, dz} = 0$, we get due to the definition of $\overline{\textbf{u}}_j$ in Lemma \ref{loesle2}, (\ref{uj}), and (\ref{intuast}) in $\Gamma^0(\delta)$
\begin{eqnarray*}
\overline{\textbf{u}}^j(x,t) & = & \frac{1}{2} \int_{\mathbb{R}}{\textbf{u}^j(z) \theta'_0 \, dz} = \tilde{\textbf{u}}^j =  \frac{1}{2} \left( \textbf{u}^+_j + \textbf{u}^-_j \right) - \frac{1}{2} \left( \textbf{u}^+_{\ast(j-1)} + \textbf{u}^-_{\ast(j-1)} \right) \,,
\end{eqnarray*}
where the last equation follows from the definition of $\tilde{\textbf{u}}^j$ in Step 3. Then the compatibility condition (\ref{loes1}) for $k=j+1$ reads on $\Gamma^0$
\begin{align}
\overline{\mu}^j(x,t) = - \sigma \Delta d^j - \frac{1}{2} \mathcal{E}^\star : \mathcal{C} \left( \nabla \textbf{u}^+_j + \nabla \textbf{u}^-_j \right) + \eta_0 d^j g^0 - \mathfrak{A}^{j-1} \, ,\label{j-ko1}
\end{align}
where we set for $(x,t) \in \Gamma^0(\delta)$
\begin{eqnarray*}
\mathfrak{A}^{j-1}(x,t) & = & \widetilde{\mathcal{A}}^{j-1}(x,t) - \frac{1}{2} \mathcal{E}^\star : \mathcal{C} \left( \nabla \textbf{u}^+_{\ast(j-1)} + \nabla \textbf{u}^-_{\ast(j-1)} \right)\!(x,t) \, .
\end{eqnarray*}
Note that $\mathfrak{A}^{j-1}$ only depends on the known functions $\mathcal{V}^0,\ldots,\mathcal{V}^{j-1}$.\\
\textbf{Step 7:} Lemma \ref{loesle4} gives us an equation for $\mu^j$ in $\mathbb{R} \times \Gamma^0(\delta)$ 
\begin{eqnarray}
\mu^j(z,x,t) = \tilde{\mu}^j(x,t) + \left(d^0 h^j+ d^j h^0\right)\!(x,t) \left(\eta(z) - \tfrac{1}{2}\right) + \mu^{j-1}_\ast(z,x,t) \,,\label{j-ko2}
\end{eqnarray}
where $\mu^{j-1}_\ast$ only depends on $\mathcal{V}^0,\ldots,\mathcal{V}^{j-1}$ and satisfies (\ref{loes5}) with $k=j$. We define $\tilde{\mu}^j$ and $h^j$ later. As in the construction of the zero order functions, it follows 
\[
\tilde{\mu}^j(x,t)= \frac{1}{2} \int_\mathbb{R}{\mu^j(z) \theta_0'(z) \, dz} = \overline{\mu}^j(x,t) \quad \forall (x,t) \in \Gamma^0(\delta) \,.
\]
The restriction of (\ref{j-ko2}) on $\Gamma^0$ and equation (\ref{j-ko1}) give us
\begin{eqnarray}
\mu^j(z,x,t)&=& - \sigma \Delta d^j - \frac{1}{2} \mathcal{E}^\star : \mathcal{C} \left( \nabla \textbf{u}^+_j + \nabla \textbf{u}^-_j \right) + d^j\left( \eta_0 g^0 + h^0 \left(\eta-\tfrac{1}{2}\right)\right)\nonumber\\
&& - \mathfrak{A}^{j-1} + \mu^{j-1}_\ast \quad \forall (x,t) \in \Gamma^0. \label{j-ko3}
\end{eqnarray}
So $\mu^j$ is uniquely determined on $\Gamma^0$.\\
\textbf{Step 8:} We consider $z \to \pm \infty$ in (\ref{j-ko3}) and use the inner-outer matching condition to get on $\Gamma^0$
\begin{eqnarray}
\left.\mu^\pm_j\right|_{\Gamma^0} = \lim_{z \to \pm \infty} \mu^j(z,.) & = & - \sigma \Delta d^j - \tfrac{1}{2} \mathcal{E}^\star : \mathcal{C} \left( \nabla \textbf{u}^+_j + \nabla \textbf{u}^-_j \right) \nonumber \\
&& + d^j \left(\eta_0 g^0 \pm \tfrac{1}{2} h^0 \right) - \mathfrak{A}^{j-1} + \mu^{\pm}_{\ast(j-1)},\label{j-ko4} 
\end{eqnarray}
where $\eta (z)=1$ and $\eta(-z)=0$ for $z>1$. Then the function $\mu^\pm_j$ is uniquely determined by the outer expansion (\ref{linHele1}) and the boundary condition in Lemma \ref{lemmamuB}, that is, $\mu^\pm_j$ is the solution to the elliptic boundary problem 
\begin{align*}
\Delta \mu^\pm_j & = \partial_t c^\pm_j && \mbox{in }Q^\pm_0 \,, \\
\mu^\pm_j & = - \sigma \Delta d^j - \tfrac{1}{2} \mathcal{E}^\star : \mathcal{C} \left( \nabla \textbf{u}^+_j + \nabla \textbf{u}^-_j \right) &&\\
& \quad + d^j \left(-\eta_0 g^0 \pm \tfrac{1}{2} h^0 \right) - \mathfrak{A}^{j-1} + \mu^{\pm}_{\ast(j-1)} && \mbox{on }\Gamma^0_t, t\geq 0 \, , \\
\tfrac{\partial}{\partial n } \mu^+_j & = G^{j-1} && \mbox{on }\partial \Omega, t\geq 0 \,. 
\end{align*}
\textbf{Step 9:} Sending $z$ in (\ref{j-ko2}) to $\pm \infty$ and using the matching condition yields
\[
\mu^\pm_j(x,t) = \overline{\mu}^j(x,t) \pm \frac{1}{2} \left(d^0 h^j + d^j h^0\right) (x,t) + \mu^\pm_{\ast(j-1)}(x,t) \quad \forall (x,t) \in \Gamma^0(\delta) \,.
\]
Hence it is necessary and sufficient to take $h^j$ and $\tilde{\mu}^j $ as
\begin{align*}
\overline{\mu}^j = \tilde{\mu}^j &:= \frac{1}{2} \left(\mu_j^+ + \mu^-_j - \mu^+_{\ast(j-1)} - \mu^-_{\ast(j-1)}\right) \quad \mbox{in } \Gamma^0(\delta),\\
h^j &:= \left\{\begin{array}{l @{\quad}l }
\frac{1}{d^0} \left(- d^j h^0 + \mu_j^+ - \mu^-_j - \mu^+_{\ast(j-1)} + \mu^-_{\ast(j-1)}\right) & \mbox{in } \Gamma^0(\delta)\backslash \Gamma^0\\
\nabla d^0 \cdot \nabla \left(- d^j h^0+ \mu_j^+ - \mu^-_j - \mu^+_{\ast(j-1)} + \mu^-_{\ast(j-1)}\right) & \mbox{on } \Gamma^0.
\end{array}\right.
\end{align*}
Consider (\ref{j-ko4}) to recognize that the numerator in the definition of $h^j$ vanishes on $\Gamma^0$. So the definition of $h^j$ is natural. Note that this definition of $\overline{\mu}^j$ coincides with (\ref{j-ko1}). To see this use (\ref{j-ko4}) again.\\
By Steps 2-9 we can determine $\mu^j, \mu^\pm_j, \textbf{u}^j, \textbf{u}^\pm_j, h^j$, and $\textbf{l}^j$ depending on $d^j$. The next step shows us how we can determine $d^j$.\\
\textbf{Step 10:} We consider the compatibility condition (\ref{loes2}) on $\Gamma^0$ for $k= j+1$. Note that $[\mu^0]=0$ and $d^0=0$ and use the definition of $h^j$ on $\Gamma^0$. Then we have
\begin{eqnarray*}
 d^j_t & = & \frac{1}{2} \Delta d^0\left[\mu^j\right] + \nabla d^0 \cdot \left[\nabla \mu^j\right] + \nabla d^j \cdot \left[\nabla \mu^0\right]\\
&& - \frac{1}{2} \nabla d^0 \cdot \nabla \left(- d^j h^0+ \mu_j^+ - \mu^-_j - \mu^+_{\ast(j-1)} + \mu^-_{\ast(j-1)}\right) - \frac{1}{2} d^j L^0 + \widetilde{\mathcal{B}}^{j-1} \\
& = & \frac{1}{2} \left(\Delta d^0 h^0 + \nabla d^0 \cdot \nabla h^0 - L^0\right) d^j + \frac{1}{2} \left(\tfrac{\partial }{\partial \nu} \mu^+_j - \tfrac{\partial }{\partial \nu} \mu^-_j\right) + d^{j-1}_{\ast},
\end{eqnarray*}
where $\nu= \nabla d^0$ is the unit outward normal of $\Gamma^0_t$ , $\left[\mu^j\right] = h^0 d^j + \mu^+_{\ast(j-1)}-\mu^-_{\ast(j-1)}$, $\nabla d^0 \cdot \left[\nabla \mu^j\right] = \frac{\partial }{\partial n} \mu^+_j - \frac{\partial }{\partial n} \mu^-_j$, $\nabla d^j \cdot \left[\nabla \mu^0\right] = \nabla d^j \cdot \nabla d^0 h^0$ on $\Gamma^0$, $\nabla d^j \cdot \nabla d^0 = - \tfrac{1}{2}\sum_{i=1}^{j-1}{\nabla d^i \cdot \nabla d^{j-i}}$, and 
\begin{eqnarray*}
d^{j-1}_\ast(x,t) & = & \frac{1}{2} \Delta d^0 \left( \mu^+_{\ast (j-1)} + \mu^-_{\ast (j-1)} \right) - \frac{1}{2} \nabla d^0 \cdot \nabla \left( \mu^+_{\ast (j-1)} - \mu^-_{\ast (j-1)} \right) \\
&& - \frac{3}{4} \sum^{j-1}_{i=1}{\nabla d^i \cdot \nabla d^{j-i}} h^0 + \widetilde{\mathcal{B}}^{j-1} \,.
\end{eqnarray*} 
Here we can see that $d^{j-1}_\ast$ only depends on the known functions $\mathcal{V}^0,\ldots,\mathcal{V}^{j-1}$.\\
Now we define the functions $a^i_{j-1}$, $i=1, \ldots, 12$, used in Subsection \ref{secapp}. 
\begin{align}
a^1_{j-1} & := E^k + A^{k-1} \mbox{ (Step 1)}, & a^{2\pm}_{j-1} & := \eta_0 g^0 \pm \tfrac{1}{2}h^0  \mbox{ (Step 8)}, \\
a^{3 \pm}_{j-1} & := - \mathfrak{A}^{j-1} + \mu^{\pm}_{\ast(j-1)} \mbox{ (Step 8)}, & a^4_{j-1} & := B_{j-1} \mbox{ (Step 5)}, \\
a^5_{j-1} & := \textbf{b}_{j-1} \mbox{ (Step 5)}, & a^6_{j-1} & := \textbf{c}_{j-1} \mbox{ (Step 5)}, \\
a^7_{j-1} & := \textbf{l}^0 \mbox{ (Step 5)}, & a^8_{j-1} & := [ \textbf{u}^\pm_{\ast (j-1)} ]_{\Gamma^0_t}  \mbox{ (Step 5)}, \\
a^9_{j-1} & := \tfrac{1}{2} (\Delta d^0 h^0 + \nabla d^0 \cdot \nabla h^0 - L^0) &&  \!\!\!\!\!\!\!\!\!\!\!\!\!\!\!\!\!\!\!\!\!\!\!\!\!\!\mbox{(Step 10)}, \\
a^{10}_{j-1} & := d^{j-1}_\ast \mbox{ (Step 10)},\\
a^{11}_{j-1} & := - \mu^+_{\ast(j-1)} \eta - \mu^-_{\ast(j-1)} (1-\eta) + \mu^{j-1}_\ast && \!\!\!\!\!\!\!\!\!\!\!\!\! \mbox{(Step 7+9)}, \\
a^{12}_{j-1} & := - \textbf{u}^+_{\ast(j-1)} \eta - \textbf{u}^-_{\ast(j-1)} (1-\eta) + \textbf{u}^{j-1}_\ast && \!\!\!\!\!\!\!\!\!\!\!\!\! \mbox{(Step 2+3)}.
\end{align}
Altogether, $(d^j,\mu^\pm_j, \textbf{u}^\pm_j)$ satisfies the problem (\ref{linHele1})-(\ref{linHele10}). Assume that $(d^j, \mu^\pm_j, \textbf{u}^\pm_j)$ solve problem (\ref{linHele1})-(\ref{linHele10}). In Lemma \ref{linHeleShawle} below we verify that $\mu^\pm_j, \mu^j, \textbf{u}^\pm_j$, and $ \textbf{u}^j $ satisfy the matching condition (\ref{match1}) and the compatibility conditions (\ref{compu}) and (\ref{loes2}) for $k=j+1$. \\
Notice that in the derivation of (\ref{linHele1})-(\ref{linHele10}) we need the inner expansions only for $(x,t) \in \Gamma^0$ where $O^+_{j-1} = O^-_{j-1} = 0$ and $\textbf{P}^+_{j-1} = \textbf{P}^-_{j-1} = 0$ by the definitions of $c^\pm_{j-1}, \mu^\pm_{j-1} ,\textbf{u}^\pm_{j-1}, O^\pm_{j-1}$, and $\textbf{P}^\pm_{j-1}$. Therefore the solution $d^1$ of (\ref{linHele1})-(\ref{linHele10}) is independent of the terms $\epsilon^2 (O^+_0 \eta^+_N + O^-_0 \eta_N^-)$ and $\epsilon^2(\textbf{P}^+_0 \eta^+_N + \textbf{P}^-_0 \eta_N^-)$ which we added in (\ref{modinnere2}) and (\ref{modinnere3}). In particular, $d^1$ is independent of the constant $N$. So we can define $N := \left\| d^1\right\|_{C^0(\Gamma^0(\delta))} + 2$, see Remark \ref{remarkN}.\\
\textbf{Step 11:} We can define $g^j$ in $\Gamma^0(\delta) \backslash \Gamma^0$ in a unique way such that the compatibility condition (\ref{loes1}) is satisfied for $k=j+1$. Since (\ref{loes1}) already holds on $\Gamma^0$, we can extend $g^j$ smoothly to $\Gamma^0$. Similarly, we can define uniquely $\textbf{K}^j$ and $L^j$ to satisfy the compatibility conditions (\ref{compu}) and (\ref{loes2}).\\
\textbf{Step 12:} By Lemma \ref{lemmauB} and Lemma \ref{lemmamuB} we immediately get $\textbf{u}^j_B$ and $\mu^j_B$.\\
After going through Step 1-12, we obtain.

\begin{lemma} \label{linHeleShawle}
Let $j \geq 1$ be an arbitrary integer and assume $\mathcal{V}^0,\ldots,\mathcal{V}^{j-1}$ are given and satisfy the matching conditions (\ref{match1}) and (\ref{matchbound1}) for all $k=0, \ldots, j-1$. Furthermore, let the compatibility conditions (\ref{compu1}), (\ref{loes1a}), and (\ref{loes2}) if $j=1$ or the compatibility conditions (\ref{compu}), (\ref{loes1}), and (\ref{loes2}) for $k=j$ if $j>1$ be satisfied. Assume that the problem (\ref{linHele1})-(\ref{linHele10}) admits a smooth solution in the time interval $[0,T]$. \\
Define $\textbf{l}^j$ as in Step 3, $h^j$ as in Step 9, and $(g^j,L^j,\textbf{K}^j)$ as in Step 11. Then, for $k=j$, $(\mu^\pm_j, \textbf{u}^\pm_j, d^j)$ defined as the solution to (\ref{linHele1})-(\ref{linHele10}) and $(c^\pm_j,c^j,c^j_B,\mu^j,\mu^j_B,\textbf{u}^j,\textbf{u}^j_B)$ defined by (\ref{defcpmj})-(\ref{systemckB}) and (\ref{Defmuj})-(\ref{DefuB}) satisfy the outer expansions equations (\ref{defcpmj}),(\ref{linHele1}), and (\ref{linHele4}), the inner expansions equations (\ref{innere3})-(\ref{innere2}), the boundary-layer expansion equations (\ref{boundODE1})-(\ref{boundODE3}) and (\ref{boundinitial1})-(\ref{boundinitial3}), the inner-outer matching condition (\ref{match1}), the outer-boundary matching condition (\ref{matchbound1}), and (\ref{djbedingung}). In addition, the compatibility conditions (\ref{compu}), (\ref{loes1}), and (\ref{loes2}) for $k=j+1$ are also satisfied.
\end{lemma}

\proof  We verify the required conditions by direct calculation. Define $(\tilde{\textbf{u}}^j, \tilde{\mu}^j)$ as in Step 3 and Step 9. Then we can verify that $\textbf{u}^j$ satisfies (\ref{uj}), $\mu^j$ satisfies (\ref{j-ko2}), $\tilde{\mu}^j = \overline{\mu}^j$, and $\tilde{\textbf{u}}^j = \overline{\textbf{u}}^j$ where $\overline{\mu}^j$ and $\overline{\textbf{u}}^j$ are defined as in Lemma \ref{loesle2}. By the interface condition (\ref{j-ko4}) we conclude that the identity for $\overline{\mu}^j$ in Step 6 coincides with the definition in Step 9.\\
\textbf{To (\ref{defcpmj}), (\ref{linHele1}), (\ref{linHele4}):} The outer expansions equations are satisfied by definition of $(c^\pm_j, \mu^\pm_j, \textbf{u}^\pm_j)$.\\
\textbf{To (\ref{innere3})-(\ref{innere2}):} The inner expansions equations are satisfied by definition of $c^j$ and by Lemma \ref{loesle5} and \ref{loesle4} since $\textbf{u}^j$ and $\mu^j$ satisfy (\ref{uj}) and (\ref{j-ko2}).\\
\textbf{To (\ref{boundODE1})-(\ref{boundODE3}) and (\ref{boundinitial1})-(\ref{boundinitial3}):} The assertions immediately follow from the definitions of $(c^j_B, \mu^j_B, \textbf{u}^j_B)$ and the results of Subsection \ref{boundaryexp}.\\
\textbf{To (\ref{match1}):} Due to Lemma \ref{loesle2} the inner-outer matching condition is satisfied for $c^j$ and $c^\pm_j$. We use $\eta(\pm z)- \frac{1}{2} = \pm \frac{1}{2}$ for $z >1$ and 
\begin{equation*}
D_x^m D_t^n D_z^l \left[(\mu^{j-1}_\ast,\textbf{u}^{j-1}_\ast)(\pm z) - (\mu_{\ast(j-1)}^\pm, \textbf{u}_{\ast(j-1)}^\pm) \right] = \mathcal{O}(e^{- \alpha z}) \quad \mbox{as } z \to \infty
\end{equation*}
for all $m,n,l \geq 0$ and $(x,t) \in \Gamma^0(\delta)$ due to Lemma \ref{loesle4}. Then the inner-outer matching conditions directly follows from the definitions of $\mu^j$ and $\textbf{u}^j$. \\
\textbf{To (\ref{matchbound1}):} Due to Lemma \ref{lemmauB} - \ref{lemmamuB}, the assertion hold.\\
\textbf{To (\ref{djbedingung}):} The equation is satisfied by the definition of $d^j$.\\
\textbf{To (\ref{compu}):} By definition of $\textbf{K}^j$ the compatibility condition holds in $\Gamma^0(\delta) \backslash \Gamma^0$. On $\Gamma^0$ the assertion follows from the interface condition $\left[ \left( \mathcal{C} \nabla \textbf{u}^\pm_j\right) \nu \right]_{\Gamma^0_t} = B_{j-1} \nabla d^j + \textbf{b}_{j-1}$ and the same calculation as in Step 4.\\
\textbf{To (\ref{loes1}):} In $\Gamma^0(\delta) \backslash \Gamma^0$ we satisfy the compatibility condition by definition of $g^j$. On $\Gamma^0$ the compatibility condition is satisfied by the definition of $\tilde{\mu}^j = \overline{\mu}^j$ in Step 9, the interface condition (\ref{j-ko4}) for $\mu^\pm_j$, and the definition of $\mathfrak{A}^{j-1}$ in Step 6. Here note that $\overline{\textbf{u}}^j = \tilde{\textbf{u}}^j$ (see Step 6).\\
\textbf{To (\ref{loes2}):} In $\Gamma^0(\delta) \backslash \Gamma^0$ we satisfy the compatibility condition by definition of $L^j$. On $\Gamma^0$ the compatibility condition is satisfied by the interface condition (\ref{linHele8}). Also we apply the inner-outer matching condition (\ref{match1}).  
\makebox[1cm]{} \hfill $\Box$\\

As consequence we obtain recursively:

\begin{theorem} \label{theoremexpansion}
Let $(\mu,\textbf{u}, \Gamma)$ be a smooth solution for the Hele-Shaw problem (\ref{sharpsystem1})-(\ref{sharpsystem7}). Then, for any fixed integer $K>0$, there exist $\mathcal{V}^0, \ldots , \mathcal{V}^K$ such that the outer expansions equations in Subsection \ref{secouterexp}, the inner expansions equations (\ref{innere3})-(\ref{innere2}), the boundary-layer expansion equations (\ref{boundODE1})-(\ref{boundODE3}) and (\ref{boundinitial1})-(\ref{boundinitial3}), the inner-outer matching condition (\ref{match1}), and the outer-boundary matching condition (\ref{matchbound1}) are satisfied for $k = 0, \ldots, K$. In addition, $(\mu^\pm_0, \textbf{u}^\pm_0, \Gamma^0)$ coincides with $(\mu,\textbf{u}, \Gamma)$.
\end{theorem}

\subsection{Proof of Theorem \ref{theoremcA-c0}} \label{subsecproof}

An important question is how good our approximate solution $(c^K_A,\mu^K_A,\textbf{u}^K_A)$ is. To this end we consider $(c^K_I,\mu^K_I,\textbf{u}^K_I)$, $(c^K_O, \mu^K_O,\textbf{u}^K_O)$, and $(c^K_\partial, \mu^K_\partial,\textbf{u}^K_\partial)$ separately. Notice that $\left| \frac{d^K_\epsilon}{\epsilon} - \frac{d^0 + \epsilon d^1}{\epsilon} \right| = \left| \sum^K_{i=2}{\epsilon^{i-1} d^i} \right| \leq 1$ for all $\epsilon$ small enough, and hence by Remark \ref{remarkN}, $\left.O^+_j \eta^+_N + O^-_j \eta^-_N\right|_{z=d^K_\epsilon/\epsilon} = 0$ and $\left. \textbf{P}^+_j \eta^+_N + \textbf{P}^-_j \eta^-_N\right|_{z=d^K_\epsilon/\epsilon} = 0$, for $j=0,\ldots,K-2$. By the inner expansion equations (\ref{innere2}) we obtain for all $(x,t) \in \Gamma^0(\delta)$ and $z=d^K_\epsilon/\epsilon$
\begin{eqnarray}
\lefteqn{((c^K_I)_t - \Delta \mu^K_I)(x,t)}\\& = & - \frac{\left|\nabla d^K_\epsilon \right|^2 - 1}{\epsilon^2} \sum^K_{i=0}{\epsilon^i \mu^i_{zz}} + \frac{1}{\epsilon} \sum_{\substack{0\leq i,j \leq K \\ i+j\geq K}}{\epsilon^{i+j} \left( c^i_z d^j_t - 2 \nabla \mu^i_z \cdot \nabla d^j - \mu^i_z \Delta d^j\right)} \nonumber \\
&& + \sum_{i=K-1}^K{ \epsilon^i \left( c^i_t - \Delta \mu^i \right)} - \epsilon^{K-2} h^K d^0 \eta'' + \frac{1}{\epsilon^2} \sum_{\substack{ 0 \leq i \leq K \\ 0 \leq j \leq K-1 \\ i+j \geq K+1}}{ \epsilon^{i+j} d^i h^j \eta''} \nonumber \\
&& - \epsilon^{K-2} L^{K-1} d^0 \eta' + \frac{1}{\epsilon} \sum_{\substack{ 0 \leq i \leq K \\ 0 \leq j \leq K-2 \\ i+j \geq K}}{ \epsilon^{i+j} d^i L^j \eta'} \quad \forall (x,t) \in\Gamma^0(\delta)\,,
\end{eqnarray}
where we have used (\ref{nabladK}). Furthermore, we have by the inner expansion equations (\ref{innere1}) for all $(x,t) \in \Gamma^0(\delta)$ and $z=d^K_\epsilon/\epsilon$
\begin{eqnarray}
\lefteqn{( \mu^K_I + \epsilon \Delta c^K_I - \epsilon^{-1} f(c^K_I) - W_{,c}(c^K_I,\mathcal{E}(\textbf{u}^K_I)))(x,t)} \nonumber \\
& = & \epsilon^K \mu^K + \epsilon \sum^K_{i=K-1}{\epsilon^i \Delta c^i} - \frac{1- \left|\nabla d^K_\epsilon \right|^2}{\epsilon} \sum^K_{i = 0}{ \epsilon^i c^i_{zz} } \nonumber \\
&& + \sum_{\substack{0 \leq i,j \leq K \\ i+j \geq K}}{\epsilon^{i+j} \left( 2 \nabla c^i_z \cdot \nabla d^j + c^i_z \Delta d^j \right)} - \epsilon^K f^K(c^0,\ldots,c^K) \nonumber \\
&& + \sum_{\substack{0 \leq i,j \leq K \\ i+j\geq K+1}}{\epsilon^{i+j-1} \mathcal{E}^\star : \mathcal{C} (\textbf{u}^i_z \otimes \nabla d^j) } + \epsilon^K \mathcal{E}^\star : \mathcal{C} \nabla \textbf{u}^K - \epsilon^K \mathcal{E}^\star : \mathcal{C} \mathcal{E}^\star c^K \nonumber \\
&& + \epsilon^{K-1} g^{K-1} d^0 \eta' - \sum_{\substack{0 \leq i \leq K-2 \\ 0 \leq j \leq K \\ i+j \geq K}}{\epsilon^{i+j} g^i d^j \eta'} + \epsilon^{K-1} k^K d^0 \eta' \nonumber \\
&& - \sum_{\substack{0 \leq i \leq K-1 \\ 0 \leq j \leq K \\ i+j \geq K+1}}{\epsilon^{i+j-1} k^i d^j \eta'} \quad\quad \forall (x,t) \in \Gamma^0(\delta) 
\end{eqnarray}
and by the inner expansion equations (\ref{innere3}) for all $(x,t) \in \Gamma^0(\delta)$ and $z=d^K_\epsilon/\epsilon$
\begin{eqnarray}
\lefteqn{( \Div \left(\mathcal{C} \nabla \textbf{u}^K_I\right) - \Div \left( \mathcal{C} \mathcal{E}^\star c^K_I \right))(x,t)} \nonumber \\
& = & \sum_{\substack{0 \leq i,j,k \leq K \\ i+j+k \geq K+1}}{ \epsilon^{i+j+k-2} (\mathcal{C}(\textbf{u}^i \otimes \nabla d^j)) \nabla d^k}  + \sum_{\substack{0 \leq i,j \leq K \\ i+j \geq K}}{\epsilon^{i+j-1}(\mathcal{C}_{klk'l'} \partial_l(\textbf{u}^i_{k'})_z \partial_{l'}d^j)_{k=1,\ldots,d}} \nonumber \\
&& + \sum_{\substack{0 \leq i,j \leq K \\ i+j \geq K}}{\epsilon^{i+j-1}(\mathcal{C}_{klk'l'} (\textbf{u}^i_{k'})_z \partial_{ll'}d^j)_{k=1,\ldots,d}} + \sum_{\substack{0 \leq i,j \leq K \\ i+j \geq K}}{\epsilon^{i+j-1} (\mathcal{C} \nabla \textbf{u}^i_z) \nabla d^j} \nonumber \\
&& + \sum_{i=K-1}^K{\epsilon^i (\mathcal{C}_{klk'l'} \partial_{ll'} \textbf{u}^i_{k'})_{k=1,\ldots,d}} - \sum_{\substack{0 \leq i,j \leq K \\ i+j \geq K}}{\epsilon^{i+j-1} c^i_z (\mathcal{C} \mathcal{E}^\star) \nabla d^j} \nonumber \\
&& - \sum_{i=K-1}^{K}{\epsilon^i (\mathcal{C} \mathcal{E}^\star) \nabla c^i} - \sum_{\substack{ 0 \leq i \leq K \\ 0 \leq j \leq K-1 \\ i+j \geq K+1}}{\epsilon^{i+j-2} d^i M \textbf{l}^j \eta''} + \epsilon^{K-2} d^0 M \textbf{l}^K \eta'' \nonumber \\
&& - \sum_{\substack{0 \leq i \leq K \\ 0 \leq j \leq K-2 \\ i+j \geq K}}{\epsilon^{i+j-1} d^i M \textbf{K}^j \eta'} + \epsilon^{K-2} d^0 M \textbf{K}^{K-1} \eta' \nonumber \\
&& - \sum_{\substack{0 \leq i \leq K \\ 0 \leq j \leq K-1 \\ i+j \geq K+1}}{\epsilon^{i+j-2} d^i \textbf{j}^j \eta''} + \epsilon^{K-2} d^0 \textbf{j}^K \eta'' \quad \forall (x,t) \in \Gamma^0(\delta) \,.
\end{eqnarray}
By definition we obtain for the outer expansions in $Q^+_0 \cup Q^-_0$
\begin{align}
(c^K_O)_t - \Delta \mu^K_O &= 0 \, , \\
 \mu^K_O + \epsilon \Delta c^K_O - \epsilon^{-1} f(c^K_O) - W_{,c}(c^K_O,\mathcal{E}(\textbf{u}^K_O)) & =  \epsilon^K \mu^\pm_K - \epsilon^K f^K(c^\pm_0,\ldots,c^\pm_K) \nonumber \\
+ \epsilon^K \mathcal{E}^\star : \mathcal{C}(\mathcal{E}&(\textbf{u}^\pm_K) - \mathcal{E}^\star  c^\pm_K) + \sum^K_{i=K-1}{\epsilon^{i+1} \Delta c^\pm_i} \, , \\
\Div \left( \mathcal{C} \mathcal{E}(\textbf{u}^K_O) \right) - \Div \left(\mathcal{C} \mathcal{E}^\star c^K_O\right) & = 0 \, . \label{fehleraussen}
\end{align}
Finally, we get for the boundary-layer expansion in $ \partial\Omega(\delta) \times (0,T)$ and $z=d_B/\epsilon$
\begin{eqnarray}
\lefteqn{((c^K_\partial)_t - \Delta \mu^K_\partial)(x,t)} \nonumber \\
& = & \epsilon^{K-1} \left( c^K_{B,z} d_{B,t} - 2 \nabla \mu^K_{B,z} \cdot \nabla d_B - \mu^K_{B,z} \Delta d_B \right)  + \sum^K_{i=K-1}{\epsilon^i \left( c^i_{B,t} - \Delta \mu^i_B \right)} \nonumber \\
&& - \epsilon^K \left( c^K_{B,t}(0) - \Delta \mu^K_B(0) \right)  \quad \forall (x,t) \in \partial\Omega(\delta) \times (0,T)\, ,
\end{eqnarray}
where $(c^K_B(0),\mu^K_B(0)) = (c^K_B,\mu^K_B)(0,x,t)$ for all $(x,t) \in \overline{ \partial \Omega(\delta) } \times [0,T]$. For the chemical potential equation we get in $ \partial\Omega(\delta) \times [0,T]$ with $z=d_B/\epsilon$ 
\begin{eqnarray}
\lefteqn{ (\mu^K_\partial + \epsilon \Delta c^K_\partial - \epsilon^{-1} f(c^K_\partial) - W_{,c}(c^K_\partial,\mathcal{E}(\textbf{u}^K_\partial)))(x,t)} \nonumber \\
& = & \epsilon^K \mu^K_B + \epsilon^K 2 \nabla c^K_{B,z} \cdot \nabla d_B + \epsilon^K c^K_{B,z} \Delta d_B + \sum^{K+1}_{i=K}{\epsilon^i \Delta c^{i-1}_{B}} \nonumber \\
&& - \epsilon^K f^K(c^0_B,\ldots,c^K_B-c^K_B(0)) + \epsilon^K \mathcal{E}^\star : \mathcal{C} \nabla \textbf{u}^K_B - \epsilon^K \mathcal{E}^\star :\mathcal{C} \mathcal{E}^\star c^K_B    \nonumber \\
&& - \epsilon^K \left( \mu^K_B(0) + \epsilon \Delta \mu^K_B(0) -  \epsilon ^{-1} f'(\theta_0) c^K_B(0) - \mathcal{E}^\star : \mathcal{C} \mathcal{E}^\star c^K_B(0) \right) 
\end{eqnarray}
for all $ (x,t) \in \partial\Omega(\delta) \times [0,T]$. For the equation of the stress tensor we obtain in $ \partial\Omega(\delta) \times [0,T]$ with $z=d_B/\epsilon$
\begin{eqnarray}
\lefteqn{( \Div \left(\mathcal{C} \nabla \textbf{u}^K_\partial \right) - \Div \left( \mathcal{C} \mathcal{E}^\star c^K_\partial \right))(x,t)} \nonumber \\
& = & \epsilon^{K-1} \left( \mathcal{C}_{iji'j'} \partial_j (\textbf{u}^K_{B,i'})_z \partial_{j'} d_B \right)_{i=1,\ldots,d} + \epsilon^{K-1} \left( \mathcal{C}_{iji'j'} (\textbf{u}^K_{B,i'})_z \partial_{jj'} d_B \right)_{i=1,\ldots,d} \nonumber \\
&& + \epsilon^{K-1} \left(\mathcal{C} \nabla \textbf{u}^K_{B,z} \right) \nabla d_B + \sum^K_{i=K-1}{\epsilon^i \Div \left( \mathcal{C} \nabla \textbf{u}^i_B \right)} - \epsilon^{K-1}\left( \mathcal{C} \mathcal{E}^\star \right) \nabla d_B c^K_{B,z} \nonumber \\
&& - \sum^K_{i=K-1}{ \epsilon^i \left( \mathcal{C} \mathcal{E}^\star \right) \nabla c^i_B} + \epsilon^K (\mathcal{C} \mathcal{E}^\star )\nabla c^K_B(0) \quad \forall (x,t) \in \partial\Omega(\delta) \times [0,T]\,.
\end{eqnarray}
Additionally, we check the boundary conditions on $\partial \Omega \times (0,T)$. Consider the boundary conditions (\ref{boundinitial1})-(\ref{boundinitial3}) and note the extra terms $\epsilon^K c^K_B(0,x,t)$ and $\epsilon^K \mu^K_B(0,x,t)$ added in the definitions of $ c^K_\partial(x,t)$ and $ \mu^K_\partial(x,t)$. Then we obtain 
\begin{align*}
\textbf{u}^K_\partial & = 0  & \mbox{and} && \frac{\partial}{\partial n} c^K_\partial &= \frac{\partial}{\partial n} \mu^K_\partial = 0 && \mbox{on } \partial\Omega \times (0,T)\,.
\end{align*}
It remains to show how good the approximate solutions $(c^K_A,\mu^K_A,\textbf{u}^K_A)$ are in the domains $\Gamma^0(\delta) \backslash \Gamma^0(\delta/2)$ and $\partial_T \Omega(\delta) \backslash \partial_T \Omega(\delta/2)$ where we have glued together the inner and outer approximate solutions and the boundary-layer and outer approximate solutions. By definition $\left|d^0(x,t)\right| \in [\delta/2, \delta)$ for $(x,t) \in \Gamma^0(\delta) \backslash \Gamma^0(\delta/2)$. So for sufficiently small $\epsilon$ the property $\left|d^K_\epsilon\right| = \left| \sum^K_{i=0}{\epsilon^i d^i} \right| \geq \delta/4$ is valid for all $(x,t) \in \Gamma^0(\delta) \backslash \Gamma^0(\delta/2)$. Applying the matching conditions (\ref{match1}) yields
\begin{eqnarray}
\lefteqn{ \left\| c^K_A - c^K_O \right\|_{C^2(\Gamma^0(\delta) \backslash \Gamma^0(\delta/2))}} \nonumber \\
& = & \left\| \zeta(d^0/\delta) \sum^K_{i=0}{\epsilon^i \left(c^i (d^K_\epsilon/\epsilon,x,t) - c^\pm_i(x,t)\right)} \right\|_{C^2(\Gamma^0(\delta) \backslash \Gamma^0(\delta/2))} \nonumber \\
& = & \mathcal{O}(\epsilon^{-2} e^{-\frac{\alpha \delta}{4 \epsilon}}) \, ,
\end{eqnarray}
and analogously we get
\begin{eqnarray}
\left\| \mu^K_A - \mu^K_O \right\|_{C^2(\Gamma^0(\delta) \backslash \Gamma^0(\delta/2))} & = & \mathcal{O}(\epsilon^{-2} e^{-\frac{\alpha \delta}{4 \epsilon}}), \label{fehleruebergang1} \\
\left\| \textbf{u}^K_A - \textbf{u}^K_O \right\|_{C^2(\Gamma^0(\delta) \backslash \Gamma^0(\delta/2))} & = & \mathcal{O}(\epsilon^{-2} e^{-\frac{\alpha \delta}{4 \epsilon}}) .\label{fehleruebergang}
\end{eqnarray}
By using the outer-boundary matching condition (\ref{matchbound1}), a similar statement holds 
\begin{eqnarray}
\lefteqn{ \left\| c^K_A - c^K_O \right\|_{C^2(\partial_T \Omega(\delta) \backslash \partial_T \Omega(\delta/2))}} \nonumber \\
& = & \left\| \zeta(d_B/\delta) \left( \sum^K_{i=1}{\epsilon^i \left(c^i_B(d_B/\epsilon) - c^+_i\right)} - \epsilon^K c^K_B(0) \right) \right\|_{C^2(\partial_T \Omega(\delta) \backslash \partial_T \Omega(\delta/2))} \nonumber \\
& = & \mathcal{O}(\epsilon^{-2} e^{- \frac{\alpha \delta}{2 \epsilon}}) + \mathcal{O}(\epsilon^K) \, . \label{fehlerrand}
\end{eqnarray}
A similar relation holds for $(\mu^K_A, \textbf{u}^K_A)$ and $(\mu^K_O, \textbf{u}^K_O)$. Since the equations (\ref{system1})-(\ref{system3}) contain second space derivatives and a first time derivative, we have used the $C^2$-norm. Therefore by (\ref{nabladK})-(\ref{fehlerrand}) the approximate solution $(c^K_A,\mu^K_A,\textbf{u}^K_A)$ fulfills the following equations
\begin{align*}
(c^K_A)_t - \Delta \mu^K_A = : e_K(x,t) & = \mathcal{O}(\epsilon^{K-2}) && \mbox{in } \Omega\times (0,T)\, , \\
\mu^K_A + \epsilon \Delta c^K_A - \epsilon^{-1} f(c^K_A) - W_{,c}(c^K_A,\mathcal{E}(\textbf{u}^K_A)) & = \mathcal{O}(\epsilon^{K-1}) && \mbox{in } \Omega \times (0,T) \, ,\\
\Div \left( \mathcal{C} \mathcal{E}(\textbf{u}^K_A) \right) - \Div \left(\mathcal{C} \mathcal{E}^\star c^K_A\right) & = \mathcal{O}(\epsilon^{K-2}) && \mbox{in } \Omega \times (0,T) \, , \\
\tfrac{\partial}{\partial n} c^K_A = \tfrac{\partial}{\partial n} \mu^K_A & = 0 && \mbox{on } \partial \Omega \times (0,T) \, , \\
\textbf{u}^K_A & = 0 && \mbox{on } \partial \Omega \times (0,T) \, .
\end{align*}
Observe that here and in the following the Landau symbols are in $C^0$-norm unless noted otherwise. In the same way as in \cite{ABC}, we modify $c^K_A$ and $\mu^K_A$ so that the error term $e_K$ vanishes. We set $c^\epsilon_A = c^K_A - \frac{1}{\left|\Omega\right|} \int^t_0{\int_\Omega{e_K(\xi,\tau) \, d\xi} \, d\tau}$ and $\mu^\epsilon_A = \mu^K_\epsilon - \hat{e}_K(x,t)$ where $\hat{e}_K(x,t)$ is the solution to the elliptic problem 
\begin{eqnarray*}
\Delta \hat{e}_K(x,t) = e_K(x,t) - \frac{1}{\left|\Omega\right|} \int_\Omega{e_K(\xi,t) \, d\xi} \quad \mbox{in }\Omega_T \,, \\
\frac{\partial}{\partial n} \hat{e}_K = 0 \quad \mbox{on } \partial_T \Omega \, , \quad \int_\Omega{\hat{e}_K(\xi,t) \, d\xi} = 0 \quad \forall t \in [0,T]\,.
\end{eqnarray*}
Note that $ \hat{e}_K  = \mathcal{O}(\epsilon^{K-2})$ since $  e_K   = \mathcal{O}(\epsilon^{K-2})$. In addition, we define $\textbf{u}^\epsilon_A = \textbf{u}^K_A$. Therefore $(c^\epsilon_A, \mu^\epsilon_A, \textbf{u}^\epsilon_A)$ satisfies

\begin{remark}
With (\ref{fehlerabsch}) and (\ref{k}) we can specify the size of $K$. 
\[
K-3 \geq \frac{pk}{2} > d +2,5\,.
\]
In particular, it is sufficient to calculate the 8th order term of the expansion in two dimensions and the 9th order term in three dimensions.
\end{remark}

We summarize the results of this subsection in the following theorem.\\

\noindent\textbf{Proof of Theorem \ref{theoremcA-c0}} The construction of an approximate solution $(c^\epsilon_A,\mu^\epsilon_A,\textbf{u}^\epsilon_A)$ satisfying (\ref{apprsol1})-(\ref{apprsol5}) is described above. \\
Due to the construction of $\mu^\epsilon_A$ and (\ref{fehleruebergang1}), it follows as $\epsilon \searrow 0$
\begin{eqnarray*}
\left\| \mu_A^\epsilon - \mu_0 \right\|_{C^0(\Omega_T \backslash \left(\Gamma^0(\delta/2) \cup \partial_T \Omega(\delta/2) \right) )} & = & \mathcal{O}(\epsilon) \, .
\end{eqnarray*}
So it remains to consider the domains $\Gamma^0(\delta/2)$ and $\partial_T \Omega(\delta/2)$. By triangle inequality it holds
\begin{eqnarray*}
\left\| \mu_A^\epsilon - \mu_0 \right\|_{C^0( \Gamma^0(\delta/2))} \leq \left\| \mu_A^\epsilon - \mu^0 \right\|_{C^0( \Gamma^0(\delta/2))} + \left\| \mu^0 - \mu_0 \right\|_{C^0( \Gamma^0(\delta/2))} \, .
\end{eqnarray*}
By definition it follows $ \left\| \mu_A^\epsilon - \mu^0 \right\|_{C^0( \Gamma^0(\delta/2))} \leq C \epsilon $. To estimate the second term, Lemma \ref{lemma0teOrd} gives us the exact definition of $\mu^0$. Therefore we have
\begin{eqnarray*}
\left\| \mu^0 - \mu_0 \right\|_{C^0( \Gamma^0(\delta/2) \cap Q^+_0)} & = & \left\| (1-\eta(d^K_\epsilon/\epsilon)) (\mu^+_0-\mu^-_0) \right\|_{C^0( \Gamma^0(\delta/2) \cap Q^+_0)} \\
& = & \left\| \chi_{\left\{d^K_\epsilon \leq \epsilon\right\}} (\mu^+_0-\mu^-_0) \right\|_{C^0( \Gamma^0(\delta/2) \cap Q^+_0)} \\
& \leq & C \epsilon \, ,
\end{eqnarray*}
where the second equality follows from (\ref{eta}) and the last inequality from $\mu^+_0 = \mu^-_0$ on $\Gamma^0$ and $\left\{ d^K_\epsilon \leq \epsilon \right\} \subset \left\{ d^0 \leq C \epsilon \right\}$ for some $C >0$. The proof for $\Gamma^0(\delta/2) \cap Q^-_0$ is done in the same way. Since $\mu^0_B = \mu^+_0 $ in $\partial_T \Omega(\delta/2)$, the construction of $\mu^\epsilon_A$ yields
\[
\left\| \mu^\epsilon_A - \mu_0 \right\|_{C^0( \partial_T \Omega(\delta/2))} \leq C \epsilon 
\]
for some $C>0$ independent of $\epsilon$. We show analogously 
\[
\left\| \textbf{u}_A^\epsilon - \textbf{u}_0 \right\|_{C^0(\Omega_T)} + \left\| c^\epsilon_A \mp 1 \right\|_{C^0(Q^\pm_0 \backslash \Gamma^0(\delta/2))} = \mathcal{O}(\epsilon) \, .
\]
To estimate the last term, we consider again the domains $\Gamma^0(\delta/2)$ and $\Gamma^0(\delta) \backslash \Gamma^0(\delta/2)$ separately. We use that $c^\epsilon_A = c^\epsilon_A - c^K_A + c^K_A$ and apply $c^K_A = c^I_K $ in $\Gamma^0(\delta/2)$ and the triangle inequality to obtain 
\begin{eqnarray}
\lefteqn{ \left\| c^\epsilon_A - \theta_0(d^0/\epsilon + d^1) \right\|_{C^0(\Gamma^0(\delta/2))}} \nonumber\\
& \leq & \left\| c^\epsilon_A - c^K_A \right\|_{C^0(\Gamma^0(\delta/2))} + \epsilon \left\| \sum^K_{i=1}{\epsilon^{i-1} c^i} \right\|_{C^0(\Gamma^0(\delta/2))} \nonumber \\
&& + \left\| \theta_0 \! \left(\frac{d^0}{\epsilon} + d^1 + \epsilon \sum^K_{i=2}{\epsilon^{i-2} d^i}\right) - \theta_0 \! \left(\frac{d^0}{\epsilon} + d^1 \right) \right\|_{C^0(\Gamma^0(\delta/2))} \nonumber \\
& \leq & C \epsilon^{K-2} +C \epsilon + C \epsilon  \label{Liptheta}
\end{eqnarray}
since $\theta_0$ is a Lipschitz function. In $\Gamma^0(\delta) \backslash \Gamma^0(\delta/2)$ we write the difference $c^\epsilon_A - \theta_0(d^0/\epsilon + d^1)$ as
\begin{eqnarray*}
\lefteqn{c^\epsilon_A - \theta_0(d^0/\epsilon + d^1)}\\
& = & \zeta(d^0/\delta) \left(\theta_0(d^K_\epsilon/\epsilon) - \theta_0(d^0/\epsilon + d^1)\right) \\
&& + \left( c^\epsilon_A - \left[\zeta(d^0/\delta) \theta_0(d^K_\epsilon/\epsilon) + (1- \zeta(d^0/\delta))(2 \chi_{Q^+_0}-1) \right]\right)\\ 
&& + (1- \zeta(d^0/\delta)) \left((2 \chi_{Q^+_0}-1) - \theta_0(d^0/\epsilon + d^1) \right) .
\end{eqnarray*}
On the right-hand side the first term can be estimated as in (\ref{Liptheta}). For the second term we use the definition for $c^K_A$ in $\Gamma^0(\delta) \backslash \Gamma^0(\delta/2)$
\begin{eqnarray*}
\lefteqn{\left\|c^\epsilon_A - \left[\zeta(d^0/\delta) \theta_0(d^K_\epsilon/\epsilon) + (1- \zeta(d^0/\delta))(2 \chi_{Q^+_0}-1) \right]\right\|_{C^0(\Gamma^0(\delta) \backslash \Gamma^0(\delta/2))}} \\
& \leq & \left\| c^\epsilon_A - c^K_A \right\|_{C^0(\Gamma^0(\delta) \backslash \Gamma^0(\delta/2))} \\
&&+ \left\| c^K_A - \left[\zeta(d^0/\delta) \theta_0(d^K_\epsilon/\epsilon) + (1- \zeta(d^0/\delta))(2 \chi_{Q^+_0}-1) \right] \right\|_{C^0(\Gamma^0(\delta) \backslash \Gamma^0(\delta/2))} \\
& \leq & C \epsilon^{K-2} + C \epsilon \,.
\end{eqnarray*}
To estimate the third term on the right-hand side, we use that $\left| d^0 + \epsilon d^1 \right| \geq \delta / 4$ in $\Gamma^0(\delta) \backslash \Gamma^0(\delta / 2)$ for all $\epsilon$ small enough. Then applying the property (\ref{proptheta2}) yield
\[
\left\| (2 \chi_{Q^+_0}-1) - \theta_0(d^0/\epsilon + d^1)  \right\|_{C^0(\Gamma^0(\delta) \backslash \Gamma^0(\delta/2))} \leq C e^{-\frac{\alpha \delta}{8 \epsilon}} 
\] 
for some constant $C>0$. This completes the proof.
\makebox[1cm]{} \hfill $\Box$\\

\section{Main Result}\label{secmain}

\begin{theorem} \label{mainresult}
Let $\Omega$ be a smooth domain and $\Gamma_{00}$ be a smooth hypersurface in $\Omega$ without boundary. Assume that the Hele-Shaw problem (\ref{sharpsystem1})-(\ref{sharpsystem7}) has a smooth solution $(\mu, \textbf{u},\Gamma)$ on a time interval $[0,T]$ such that $\Gamma_t \subset \Omega$ for all $t \in [0,T]$ where $\Gamma: = \bigcup_{0 \leq t \leq T} (\Gamma_t \times \left\{ t \right\})$. Then there exists a family of smooth functions $\left\{ c^\epsilon_0(x) \right\}_{0 < \epsilon \leq 1}$ which are uniformly bounded in $\epsilon \in (0,1]$ and $x \in \overline{\Omega}$, such that if $(c^\epsilon, \textbf{u}^\epsilon)$ satisfies the Cahn-Larché equation 
\begin{align}
c^\epsilon_t - \Delta \left( - \epsilon \Delta c^\epsilon + \epsilon^{-1} f(c^\epsilon) +  W_{,c}(c^\epsilon, \mathcal{E}(\textbf{u}^\epsilon)) \right) & = 0 & \mbox{in }& \Omega_T, \\
\mathrm{div} \, W_{, \mathcal{E}}(c^\epsilon, \mathcal{E}(\textbf{u}^\epsilon)) & = 0 &\mbox{in }&\Omega_T, \\
\tfrac{\partial}{\partial n}c^\epsilon = \tfrac{\partial}{\partial n} \left( - \epsilon \Delta c^\epsilon + \epsilon^{-1} f(c^\epsilon) + W_{,c}(c^\epsilon, \mathcal{E}(\textbf{u}^\epsilon)) \right) & = 0 &\mbox{on }& \partial_T \Omega,\\
\textbf{u}^\epsilon & = 0 &\mbox{on }& \partial_T \Omega,\\
\left. c^\epsilon \right|_{t=0}  & = c^\epsilon_0 &\mbox{on }& \Omega  \,, 
\end{align}
then
\begin{align*}
& \lim_{\epsilon \to 0} c^\epsilon(x,t) = \left\{\begin{matrix}
-1 & \mbox{if } (x,t) \in Q^-\\
1 & \mbox{if } (x,t) \in Q^+ 
 \end{matrix}\right.  \mbox{ uniformly on compact subsets,} &\\
& \lim_{\epsilon \to 0}\left( - \epsilon \Delta c^\epsilon + \epsilon^{-1} f(c^\epsilon) +  W_{,c}(c^\epsilon, \mathcal{E}(\textbf{u}^\epsilon)) \right)(x,t) = \mu(x,t)  \mbox{ uniformly on } \overline{\Omega_T} \,,&\\
& \lim_{\epsilon \to 0} \textbf{u}^\epsilon(x,t) = \textbf{u}(x,t)  \mbox{ uniformly on } \overline{\Omega_T} \,,&
\end{align*}
where $Q^+ $ and $Q^-$ are respectively the exterior (in $\Omega_T$) and interior of $\Gamma$. 
\end{theorem}

\proof An existence and uniqueness result for the Cahn-Larché equation can be found in \cite{GarckeCH}. Let $(c^\epsilon_A, \mu^\epsilon_A, \textbf{u}^\epsilon_A)$ be the approximate solution constructed in Theorem \ref{theoremcA-c0}. Then by Theorem \ref{hoeregR} and \ref{theoremcA-c0}, we obtain
\begin{eqnarray*}
\lim_{\epsilon \to 0}\left\|c^\epsilon \mp 1\right\|_{C^0(Q^\pm \backslash \Gamma(\delta/2))} & = & 0 \,, \\
\lim_{\epsilon \to 0}\left\|\mu^\epsilon - \mu \right\|_{C^0(\Omega_T)} & = & 0 \,, \\
\lim_{\epsilon \to 0} \left\|\textbf{u}^\epsilon - \textbf{u} \right\|_{C^0(\Omega_T)} & = & 0 
\end{eqnarray*}
for every $\delta >0$ small enough, as long as $\Phi^\epsilon_t(.) = c^\epsilon_A(.,t)$ has the form (\ref{bedingungphi2}).\\
Hence we have to check that $\Phi^\epsilon_t(.) = c^\epsilon_A(.,t)$ has the form (\ref{bedingungphi2}) where $r = r_t(x)$ is the signed distance function to $\Gamma^{\epsilon K }_t := \left\{ x \in \Omega : d^K_\epsilon(x,t) = 0 \right\}$. This can be shown in the same way as \cite[Theorem5.1]{ABC}. In particular, we can prove that $c^1(z,x,t) = \Delta d^0(x,t) \theta_1(z)$ for all $(z,x,t) \in \mathbb{R} \times \Gamma^0$, where $\theta_1$ satisfies (\ref{theta1}). This can be seen as follows. To this end we solve for $c^1$. For $(x,t) \in \Gamma^0$ the equation for $c^1$ in (\ref{innere1}) reads
\begin{eqnarray}
c^1_{zz} - f'(\theta_0) c^1 & = & - \mathcal{E}^\star : \mathcal{C} ((\textbf{u}^0_\ast)_z \otimes \nabla d^0) - \mu^0 - \Delta d^0 c^0_z \nonumber \\
&& - \mathcal{E}^\star : \mathcal{C} (\nabla \textbf{u}^0 - \mathcal{E}^\star \theta_0)- z k^0 \eta' \nonumber \\
& = & - \mathcal{E}^\star : \mathcal{C} ((\textbf{u}_\ast^0)_z \otimes \nabla d^0) + \sigma \Delta d^0 + \tfrac{1}{2} \mathcal{E}^\star : \mathcal{C} (\nabla \textbf{u}^+_0 + \nabla \textbf{u}^-_0) \nonumber \\
&& - \Delta d^0 \theta'_0 - \mathcal{E}^\star : \mathcal{C} (\nabla \textbf{u}^+_0 \eta + \nabla \textbf{u}^-_0(1- \eta ) - \mathcal{E}^\star \theta_0) \nonumber \\
&& - z \mathcal{E}^\star : \mathcal{C}(\textbf{l}^0 \otimes \nabla d^0) \eta' \nonumber \\
& = & - \mathcal{E}^\star : \mathcal{C} ((\textbf{u}_\ast^0)_z \otimes \nabla d^0) + \sigma \Delta d^0 - \Delta d^0 \theta'_0 \nonumber \\
&& + \tfrac{1}{2} \mathcal{E}^\star : \mathcal{C} (\nabla \textbf{u}^+_0 - \nabla \textbf{u}^-_0) - \mathcal{E}^\star : \mathcal{C} (\nabla \textbf{u}^+_0 - \nabla \textbf{u}^-_0) \eta  \nonumber \\
&& + \mathcal{E}^\star : (\mathcal{C} \mathcal{E}^\star) \theta_0 - z \mathcal{E}^\star : \mathcal{C}(\textbf{l}^0 \otimes \nabla d^0) \eta' \, , \label{c1onGamma0}
\end{eqnarray}
where we have used the definitions of $(c^0,\mu^0,\textbf{u}^0,k^0)$. To handle the term on the right-hand side, we calculate $(\textbf{u}^0_\ast)_z$. For that we use the ordinary differential equation (\ref{innere3}) on $\Gamma^0$. By (\ref{uast}) we obtain $\lim_{z \to -\infty} \partial_z \textbf{u}^0_\ast(z) = 0 $, and therefore we get
\begin{eqnarray*}
(\textbf{u}^0_\ast)_z & = & \int^z_{-\infty}{(\textbf{u}^0_\ast)_{zz}(s) \, ds} \\
& = & \int^z_{-\infty}M^{-1} \left[ - (\mathcal{C}_{iji'j'} \partial_j(\textbf{u}^0_{i'})_z \partial_{j'} d^0)_{i=1, \ldots ,d} - (\mathcal{C}_{iji'j'} (\textbf{u}^0_{i'})_z \partial_{jj'} d^0)_{i=1, \ldots ,d} \right.\\
&& \left. - (\mathcal{C} \nabla \textbf{u}^0_z) \nabla d^0 + \theta'_0 (\mathcal{C} \mathcal{E}^\star) \nabla d^0  \right] - s \textbf{l}^0 \eta'' \, ds \,.
\end{eqnarray*}
The matching condition $\lim_{z \to -\infty} (\theta_0(z), \textbf{u}^0(z)) = (-1,\textbf{u}^-_0)$ (see (\ref{match1})) and the equation $\int^z_{-\infty}{s \eta''(s) \, ds} = z \eta'(z) - \eta(z)$ yield
\begin{eqnarray*}
(\textbf{u}^0_\ast)_z & = & M^{-1} \left[ - (\mathcal{C}_{iji'j'} \partial_j(\textbf{u}^+_0 \eta + \textbf{u}^-_0(1-\eta))_{i'} \partial_{j'} d^0)_{i=1, \ldots ,d} \right. \\
&& + (\mathcal{C}_{iji'j'} \partial_j\textbf{u}^-_{0,i'} \partial_{j'} d^0)_{i=1, \ldots ,d} \\
&& - (\mathcal{C}_{iji'j'} (\textbf{u}^+_0 \eta + \textbf{u}^-_0(1-\eta))_{i'} \partial_{jj'} d^0)_{i=1, \ldots ,d} \\
&&  + (\mathcal{C}_{iji'j'} \textbf{u}^-_{0,i'} \partial_{jj'} d^0)_{i=1, \ldots ,d} - (\mathcal{C} \nabla (\textbf{u}^+_0 \eta + \textbf{u}^-_0(1 - \eta))) \nabla d^0 \\
&& \left.+ (\mathcal{C} \nabla \textbf{u}^-_0) \nabla d^0 + \theta_0(\mathcal{C} \mathcal{E}^\star) \nabla d^0 + (\mathcal{C} \mathcal{E}^\star) \nabla d^0 \right] - z \eta' \textbf{l}^0 + \eta \textbf{l}^0 \,,
\end{eqnarray*}
where we have used the definition of $\textbf{u}^0$. Since $\left[\mathcal{S} \nu\right]_{\Gamma^0_t} = \left[\textbf{u}^\pm_0\right]_{\Gamma^0_t} = 0 $, we can conclude
\begin{eqnarray*}
(\textbf{u}^0_\ast)_z & = & M^{-1} \left[ - (\mathcal{C}_{iji'j'} (\nabla \textbf{u}^+_0 - \nabla \textbf{u}^-_0)_{i'j} \partial_{j'} d^0)_{i=1, \ldots ,d}\, \eta \right. \\
&&  - (\mathcal{C} (\nabla \textbf{u}^+_0 - \nabla \textbf{u}^-_0)) \nabla d^0 \eta + \tfrac{1}{2} \theta_0 (\mathcal{C}(\nabla \textbf{u}^+_0 - \nabla \textbf{u}^-_0)) \nabla d^0 \\
&& \left. + \tfrac{1}{2} (\mathcal{C}(\nabla \textbf{u}^+_0 - \nabla \textbf{u}^-_0)) \nabla d^0 \right] - z \eta' \textbf{l}^0 + \eta \textbf{l}^0 \quad \mbox{on }\Gamma^0\,.
\end{eqnarray*}
The interface condition $\left[\textbf{u}^\pm_0 \right]_{\Gamma^0_t} = 0$ also yields $\nabla \textbf{u}^+_0 - \nabla \textbf{u}^-_0 = \frac{\partial}{\partial \nu} (\textbf{u}^+_0 - \textbf{u}^-_0) \otimes \nabla d^0$ on $\Gamma^0$. Furthermore, note that $M = (\mathcal{C}_{iji'j'} \partial_j d^0 \partial_{j'} d^0)^d_{i,i'= 1}$. So we obtain
\begin{eqnarray*}
(\textbf{u}^0_\ast)_z(z) & = & - 2 \eta(z) \frac{\partial}{\partial \nu} (\textbf{u}^+_0 - \textbf{u}^-_0) + \frac{1}{2} \theta_0(z) \frac{\partial}{\partial \nu} (\textbf{u}^+_0 - \textbf{u}^-_0) \\
&& + \frac{1}{2} \frac{\partial}{\partial \nu} (\textbf{u}^+_0 - \textbf{u}^-_0) - z \eta'(z) \textbf{l}^0 + \eta \textbf{l}^0 \\
& = & - \eta(z) \frac{\partial}{\partial \nu} (\textbf{u}^+_0 - \textbf{u}^-_0) + \frac{1}{2} \theta_0(z) \frac{\partial}{\partial \nu} (\textbf{u}^+_0 - \textbf{u}^-_0) \\
&& + \frac{1}{2} \frac{\partial}{\partial \nu} (\textbf{u}^+_0 - \textbf{u}^-_0) - z \eta'(z) \textbf{l}^0 \quad \mbox{on }\Gamma^0\,,
\end{eqnarray*}
where we have used the definition of $\textbf{l}^0$ in the last equation. We insert this equation into (\ref{c1onGamma0}) and use $\nabla \textbf{u}^+_0 - \nabla \textbf{u}^-_0 = \frac{\partial}{\partial \nu} (\textbf{u}^+_0 - \textbf{u}^-_0) \otimes \nabla d^0$ on $\Gamma^0$ to obtain
\begin{eqnarray*}
c^1_{zz} - f'(\theta_0) c^1 = \sigma \Delta d^0 - \theta'_0 \Delta d^0 + \theta_0 \left( \mathcal{E}^\star : (\mathcal{C} \mathcal{E}^\star) - \tfrac{1}{2} \mathcal{E}^\star : \mathcal{C} (\nabla \textbf{u}^+_0 - \nabla \textbf{u}^-_0 ) \right) .
\end{eqnarray*}
Since $\left[\mathcal{S} \nu\right]_{\Gamma^0_t} = 0$ and $\Div \mathcal{S} = 0$ in $Q^\pm_0$, we obtain $\mathcal{C} \mathcal{E}^\star - \tfrac{1}{2} \mathcal{C} (\nabla \textbf{u}^+_0 - \nabla \textbf{u}^-_0 ) = \left[\mathcal{S}\right]_{\Gamma^0_t} = 0 $. Therefore $c^1$ satisfies the following ordinary differential equation on $\Gamma^0$
\[
c^1_{zz} - f'(\theta_0) c^1 = \sigma \Delta d^0 - \theta'_0 \Delta d^0 \,,
\]
and so $c^1(z,x,t) = \Delta d^0(x,t) \theta_1(z)$ on $\Gamma^0$ where $\theta_1$ satisfies
\[
\theta''_1 - f'(\theta_0) \theta_1  = \sigma - \theta'_0 \mbox{ in } \mathbb{R}\,, \quad \theta_1(0) = 0 \,, \quad \theta_1 \in L^\infty(\mathbb{R})\,.
\]
Since $\int_\mathbb{R}{(\sigma - \theta'_0) \theta'_0} = 0$ by definition of $\sigma$, \cite[Lemma 4.1]{ABC} yields that $\theta_1$ exists and is unique. By integration by parts we show that $\theta_1$ satisfies the property (\ref{theta1}) 
\begin{align*}
0 & = \int_\mathbb{R}{\theta''_0 (\sigma - \theta'_0) \, dz} = \int_{\mathbb{R}}{\theta''_0 \left(\theta''_1 - f'(\theta_0) \theta_1 \right)\, dz} \nonumber  \\
& = - \int_\mathbb{R}{\theta_1'\left( \theta''_0 - f(\theta_0) \right)'dz} + \int_\mathbb{R}{f''(\theta_0) (\theta'_0)^2 \theta_1 \, dz} = \int_\mathbb{R}{f''(\theta_0) (\theta'_0)^2 \theta_1 \, dz} \,. \label{eigentheta1}
\end{align*}
As mentioned above the rest of the proof is done in the same way as \cite[Theorem 5.1]{ABC}.
\makebox[1cm]{} \hfill $\Box$\\

\begin{remark}
The initial value $c^\epsilon_0(x)=c^\epsilon_A(x,0)$ is independent of the solutions for the modified Hele-Shaw problem and the linearized Hele-Shaw problem. This can be seen as follows. By solving the first order partial differential equation (\ref{djbedingung}) with Cauchy data $d^k(x,0)= 0$ on $\Gamma_{00}$, we can directly determine $d^k$ for all $k \in \mathbb{N}$. Hence one can find $\mathcal{V}^k$ for $t=0$ and for all $k \in \mathbb{N}$.
\end{remark}

\medskip

\noindent
{\bf Acknowledgments:} The authors acknowledge supported by the SPP 1506 "Transport Processes
at Fluidic Interfaces" of the German Science Foundation (DFG) through grant AB285/3-1, AB285/4-1 and AB285/4-2.

%\addcontentsline{toc}{chapter}{Bibliography}

%\input{biblio}

%\printbibliography

%\bibliography{literaturverzeichnis}
%\bibliographystyle{plain}

%\addcontentsline{toc}{chapter}{Bibliography}

\end{document}